\documentclass[publmath,final]{svjour}  
\usepackage{publmath} 
\usepackage{graphicx}          
\usepackage{epsfig}
\usepackage[all]{xy}
\newcommand{\Cone}{{\rm Cone}}
\newcommand{\Ang}{{\rm Ang}}
\newcommand{\MaxAng}{{\rm MaxAng}}
\spnewtheorem{thm}{Theorem}[section]{\it}{\rm}
\spnewtheorem{defn}[thm]{Definition}{\it}{\rm}
\spnewtheorem{lem}[thm]{Lemma}{\it}{\rm}
\spnewtheorem{cor}[thm]{Corollary}{\it}{\rm}
\spnewtheorem{prop}[thm]{Proposition}{\it}{\rm}
\spnewtheorem{ex}[thm]{Example}{\it}{\rm}

\spnewtheorem{remarknum}[thm]{Remark}{\it}{\rm}
\spnewtheorem{hypremark}[thm]{Remark (Hyperbolic case)}{\it}{\rm}
\spnewtheorem*{notation}{Notation}{\it}{\rm}
\spnewtheorem*{convention}{Convention}{\it}{\rm}
\spnewtheorem{examples}[thm]{Examples}{\it}{\rm}
\spnewtheorem{thmA}{Theorem}{\it}{\it}
\spnewtheorem{corA}[thmA]{Corollary}{\it}{\it}
\def\R {\mathbb R}

\begin{document}

%%%%%%%%%%%%%%%%%%%%%%%%%%%%%%
%%% FRONT PAGE FOR PUBLISHER
%%%%%%%%%%%%%%%%%%%%%%%%%%%%%%

\title{THE ISOMORPHISM PROBLEM FOR TORAL RELATIVELY HYPERBOLIC GROUPS}
\titlerunning{ISOMORPHISM PROBLEM FOR RELATIVELY HYPERBOLIC GROUPS}

\author{\firstname{Fran\c{c}ois} DAHMANI \and 
\firstname{Daniel} GROVES\thanks{The second author was supported in part by NSF
Grant DMS-0504251.  This work was undertaken whilst the second author
was a Taussky-Todd Instructor at the California Institute of Technology.}}

\institute{F. D. \\
Institut de Math\'ematiques de Toulouse, \\
Universit\'e Paul Sabatier, Toulouse III, \\
31062, Toulouse, cedex 9, France\\
francois.dahmani@math.univ-toulouse.fr
\and
D. G. \\
MSCS UIC 322 SEO, \textsc{M/C} 249\\
851 S. Morgan St.\\
Chicago, IL 60607-7045, USA\\
groves@math.uic.edu
}

\date{??/??/200?}

%%%%%%%%%%%%%%%%%%%%%%%%%
%%%  END/ FRONT PAGE PUBLISHER
%%%%%%%%%%%%%%%%%%%%%%%%%

%%%%%%%%%%%%%%%%%%%%%%%%%
%% FRONT PAGE FOR ARXIV
%%%%%%%%%%%%%%%%%%%%%%%%%

%\title[ISOMORPHISM PROBLEM FOR RELATIVELY HYPERBOLIC GROUPS]{THE ISOMORPHISM PROBLEM FOR TORAL RELATIVELY HYPERBOLIC GROUPS}
%%\author[Fran\c{c}ois Dahmani]{Fran\c{c}ois Dahmani}
%\address{Fran\c{c}ois Dahmani\\
%Institut de Math\'ematiques de Toulouse, \\
%Universit\'e Paul Sabatier\\
%F-31062 Toulouse, cedex 9, France}
%\email{dahmani@math.univ-toulouse.fr}

%\author[Daniel Groves]{Daniel Groves}
%\address{Daniel Groves\\
%MSCS UIC 322 SEO, \textsc{M/C} 249\\
%851 S. Morgan St.\\
%Chicago, IL 60607-7045, USA}
%\email{groves@math.uic.edu}

%\date{November 2007}
%\subjclass[2000]{20F10, (20F65)}

%\thanks{The second author was supported in part by NSF
%Grant DMS-0504251.  This work was undertaken whilst the second author
%was a Taussky-Todd Instructor at the California Institute of Technology.}
%
%
%
%

%%%%%%%%%%%%%%%%%%%%%%%%%
%% END FRONT PAGE ARXIV
%%%%%%%%%%%%%%%%%%%%%%%%%

\maketitle
\begin{abstract}
We provide a solution to the isomorphism problem for torsion-free relatively hyperbolic groups with abelian parabolics.  As special cases we recover solutions to the isomorphism problem for: 
(i) torsion-free hyperbolic groups (Sela, \cite{SelaIso} and unpublished); and (ii) finitely generated fully residually free groups (Bumagin, Kharlampovich and Miasnikov \cite{BKM}). 
We also give a solution to the homeomorphism problem for finite volume hyperbolic $n$-manifolds, for $n \ge 3$. 
In the course of the proof of the main result, we prove that a particular JSJ decomposition of a freely indecomposable torsion-free relatively hyperbolic group with abelian parabolics is algorithmically constructible.
\end{abstract}

\setcounter{tocdepth}{1}
\tableofcontents

 \section{Introduction} 
 
The {\em isomorphism problem} is the hardest of the three fundamental problems 
described by Dehn in 1912 (the other two are the {\em word problem} and the 
{\em conjugacy problem}: see \cite{Dehn1} and \cite{Dehn2}).  
The isomorphism problem asks for a general algorithm which will, given two finite group presentations, 
decide whether or not the presentations define isomorphic groups.  
For finitely presented groups in general, Adian \cite{Adian} and Rabin \cite{Rabin} 
proved that there is no such algorithm.  One can then ask whether there is a solution 
within a class $\mathcal C$ of groups.  Namely, is there an algorithm which, given 
two finite group presentations and the knowledge that they define groups in $\mathcal C$, 
decides whether or not the presentations define isomorphic groups?
Ideally, a positive solution to this question should not require that
the presentations are given along with a proof that the groups defined
lie in $\mathcal C$, merely the knowledge that they do should suffice.
This is the approach we take in this paper (see Remark \ref{Rem:FindRH}).

For interesting classes of groups, the expected answer to the above
question is `no'.  In fact,
there are very few large classes of groups for which the isomorphism problem is known to be solvable.  For finite groups it is easy to see that there is such an algorithm, for given enough time one can write down the multiplication table from the group presentation, and then check if the two tables define the same group.  For abelian groups, the algorithm is similarly straightforward.  Other classes of groups for which there is a positive answer include nilpotent groups (Grunewald and Segal \cite{GrunSeg}) and  polycyclic-by-finite groups (Segal \cite{Segal}). 

In recent years, geometric group theory has provided solutions to the isomorphism problem for a few more classes of groups.  We mention two classes.  Sela \cite{SelaIso} solved the isomorphism problem for torsion-free hyperbolic groups which do not admit a small essential action on an $\R$-tree, and has an (unfortunately unpublished) proof for arbitrary torsion-free hyperbolic groups.  The second class we mention is finitely generated fully residually free groups (also known as `limit groups'), for which a solution to the isomorphism problem was provided by Bumagin, Kharlampovich and Miasnikov \cite{BKM}.

We mention two negative results about the isomorphism problem.  Baumslag, Gildenhuys and Strebel \cite{BGS} proved that the isomorphism problem is unsolvable for solvable groups of derived length $3$.  In \cite{Miller}, Miller studied many algorithmic questions.  By starting with a group $G$ with unsolvable word problem, one can construct a series of groups $G_w$ (indexed by words in the generating
set of $G$) so that $G_w$ has solvable word problem if and only
if $w =_G 1$ (see \cite[Chapter V]{Miller}).  By then applying the 
construction from \cite[Chapter III]{Miller}, we obtain a free-by-free
group $H_w$ which has solvable conjugacy problem if and only if
$w =_G 1$.  Thus one cannot decide if $H_w \cong H_1$, and the
isomorphism problem is unsolvable for free-by-free groups.  It is
not difficult to see that the groups $H_w$ all have a cubic isoperimetric
function.  The isomorphism problem for the class of groups which satisfy a quadratic isoperimetric inequality remains open.
See \cite{Miller:dec} for a discussion of decision problems
in group theory, and many more examples.

For definitions of relatively hyperbolic groups, see Section \ref{Prelim} below.  A relatively hyperbolic group is called {\em toral} if it is torsion-free and its parabolic subgroups are all (finitely generated) abelian.  The main result of this paper is the following

\begin{thmA} \label{Toral}
The isomorphism problem is solvable for the class of toral relatively hyperbolic groups.
\end{thmA}

As special cases of Theorem \ref{Toral} we recover the above-mentioned results of Sela and of Bumagin, Kharlampovich and Miasnikov.

\begin{corA} [Sela, \cite{SelaIso} and unpublished] \label{hyp}
The isomorphism problem is solvable for the class of torsion-free hyperbolic groups.
\end{corA}
Sela has a proof of Corollary \ref{hyp}, but it remains unpublished.
The published case in \cite{SelaIso} is a particular special (but key)
case.

It is worth remarking that in the case of torsion-free hyperbolic 
groups which admit no essential small action on an $\R$-tree (this is the
case from \cite{SelaIso}), our methods provide significant
simplifications.  The major innovation in our approach is the use of equations
with {\em rational constraints} (see Sections \ref{SS:Algo} and \ref{List} below for
definitions and more discussion on this).  This greatly streamlines the solution to 
the isomorphism problem (this point is developed in Section \ref{Strategy},
before Remark \ref{rem;expo}).

Another immediate corollary of Theorem \ref{Toral} is the following
result.

\begin{corA} [Bumagin, Kharlampovich and Miasnikov \cite{BKM}] \label{limit}
The isomorphism problem is solvable for the class of finitely generated fully residually free groups.
\end{corA}

Finitely generated fully residually free groups are also known as `limit groups', and 
are central to the recent work on the Tarski problem by Sela \cite{SelaDio1,SelaDio27}
and also by Kharlampovich and Miasnikov \cite{KM2}.
See \cite{DahComb} and \cite{Alibegovic} for proofs that limit groups are toral relatively hyperbolic.  More generally, groups acting freely on $\R^n$-trees (see \cite{Guirardel})  form a wider class of  toral relatively hyperbolic groups.

The proof of Corollary \ref{limit} in \cite{BKM} relies heavily on the results from \cite{KM}, a very long paper in which it is shown that
the Grushko and JSJ decompositions of finitely generated
fully residually free groups can be effectively computed.  In turn, the 
results of \cite{KM} rely on previous work of Kharlampovich and Miasnikov.  The geometric methods in this paper (and in \cite{DahG_freeprod}) recover the main result of \cite{KM} (namely
\cite[Theorem 1, p.3]{KM} -- the effective construction of the JSJ 
decomposition of a freely indecomposable limit group; see Theorem
\ref{FindJSJ} below)
as well as Corollary \ref{limit} (in the more general setting of
toral relatively hyperbolic groups).

As hinted at in the paragraph above,
a key step in our proof of Theorem \ref{Toral} is the algorithmic construction of the (primary) JSJ decomposition of a freely indecomposable toral relatively hyperbolic group:

\begin{thmA} \label{FindJSJ}
There is an algorithm
which takes a finite presentation for a freely indecomposable
toral relatively hyperbolic group, $\Gamma$ say, as input
and outputs a graph of groups which is a primary JSJ decomposition 
for $\Gamma$.
\end{thmA}

(See Section \ref{JSJ-RH} for a discussion of JSJ decompositions,
and Theorem \ref{JSJTheorem} for a statement of the properties
of the primary JSJ decomposition.  In the case of torsion-free hyperbolic
groups, the primary JSJ decomposition is just the essential JSJ decomposition,
as defined in \cite{SelaGAFA}.)

For torsion-free hyperbolic groups, such an algorithm (to find
the essential JSJ decomposition) is due to Sela (unpublished).  
This result is of independent interest and should be useful for many
other applications.  For example, the automorphism group of a
hyperbolic group can
be calculated from the JSJ decomposition (see Sela \cite{SelaGAFA} and Levitt \cite{Levitt}), and a similar analysis applies to toral
relatively hyperbolic groups.  Also, the JSJ decomposition is one of the key
tools in the above-mentioned work on the Tarski problem.  Thus to be able
to effectively find the JSJ decomposition is an important first step for many
algorithmic questions about the elementary theory of free (and possibly
torsion-free hyperbolic or toral relatively hyperbolic) groups.

One of the most important classes of relatively hyperbolic groups (and one which is not covered by either Corollary \ref{hyp} or Corollary \ref{limit}) consists of the fundamental groups of finite-volume hyperbolic manifolds.  In this case the parabolic subgroups are virtually abelian, but not necessarily actually abelian, so this class is not covered by Theorem \ref{Toral} either.  However, these groups have well behaved finite index subgroups, and we can use their properties along with Theorem \ref{Toral} to solve the isomorphism problem for this class (see Section \ref{HypSection}).  Using Mostow-Prasad Rigidity, this implies

\begin{thmA} \label{HomeoHyp}
The homeomorphism problem is solvable for finite volume hyperbolic $n$-manifolds, for $n \ge 3$.
\end{thmA}

Another natural class of relatively hyperbolic groups is the class of fundamental groups of finite-volume manifolds with pinched negative curvature.  In this case the Margulis Lemma implies that the parabolic subgroups are virtually nilpotent.

Although the isomorphism problem is solvable for nilpotent groups \cite{Segal}, the universal theory is not in general decidable \cite{Romankov}.  This is a significant barrier to implementing our method for these groups.

The presence of torsion seems to be a real problem.  In particular
the isomorphism problem for hyperbolic groups with torsion remains
open.

\begin{acknowledgement}
  We would like to thank Zlil Sela, Martin Bridson, Vincent Guirardel, Gilbert Levitt, David Fisher,
Igor Belegradek, Andrzej Szczepa\'nski  and Henry Wilton for useful conversations about this paper.
We would also like to thank the anonymous referees whose careful reading helped 
improve the exposition of this paper.
\end{acknowledgement}

\section{Strategy}\label{Strategy}

In this section, we briefly summarize our strategy towards the
isomorphism problem.  In outline, it is similar to Sela's approach
for torsion-free hyperbolic groups, but we make several significant
simplifications and generalizations.

 Each finitely generated group admits a Grushko decomposition
as a free product, in which the number of factors and their isomorphism
types is determined by the group.  Thus, to solve the isomorphism
problem it is sufficient to compute a Grushko decomposition 
(an algorithm for this was proposed by Gerasimov for hyperbolic groups, and in \cite{DahG_freeprod} for toral relatively hyperbolic groups), and then
solve the isomorphism problem for freely indecomposable groups. 
In this light, the main result of this paper is the following (see Definition \ref{def:RH} for a definition
of toral relatively hyperbolic groups):

\begin{thmA}
The isomorphism problem is solvable for freely indecomposable toral relatively
hyperbolic groups.
\end{thmA}
Henceforth, we will always assume that our groups are freely indecomposable.  Since the two-ended case
is straightforward, the interesting case is when the groups are one-ended.

So, let $H_1$ and $H_2$ be one-ended toral relatively hyperbolic groups.
Given a finite presentation for each of them, we want to compute 
canonical JSJ splittings for each, and to compare them as graphs. Then,
for each isomorphism of graphs, we want to compare the vertex groups that are matched by the graph isomorphism,
relative to the  subgroups of their adjacent edges (marked by the inclusion maps). If we can find an isomorphism of graphs, 
and isomorphisms between matched vertex groups, which respect the isomorphisms on the subgroups corresponding to the groups of adjacent edges, then the two groups are isomorphic. 
On the other hand, if we can check that for any isomorphism of graphs, at least one pair of matched vertex groups are not isomorphic (relative to marked adjacent edge groups), then by the properties of the JSJ decomposition, the groups are not isomorphic.

 There are already subtleties hidden here (that may be skipped on
first reading). We need to make clear which JSJ we will be aiming at.
For torsion free hyperbolic groups, it is the essential JSJ
decompositions studied by Sela in \cite{SelaGAFA}.\footnote{Although he never published his solution to the isomorphism problem for torsion-free hyperbolic groups, the original purpose of the essential JSJ decomposition
was this solution.} Such a decomposition can be
described as follows. 
We say that a subgroup is {\em elliptic} in a splitting if it fixes a point in the Bass--Serre tree of the splitting. 
An {\em essential JSJ splitting} is a maximal reduced universally elliptic splitting over
maximal cyclic subgroups (universally elliptic means that every edge
group is elliptic in any other splitting over a maximal
cyclic subgroup). This is different from ``the'' canonical JSJ
decomposition of a hyperbolic group, as constructed by Bowditch \cite{B_JSJ},  that takes into account every universally elliptic
splitting over (virtually) cyclic subgroups, not necessarily maximal.  The reason we discard the splittings
of the form $G*_C Z$ in which $C$ is of finite index in $Z$ is that
there is no Dehn twist over them, and that, in some sense, these splittings
do not  produce many outer automorphisms (we return to this issue later).
  For toral relatively hyperbolic groups, we use 'primary' JSJ
decompositions, that can be similarly described as maximal, reduced,
universally elliptic splittings over direct factors of maximal abelian
subgroups  (see Definition \ref{PrimarySplitting} for the definition of `primary' splittings, and Section \ref{JSJ-RH} for more discussion, results and proofs about the primary JSJ decomposition).  This primary JSJ decomposition is very similar to the abelian JSJ decomposition of limit groups constructed by Sela in \cite[Section 2]{SelaDio1}. 
 Another difference from Bowditch's JSJ decomposition is that the JSJ we use might not be unique.  However, any two of them are related by
boundedly many moves. Finally a major feature of a JSJ decomposition is
that it gives a decomposition of the group into rigid and ``flexible"
vertex groups. For Bowditch's canonical JSJ, the flexible vertex groups are (in
the torsion free case)  fundamental groups of surfaces with boundaries
(the groups of the boundary components being exactly the adjacent edge groups). Here,
because we disregard splittings over non-maximal cyclic groups, the flexible groups are larger. 
Namely, to each boundary component of a
surface, one can possibly find a larger cyclic group amalgamated. The flexible groups
are not surface groups, but "socket" groups, {\it i.e. } surfaces to
which, for every (conjugacy class of) element representing a boundary
component,  a certain root has been added. 
The advantage of the JSJ decomposition which we use over Bowditch's is that it is well suited
to algorithmic questions.  

We now briefly
discuss how to algorithmically find a primary JSJ splitting of a one-ended toral relatively
hyperbolic group, $H_1$ say.  

First, we want to compute a maximal primary splitting, which is a
splitting over direct factors of maximal abelian subgroups, in which any
vertex group admits no further refinement (no primary splitting in which
its adjacent edge groups fix points in the Bass--Serre tree of the new splitting). 

This maximal splitting is certainly larger than
necessary for the JSJ (it has more edges) but this is our only way to
be sure that we collect all the universally elliptic primary splittings
(among probably many that are not universally elliptic). 
How do we compute a maximal primary splitting? It is easy to
enumerate the graph of groups decompositions of $H_1$ by Tietze transformations on
its presentation. To recognize when a given splitting is primary is not
hard, but involves solutions to several algorithmic problems in $H_1$
such as the word problem, the root problem, and satisfiability of
equations (see Theorem \ref{t:FindSplit}).  It only remains to  get a certificate that  a
particular splitting is maximal, that is, a certificate that vertex groups have no further primary splittings compatible with
their adjacent edge groups.  It is worth remarking that the accessibility results of 
Bestvina and Feighn \cite{BFAccess} prove that there exists a maximal primary splitting.

To this end, recall that, if $G$ is  torsion free hyperbolic, then there is a
now well known equivalence (essentially due to Bestvina, Paulin and Rips, see
  \cite[Theorem 9.1]{SelaIso} and Remark \ref{rem;BPR} below) 
between the following properties:   

\begin{enumerate}
\item\label{BPR:1} $G$ {\it does not } have any nontrivial essential splitting,
\item\label{BPR:2} $Out(G)$ is finite, 
\item\label{BPR:3} there is a finite subset $\mathcal{B}
\subset G$ on which only finitely many non-conjugate endomorphisms of $G$
are injective. 
\end{enumerate}
There are similar statements in the presence of collections of subgroups, and for
toral relatively hyperbolic groups (Theorem \ref{H1Splits} below).

Here is a result that is proved using our main simplification  of Sela's approach. 
Using the characterizations above, it allows us to prove that a toral relatively hyperbolic
group does not admit a nontrivial primary splitting.

\begin{thmA} \label{BasicList}
There is an algorithm with the following properties. 

It takes as input finite presentations for toral relatively hyperbolic
 groups  $G$ and $\Gamma$.
 
 It terminates if and only if there is some finite subset $\mathcal B$ of $G$ so that there are only finitely many
 non-conjugate  homomorphisms from $G$ to $\Gamma$ which are
 injective on $\mathcal B$. 

 In case it terminates, it provides a finite list of homomorphisms which contains a representative of every conjugacy
 class of monomorphism from $G$ to $\Gamma$.
\end{thmA}

 See also Theorem \ref{theo;algo_comp}, which is somewhat stronger.   In particular, the groups are given with a collection of abelian subgroups (direct factors of maximal abelian ones), and the homomorphisms under consideration are only those that are compatible with these collections (see definitions in the next section).

We will comment upon the proof of Theorem \ref{BasicList} later, but let us just mention that this
justifies our interest in splittings producing many Dehn twists, and
therefore for essential or primary JSJ splittings instead of canonical
ones: the absence of a nontrivial compatible primary splitting of $G$ 
characterizes the termination of the algorithm.
 
 For the time being, using the equivalence \eqref{BPR:3} $\Leftrightarrow$ \eqref{BPR:1} of Bestvina, Paulin and Rips
 from above, Theorem \ref{BasicList}
can be used to produce the certificate
we were looking for: that in a given primary splitting, the vertex groups have no further compatible splitting (hence a certificate of maximality of the splitting, among primary splittings).   
To see this, we apply Theorem \ref{BasicList} (or more precisely its version with collections of subgroups) to each vertex
group $V_i$ of our splitting in turn, using $G=\Gamma = V_i$, and using the collection of subgroups consisting of the adjacent edges groups.  
Thus, we can compute a maximal primary splitting.  This is the content of
Theorem \ref{MaxSplit} below.

So far  we have obtained a maximal primary splitting. However this is not yet the
primary JSJ splitting. This decomposition is not
canonical: it is obtained from the  JSJ decomposition by cutting the
pieces corresponding to surfaces or sockets in an arbitrary maximal way.
We need to recognize where the sockets are, as subgraphs of the graph of
groups. This is done by noticing that each vertex group coming from a
maximal primary splitting of a socket must be free of rank $2$,
and associated to a trivalent vertex of the ambient graph, and that the adjacent edge groups are
generated by elements $a,b, a^nb^m$ for some basis $a,b$, and some
integers $n,m$ (we say that the group is a `basic' socket, defined just before Proposition \ref{prop;find_sockets}). 
To picture it, let us just say that it must correspond 
to a pair of pants with boundary components $A, B$ and $AB$, in which $A$ and $B$ are given 
respectively $n$-th and $m$-th roots $a$ and $b$ in the group. 
 We can detect which vertex groups are free, and by
solving well chosen equations in them (see Proposition \ref{prop;find_sockets}), we can detect which of
these are basic sockets. We prove in  Proposition \ref{HaveaJSJ} that, once we have collapsed
the subgraphs of basic sockets, the splitting obtained is a primary JSJ
decomposition of $H_1$.  This proves Theorem \ref{FindJSJ}.

As in the hyperbolic case (see \cite{SelaGAFA}), primary JSJ splittings are
not unique, but can be obtained from one another by a sequence of
boundedly many computable moves (see Remark \ref{r:BoundedMovesonJSJ}). 
Thus we can get every primary
JSJ splitting of $H_1$, and similarly of $H_2$. It remains to compare
them pairwise.

For simplicity let us assume that we have only one primary JSJ decomposition for
each group $H_1, H_2$. Clearly, if the underlying graphs are different, the groups are different. 
But if now they are isomorphic, we
need to check, for each possible isomorphism of graph,  whether one can
extend it to an isomorphism of groups. We are facing the following final
problems:

First, given two socket groups, with a collection of conjugacy classes
of embeddings of cyclic groups, decide whether they are isomorphic by an
isomorphism preserving the collections.

Secondly, given two toral relatively hyperbolic groups with a collection
of conjugacy classes of embeddings of abelian groups, and with no further
compatible primary splitting,  decide whether they are isomorphic by an
isomorphism preserving the collections.

The first point is rather easy, since it is a matter of recognizing the
surface, and the orders of the roots added at each boundary component.

The second is finally the analogue of Sela's original result
\cite{SelaIso} (although with peripheral structures, and for toral relatively
hyperbolic groups). 
Bestvina--Paulin--Rips' characterization says that for these groups, without
compatible splittings, 
the algorithm of  Theorem \ref{BasicList} will terminate for morphisms in both directions. 
This allows us to obtain
two finite lists of morphisms in both directions, both containing
representatives of every compatible monomorphism up to conjugacy. 
It remains to check
whether  pairs of morphisms are inverses of one another, up to conjugacy.
 This is easily done, by testing whether the compositions are conjugate to the
identities (with for instance a solution to the simultaneous conjugacy problem).

Now we would like to come back to the key Theorem \ref{BasicList} and comment upon
it. This theorem, even without considering the families of subgroups, is in fact the heart of the algorithm when one is only
interested in a simpler version of Sela's original result in \cite{SelaIso} (for torsion-free
hyperbolic groups with trivial essential
JSJ splitting). Indeed, Theorem \ref{BasicList} allows us to verify
the absence of essential splittings, and also to solve the isomorphism
problem, as explained just above.  Thus, for Sela's published result, Section
\ref{List} and Remark \ref{rem;BPR} are sufficient. 

Let us now comment upon the proof of Theorem \ref{BasicList}. The theorem asks for a certificate (in an algorithmic way) of
the absence of  compatible $G\to \Gamma$ injective on a given finite set $\mathcal{B}$
and not conjugate to  a map in  a finite given family.  

Let us consider a simpler problem. Certifying the absence  of  compatible $G\to \Gamma$ injective on a given finite set $\mathcal{B}$
and not \emph{in}  a finite given family of maps,    is a matter of deciding that a
certain 
finite system of equations and inequations in $\Gamma$ has no solution. Indeed, in
this system, the unknowns are the
images of the generators of $G$, the equations ensure that they satisfy
(in $\Gamma$)  the  defining relations of $G$, and that the morphism is compatible,
then some inequations ensure the
injectivity on $S$ and other inequations ensure the map is different from
 any of the given maps. Such a problem (deciding whether a system of
equations and inequations in $\Gamma$ has no solution) is solvable by an algorithm
given in \cite{Dah_eq}. 
 
 However, doing that,  we miss the issue about not being \emph{conjugate} to a map in the given list, instead of just
 not being {\em in} the list.  This additional
 requirement is necessary in our approach, because, if 
 there is a monomorphism $G\to \Gamma$,  there are infinitely many different ones (its conjugates), 
 and this fact would prevent us from making a list of them. There seems to be no
 easy way to write equations and inequations whose solutions are morphisms that
 are not conjugate to something in the given list (it might be impossible).

 Our solution to this problem is to look  only for morphisms that are ``short'' in
 their conjugacy class (among those injective on $\mathcal{B}$), in order to
 count  only finitely many morphisms by conjugacy class.  By `short' we mean that
 the maximal length of the image of the generators of $G$ in $\Gamma$ is as small
 as possible in the conjugacy class.  This notion of short may not
  be expressible in terms of equations and inequations, but we can formulate it
 in terms of \emph{rational constraints} that are, in some sense, conditions recognized by a finite
state automaton (see Subsection \ref{SS:Algo}). In \cite{Dah_eq} the first author studied solvability of finite systems
of equations, inequations and certain constraints in toral relatively hyperbolic
groups.
 The main difficulties are then to express the property of being `short' in its
conjugacy class, for a morphism, by the property of membership to certain 
languages, to prove that these languages are rational, and to prove that
there are short representatives of every conjugacy class, but only
finitely many.  The condition for a homomorphism to be considered `short'
is listed in Proposition \ref{prop;Omega}.  
 The proof of the required properties is unfortunately much more
difficult in the relative case than in the hyperbolic case.  Therefore the proof
for the hyperbolic case is in Section \ref{List} while the proof in the toral
relatively hyperbolic case is deferred to Section \ref{part;proof_gene}.

The main step in proving Theorem \ref{BasicList} is then the technical
Proposition \ref{prop;Omega}, in which we give a list of conditions 
  $\Omega$, expressed in rational languages, on a homomorphism so that
  \begin{itemize}
  \item[$\bullet$] In each conjugacy class of compatible monomorphism, at least one homomorphism  satisfies $\Omega$; but
  \item[$\bullet$] Only finitely many homomorphisms in each conjugacy
  class of compatible homomorphism satisfy $\Omega$.
  \end{itemize}
  
We explain the idea behind the condition $\Omega$ when we introduce it
in Proposition \ref{prop;Omega}.

It is worth mentioning that, in \cite{SelaIso}, Sela did not use
solvability of inequations (let alone rational constraints), which was
not available at that time. He had to make an intricate construction
using free monoids with paired alphabet, and could not use the standard
construction of the limiting $\mathbb{R}$-tree from the argument of
Bestvina and Paulin (this is his constructions of ``states''). Using rational
constraints, we bypass this construction and instead always use the standard
limiting construction. We believe that our contribution provides a neat
simplification of his method.

\begin{remarknum}[Commentary on the exposition] \label{rem;expo}
Although we make some major simplifications and generalizations,
we believe the proof of Theorem \ref{Toral} largely follows
the outline of Sela's unpublished proof of Corollary \ref{hyp}.
However, there are serious technical difficulties to our 
generalization to the relatively hyperbolic setting.

Most striking is the proof of Theorem \ref{theo;algo_comp}, which is the main tool
used to prove the certain primary splittings  are maximal.
In the toral case, the proof is twice as long as the 
proof in the hyperbolic case.   Also, the general proof obscures the hyperbolic case.  Thus,
we have chosen to include the proof of the hyperbolic case in Section \ref{List} inside the main body of the text, and leave
the proof of the general case until Section \ref{part;proof_gene}.

Similarly, as remarked above, the JSJ decomposition we use 
(see Theorem \ref{JSJTheorem}) is exactly the {\em essential} 
JSJ decomposition in the hyperbolic case, which was constructed in \cite{SelaGAFA}.  By now,
this is standard, and well-known to those familiar with JSJ decompositions.  
There are a few differences in the toral relatively hyperbolic case,
and the existence of our JSJ decomposition must be proved.  Thus
we prefer to leave the proof of Theorem \ref{JSJTheorem} in the
relatively hyperbolic setting, as well as a discussion of JSJ decompositions, to Section \ref{JSJ-RH}.

We justify this slightly unusual method of organizing the paper
as follows.  The solution to the isomorphism problem for torsion-free hyperbolic groups is itself a very important result which, although 
known to Sela for many years, was never made public and
has not appeared in print before now.  
Thus, for the reader interested only in this case, 
Sections \ref{Prelim}--\ref{ProofMainThm} provide a complete
solution (for freely indecomposable groups).

However, in broad outline, our solution to the isomorphism problem
for toral relatively hyperbolic groups is almost unchanged when
restricted to the torsion-free hyperbolic case.  Therefore, 
Sections \ref{Prelim}--\ref{ProofMainThm} also provide a complete
solution to the isomorphism problem for toral relatively hyperbolic
groups, {\em except} for the proofs of Theorem \ref{theo;algo_comp}
and Theorem \ref{JSJTheorem}.

The middle sections of this paper are therefore somewhat schizophrenic 
in nature.  We 
state most results in the (more general) setting of toral relatively
hyperbolic groups, and make numerous comments about the
hyperbolic case.  If the reader is only interested in the hyperbolic
case, they should always keep in mind that torsion-free hyperbolic
groups are, in particular, toral relatively hyperbolic.
\end{remarknum}

\section{Toral relatively hyperbolic groups} \label{Prelim}

In this section we gather definitions and basic tools.  Subsection
\ref{ss:RH} contains basic definitions and the definition of (toral) relatively
hyperbolic groups.  Subsection \ref{SS:Algo}
collects some known results about algorithms for hyperbolic and relatively
hyperbolic groups.  Subsection \ref{ss:splittings} contains the basic terminology
of splittings of groups, and the definition of {\em primary} splittings, 
the class of splittings we consider in this paper.  Subsection \ref{ss:vertexRH} proves a 
simple result about the vertex groups of primary splittings of toral relatively hyperbolic 
groups.

\subsection{Relatively hyperbolic groups} \label{ss:RH}
Relatively hyperbolic groups were first introduced by Gromov in \cite{Grom}.
An alternative definition was given and studied by Farb in \cite{Farb}.  
In \cite{B_rh}, Bowditch gave further
definitions and studied further aspects.  By now, it is known that all definitions
are equivalent (see \cite[Appendix]{Dah_thesis} for example).  We adopt the
now-standard convention that `relatively hyperbolic' means in the sense of Gromov
(which is equivalent to Farb's `relatively hyperbolic with BCP').  The definition we
give below in Definition \ref{def:RH} is due to Bowditch \cite{B_rh}.

\begin{convention}
Throughout this paper, we will assume that graphs are given a metric where
edge lengths are $1$. The paths are always continuous.
\end{convention}

\begin{defn}
Let $X$ be a graph, $v$ a vertex of $X$ and $e_1$, $e_2$ edges
of $X$ which contain $v$.  The {\em angle at $v$ between $e_1$ and 
$e_2$} is the distance between (the remaining parts of) $e_1$ and $e_2$ in 
$X \setminus \{ v \}$ (endowed with its path metric).  We allow (though will
never be interested in) the possibility that this distance is $\infty$.
  
Denote this quantity by  $\Ang_v(e_1,e_2)$.  If $\gamma$ is a path
containing a vertex $v$ (which intersects only two edges of $\gamma$) 
then $\Ang_v(\gamma)$ is the angle at $v$ between the edges of 
$\gamma$ intersecting $v$.

If $\gamma$ is a path, then $\MaxAng(\gamma)$ denotes the maximum
angle of $\gamma$ at all its vertices.
\end{defn}

\begin{defn}
The {\em cone} centered at an un-oriented edge $e$ of radius $r$ and angle
$\theta$ is the set of vertices $w$ such that there is a 
path whose first edge is $e$ (with any orientation) and that ends at $w$, which is of length at most $r$
and has maximal angle at most $\theta$.     

The {\em conical neighborhood} of a path $p$, denoted 
$ConN_{r,\theta}(p)$, is the union of the cones of radius $r$ and 
angle $\theta$, centered on the consecutive edges of $p$.
\end{defn}

Here is an easy expression of hyperbolicity:

\begin{prop}(Conically thin triangles)

  If the graph $X$ is $\delta$-hyperbolic, then for any geodesic triangle,
  every side is contained in the union of the conical neighborhood of radius
  and angle $50\delta$ of the two others. 
\end{prop}

We briefly  sketch the proof.
Consider an edge $e$ on one side of the triangle, far from the vertices and from the center of the triangle. Consider the loop obtain by moving $5\delta$ forward on the given side, going by a path of length $\delta$ on another side, moving backward  on this new side for $10\delta$, and going back to the side of $e$, by a path of length $\delta$, and finally, going back to $e$ (by triangular inequality, this last path has length at most $7\delta$). If one can extract from this loop a simple loop containing $e$, then it has length at most $24\delta$, ensuring that all its angles are at most $24\delta$, and that $e$ is in fact in the cone of radius and angle $24\delta$ of every of its edges. It is then easy to check that this loop must contain an edge of the second side of the triangle. 
 If one cannot extract a simple loop containing $e$, this means that $e$ appears twice in our loop, and by triangular inequality, it must be on the other side of the triangle. The conclusion follows. The same argument can be adapted if $e$ is close to a vertex of the triangle, or to the center.

\begin{defn}
For numbers $L, \lambda,\mu$, an {\em $L$-local, $(\lambda,\mu)$-quasi-geodesic} in $X$ is a
path parameterized by arc length $p:[0,T] \to X$ such that for all $x,y \in 
[0,T]$, with $|x-y|\leq L$, we have 
$d(p(x),p(y)) \geq \frac{|x-y|}{\lambda} - \mu$. The path
is a {\em $(\lambda,\mu)$-quasigeodesic} of $X$ if the inequality holds for all $x,y$.
\end{defn}

We also need a mild assumption on backtracking, stated in terms of `detours', which
we now define.

\begin{defn} \label{d:detour}
A subpath $w$ of a path $p=p_1wp_2$ is an 
        $r$-\emph{detour},  if  (i) it is a (not necessarily simple) 
        loop in $X$;  
        (ii) $p_1$ and $p_2$ are not 
        empty paths; and (iii) the angle between the last edge of $p_1$ and the
        first of $w$ (resp. the last of $w$ and the first of $p_2$) is 
        greater than $r$. 
\end{defn}

\begin{remarknum} \label{r:detour}
Note that if a path $p$ in a  graph is a
$(\lambda,\epsilon)$-quasigeodesic, then $\epsilon$ is an explicit bound
on the length of a simple loop which is a subpath of $p$.  
The condition that $p$ has no $r$-detours can thus be detected locally.
\end{remarknum}

        \begin{prop}(Conical stability of  quasi-geodesics without detours)
        \cite[Proposition  2.7]{Dah_eq} \label{prop;stab_con} 
          Let $X$ be a hyperbolic graph. 
          Given $r, \lambda, \mu >0$, there is an explicit 
          constant $\epsilon$ such that for any 
          $(\lambda,\mu)$-quasi-geodesic  $\rho$ which is 
          without $r$-detours, and any geodesic $[x,y]$ joining
          the endpoints of $\rho$,  if $x\neq y$,  the path $\rho$
           is contained in the $(\epsilon,\epsilon)$-conical
          neighborhood of $[x,y]$.

          Moreover, if $[x,y]$ contains a vertex $v$ so that $\Ang_v([x,y]) \ge 4\epsilon$
          then
          \begin{enumerate}
          \item $v \in \rho$;
          \item If $e_1$ (respectively $e_2$) is the first edge in $\rho$ (resp. $[x,y]$) which 
          contains $v$ then $\Ang_v(e_1,e_2) \le \epsilon$; and
          \item The path $\rho$ contains a pair of consecutive edges intersecting in $v$
          which make an angle at most $\epsilon$ with the first and second edges of
          $[x,y]$ which contain $v$.
          \end{enumerate}
        \end{prop}

        The second assertion is a consequence of the first, since if
        $\Ang_v[x,y] >4\epsilon$, then it is easy to check that by triangular inequality, 
        $ConN_{\epsilon,\epsilon}[x,v] \cap  ConN_{\epsilon,\epsilon}[v,y] =
        \{v\}$.

\begin{defn}\label{def;fine}
A graph is {\em fine} if every cone is finite.
\end{defn}

Although this is not Bowditch's original definition, formulated in \cite{B_rh} (which ask, for each edge $e$, for only finitely many edges $e'$ with prescribed angle with $e$), this is equivalent to it, by Koenig's lemma  (see \cite[Appendix]{Dah_thesis} for instance).

We turn now to groups.

\begin{defn} (Coned-off graph)

Given a  group $H$, with a finite generating set $S$, and finitely
generated subgroups $G_1, \dots, G_p$,
the {\em coned-off graph of $H$ relative to $G_1, \dots, G_p$, with
respect to $S$}, which we denote by $\widehat{Cay}(H)$,  
is the graph obtained from the Cayley graph of $H$ (with generating set $S$) by adding,
for each left coset of each $G_i$, a vertex $v$, and edges from $v$ to every
element of the coset (see \cite{Farb}). 
The coset corresponding to a vertex $v$ is then denoted $Coset(v)$.
\end{defn}
  
  \begin{remarknum}
  In coned-off graphs, large angles will occur only on vertices of infinite valence.
  Moreover, suppose that $v$ is an infinite valence vertex in the coned-off graph,
  and that $e_1, e_2$ are edges adjacent to $v$.  Then the word
  distance between the finite valence vertices of $e_1$ and $e_2$ gives
  an upper bound to the angle between $e_1$ and $e_2$ at $v$.
\end{remarknum}

\begin{defn} \label{def:RH} [Relatively hyperbolic groups]

  A finitely generated group $H$, endowed with a family 
  $\{ G_1, \dots, G_p \}$, of finitely generated subgroups,  
  is {\em hyperbolic relative to $\{ G_1, \ldots , G_p \}$}
   if some (hence any) associated coned-off graph is 
  hyperbolic and fine.   

  A finitely generated group $H$ is {\em toral relatively
  hyperbolic}, if it is torsion-free, and hyperbolic relative to a family
  of finitely generated, noncyclic, free abelian subgroups.

\end{defn}

  In case $H$ is a toral relatively hyperbolic group, 
  the family of maximal parabolic subgroups is canonical (up to conjugacy).

  The definition  of relatively hyperbolic given here is equivalent to
  the other classical definitions of the literature (see
  \cite[Appendix]{Dah_thesis} for instance).

\begin{examples}
\begin{enumerate}
\item Hyperbolic groups are hyperbolic relative to an empty collection
of parabolics, and are toral if torsion-free;
\item Limit groups are toral relatively hyperbolic (see \cite{Alibegovic} 
and \cite{DahComb}), and more generally, groups acting freely on $\R^n$-trees (see \cite{Guirardel});
\item Fundamental groups of geometrically finite hyperbolic manifolds
are hyperbolic relative to their cusp subgroups, and are toral if the
cusps are all homeomorphic to $T \times \R^+$, for a torus $T$.
\end{enumerate}
See \cite{Farb}, \cite{B_rh}, \cite{DS}, \cite{HK} and \cite{Osin} for further examples.
\end{examples}

\begin{remarknum} \label{rem;hyprelhyp}
Suppose that $G$ is hyperbolic relative to a family $\mathcal{P}$ of subgroups, and that some of the subgroups in $\mathcal{P}$ are hyperbolic.  Let $\mathcal{P}'$ be the
non-hyperbolic groups in $\mathcal{P}$.  Then $G$ is also hyperbolic relative to
$\mathcal{P}'$   (see for instance \cite[Theorem 2.37]{Osin}).   

Therefore, it is no restriction in the definition of toral to assume that the abelian 
parabolics are noncyclic.  It does, however, make certain arguments easier, and
so it is a useful assumption.
\end{remarknum}

  \subsectionind{Algorithmic tools} \label{SS:Algo}

    \subsubsection{General background}
     We recall a few existing algorithms related to hyperbolic groups
     and relatively hyperbolic groups. We will use many of these algorithms
 in the sequel, often without mention.  Of particular importance for our work is
 the discussion of the `uniformity' of the algorithms in Subsection \ref{ss:uniform}.

     First of all, the first and the second Dehn problems 
     (the word problem and the conjugacy 
    problem) are 
    decidable in a hyperbolic group,  
    (this is due to Gromov, see \cite{Grom}, or \cite{CDP}).
 
    In a toral relatively hyperbolic group, these problems are also decidable
    (the word problem is due to Farb, \cite{Farb} and the conjugacy problem to Bumagin \cite{Bumagin}, both of whom proved more general results; see also 
Rebbechi \cite{Rebbechi}).

    There is  an algorithm deciding whether a given 
    element is primitive 
    ({\it i.e.} not a proper power), given by Lysenok \cite{Ly}  for hyperbolic groups.
    For toral relatively hyperbolic groups, Osin \cite{Osin} proved that it is possible
    to decide if a given element which is not conjugate into a parabolic subgroup has
    a nontrivial root.

    An essential tool for our work is the decidability of the existence 
    of solutions to equations, and systems of equations. By a {\em system
    of equations} over a group $G$, we mean a family of words with letters
    which are either (i) elements of $G$ (known as {\em coefficients}); or (ii)
    elements of some fixed list of unknowns $x_1,x_2, \ldots$.
    
    A {\em solution} to a system $\Sigma$ of equations is a collection of elements
    $g_1,g_2, \ldots$ in $G$ so that the substitution $g_i = x_i$ yields $\sigma = 1$
    for all $\sigma \in \Sigma$.  We will only be concerned with finite systems of
    equations, so that only finitely many unknowns appear in $\Sigma$.
A {\em system of equalities and inequalities} is a pair $\Sigma_1, \Sigma_2$, and
a {\em solution} is a substitution which yields $\sigma = 1$ for all $\sigma \in \Sigma_1$
and $\eta \neq 1$ for all $\eta \in \Sigma_2$.

Grigorchuk and Lysenok \cite{GL}, \cite{Ly} proved that it is possible to decide
if quadratic equations over a hyperbolic group have a solution.
Rips and Sela \cite{RS} extended this to arbitrary finite systems of 
     equations (with coefficients)
    over a torsion-free hyperbolic group.
    The analogous result for toral relatively hyperbolic groups is proved in
    \cite{Dah_eq}.

     Recently, there have been improvements of this algorithm. 
     It is proved by Sela \cite{SelaDio8}, and by the first 
     author \cite{Dah_eq}  
     with another method, that the existence of solutions of 
     finite systems of equations and inequations, 
     with coefficients, over a torsion-free hyperbolic group is decidable.
This is equivalent to saying that the universal theory of a torsion-free
hyperbolic group is decidable.  Note that Makanin \cite{Makanin} proved
that the universal theory of a free group is decidable.
     In \cite{Dah_eq}, it is proved that toral relatively hyperbolic
     groups have a decidable universal theory.

In fact, this is not enough for our needs, we will need one more refinement  
     of this algorithm. We want to  allow  
     \emph{constraints} together with equations and inequations. We
     briefly describe what this means.
     
     Let $H$ be a group generated by an ordered finite set $S$.  The set of
     {\em normal forms} of $H$ consists of a choice of word in $(S)^*$ for each
     $h \in H$ which is of minimal length amongst all representatives and is ShortLex 
	({\it i.e.} least amongst those shortest length representatives).  This notion is most useful for free
     groups, and free products of free abelian groups, which is where we will use it.
     See Remark \ref{rem:lang} below.

\begin{defn}
     A {\em normalized rational language} in  $H$ is a subset $A \subset H$ such that 
     if $\mathcal L$ is the set of normal forms for the elements of $A$ then
     $\mathcal L$ is a regular language
     (that is, a language recognized by a finite state automaton).  
     A {\em normalized rational constraint} on an 
     unknown is the requirement that this unknown belongs to some 
     specified rational language (see for instance \cite{DGH} and \cite{Diekert_Muscholl} 
        for more on this subject).
        \end{defn}

     It is decidable whether a finite system 
     of equations and inequations, with rational constraints, 
     in a free group admits a solution 
     (V.~Diekert, C.~Guti\'errez, 
     C.~Hagenah, \cite{DGH}), and also in a free product of 
     abelian groups (V.~Diekert, and A.~Muscholl \cite{Diekert_Muscholl}).    
     This problem is not known to be decidable 
     in a torsion free hyperbolic group, not 
     to mention toral relatively hyperbolic groups. 
     Nevertheless, some weaker form of this is known, as we now explain.

     \subsubsection{Decidability with constraints}  \label{para;decision}
In this paragraph we recall one of the main results of \cite{Dah_eq}.  At the
beginning of Section \ref{List}, we will describe briefly how this result is used
in this paper.

     Let $H$ be a toral relatively hyperbolic group (for example: a torsion
     free hyperbolic group), endowed with some finite presentation, whose
     finite symmetrical generating set is denoted $S$. Let $G_1,\dots G_p$ be
     representatives of conjugacy classes of maximal parabolic subgroups. We
     denote by $F_S$ the free group on $S$ (respecting inversion), and by $F$
     the free product $F= F_S * G_1 * \dots * G_p$. In the case of a
     hyperbolic group, of course, $F= F_S$. We endow $F$ with a natural
     generating set, denoted $\mathcal A$, which is the union of $S$ with a basis for each of the abelian subgroups $G_i$.
   
   \begin{remarknum}  \label{rem:lang}
The set of words in $\mathcal A^*$ which are normal forms of elements of $F$ is easily
seen to be a regular language, which we write $\mathcal L_{N}$.  This is equivalent
to saying that the set $F$ is a normalized rational language inside itself.
There is a natural one-one correspondence between elements of $F$ and elements
of $\mathcal L_N$.  Therefore, when we define a sub-language of $\mathcal L_N$,
we will speak of a language of elements of $F$.  In case this sub-language
is regular (as a word language), it will be a normalized rational language (as a subset of the group).
\end{remarknum}

\begin{remarknum} \label{rem;labels}  There is a natural quotient map $\nu: F \to H$, but we can be more precise.
     Any element $f$ in $F$ \emph{labels} a (continuous) path in $\widehat{Cay}(H)$, from $1$ to $\nu(f)$ 
	(seen as vertices of 
     $\widehat{Cay}(H)$),  by the
     mean of its normal form in the free product, in such a way that 
     any subword in $F_S$ labels
     consecutive edges in the Cayley graph, and any maximal subword  %
	with letters in a basis of some  abelian
     factor group $G_i$ labels two consecutive edges, with a common  vertex of 
     infinite valence, corresponding to the relevant coset of $G_i$. 
\end{remarknum}

\begin{ex}
Assume $s_1,s_2,s_3 \in S$, and $s_1s_2s_3= g\in
G_1$ a parabolic subgroup of $H$, and $g\neq 1$.  Then, in $F$,  the element $s_1s_2s_3g^{-1} \in F$ is non-trivial (the $4$ letters word given is its normal form) and labels a closed loop in $\widehat{Cay}(H)$. The three first letters $s_1,s_2$, and $s_3$  label three edges in $Cay(H)$, from $1$ to $(s_1s_2s_3)$,  and the last letter $g^{-1}\in G_1$ labels two consecutive edges, from $(s_1s_2s_3)$ to  $(s_1s_2s_3g^{-1})=1$ via the vertex of the left coset $s_1s_2s_3 G_1$. 
  \end{ex}

     Let $\delta$ be the hyperbolicity constant of the coned-off Cayley graph
     $\widehat{Cay} (H)$ (in case $H$ is hyperbolic, this is just $Cay( H)$).  See Subsection \ref{ss:uniform} below for a discussion of
     algorithmically finding $\delta$.

     Recall that in $\widehat{Cay}(\Gamma)$, cones are finite (Definition \ref{def;fine}), and that there are finitely many orbits of edges. 
     We define the constants $L_1$ and $L_2$, as in \cite{Dah_eq}, 
     to be $L_1=10^4 \delta M$, and
     $L_2 = 10^6\delta^2 M$, where $M$ is a bound to the 
     cardinality of the cones of radius and
        angle $100\delta$ in 
        $\widehat{Cay}(\Gamma)$ (the ball of radius $100\delta$ in the
     case of hyperbolic group). Let $r= 10^9 D$, where $D$ is a 
     fellow-traveling constant between $1000\delta$-quasi-geodesics in 
     $\widehat{Cay}(\Gamma)$.  Note that a suitable value of $D$ can be computed
     explicitly in terms in $\delta$.

We will consider `lifting' equations in $H$ to equations in $F$.  Therefore, the
following definition will be useful.
\begin{defn}
Suppose that $h \in H$.  A {\em representative} of $h$ is an element of $F$
(considered as a normal form), so that $\nu(f) = h$.
\end{defn}

\begin{remarknum} \label{rem:L1prime}
     By a classical property of hyperbolic spaces, (see for instance
     \cite[Chapter 3]{CDP}), there are explicit constants $L$, 
     $L'_1$ and $L'_2$, such that any $L$-local, $(L_1,L_2)$-quasigeodesic in
     $\widehat{Cay}( H)$ is a global $(L'_1,L'_2)$-quasigeodesic. 
\end{remarknum}
     
\begin{defn} \label{def;lang_L}
Define $\mathcal{L} \subset F$, to be the language of the elements of $F$
     which label paths in $\widehat{Cay} (H)$ that are $L$-local,
     $(L_1,L_2)$-quasigeodesic, and do not have any $r$-detours.
\end{defn}
As noted in Remark \ref{r:detour} the existence of $r$-detours can
be recognized locally (on subpaths of length $L'_2$).  Of course, this is also true of being an $L$-local
$(L_1,L_2)$-quasigeodesic.  Thus, the properties defining $\mathcal{L}$ 
are purely local, and this makes $\mathcal{L}$ a
     normalized rational language (this is proved in 
     \cite[Proposition 2.5]{Dah_eq}).
The interest of the language $\mathcal{L}$ is that all its elements satisfy
the property of being globally $(L_1',L_2')$-quasigeodesic, for constants
$L_1'$ and $L_2'$ that are explicit in $\delta$, $L_1$, and $L_2$.

     \begin{defn}
     Define the language $\mathcal{L}_0 \subset \mathcal{L}$ to
     consist of the
     elements of $\mathcal{L}\subset F$ that, after projection to $H$,
     represent the trivial element of $H$.
     \end{defn}

     The language $\mathcal{L}_0$ is a finite subset of $\mathcal{L}$, and is computable. It
     is a normalized rational language, and so is $\mathcal{L}\setminus
     \mathcal{L}_0$  (we refer to \cite{Dah_eq} for proofs).

                \begin{thm}\label{theo;decision} 
                \cite[Proposition 5.3]{Dah_eq}
        
                  Let $H$ be a toral relatively  hyperbolic group, and $F$ as described above.
                  There is an algorithm that, 
\begin{itemize}
\item[$\bullet$]                  takes as input a system of equations with coefficients in $H$, 
                  and unknowns
                  $x_1,x_2,\dots ,x_n$, and normalized rational  languages
                  $\mathcal{L}_i \subset \mathcal{L} \subset F$ 
                  (with notations as above),
                  $i=1,\dots,n$ ;
\item[$\bullet$]                   always terminates; 
\item[$\bullet$]                  answers `yes'   if  there 
                  is a  solution such that every representative in 
                  $\mathcal{L}$ of $x_i$ is in $\mathcal{L}_i$, 
                  for all $i$; and 
\item[$\bullet$]                  answers `yes' only if there exists some solution with, 
                  for all $i$, some 
                  representative in $\mathcal{L}_i$.
\end{itemize}
                \end{thm}
Note that in case the algorithm says `no', it could still be that there is some solution
which has a representative in $\mathcal{L}_i$ for each $i$ (we just know that there
is no solution with {\em every} representative in $\mathcal{L}_i$).  This technical
point leads to certain complications in our work in Sections \ref{List} and \ref{part;proof_gene}.

         The interest of Theorem \ref{theo;decision} is that, for example, one can choose the languages $\mathcal{L}_i =
        \mathcal{L}\setminus \mathcal{L}_0$, for certain $i$,  to 
        embody  possible inequations.  To illustrate our basic approach, we recall the
        proof of this from \cite{Dah_eq}.
        
\begin{thm} \cite[Theorem 0.1]{Dah_eq}
The universal theory of a toral relatively hyperbolic group is decidable.
\end{thm}
\begin{proof}
Let $\Sigma(\underline{x})$ and $\Lambda(\underline{x})$ be a pair of systems of
equations with coefficients over a toral relatively hyperbolic group $\Gamma$.  We consider whether the system $(\Sigma =1) \wedge (\Lambda \neq 1)$ has a solution.

By adding some extra unknowns and equations, we may assume without loss of generality that each inequation in $\Lambda$ has the form $x_i \neq 1$ for some variable
$x_i$.  Let $T$ be the set of indices so that $x_i \neq 1$ is in $\Lambda$.

Choose languages as follows:  if $i \in T$ then let $\mathcal{L}_i = \mathcal{L} \setminus \mathcal{L}_0$.
Otherwise, let $\mathcal{L}_i = \mathcal{L}$.  Since $\mathcal{L}\setminus \mathcal{L}_0$ is precisely the set of elements which do not project onto the trivial element of $\Gamma$, the algorithm from Theorem \ref{theo;decision} applies to decide whether
the given system of equations and inequations has a solution or not. 
\qed \end{proof}

        However, this is not the only way we will use Theorem \ref{theo;decision}. In fact,
        we need the full force of Theorem \ref{theo;decision}, and one of our difficulties
        is in finding the correct languages $\mathcal{L}_i$ (this is much more difficult
        in the relatively hyperbolic case than the hyperbolic case).  
        See the beginning of Section
        \ref{List} for a discussion of our application of Theorem \ref{theo;decision}.

\subsubsection{Uniformity} \label{ss:uniform}
We are working in the context that we are given a finite presentation
for a group $H$. It is crucial that each algorithm we use can be implemented with only the knowledge of this presentation. This often means that the hyperbolicity constant of the relevant graph should be computable from the presentation.

For hyperbolic groups, this computability was pointed out by  M.~Gromov \cite{Grom}  
(see for instance   P.~Papasoglu \cite{Papa} for a detailed study). For toral
relatively hyperbolic groups, what we need was studied in  \cite{Dah_find}.

        The following result is applied many times throughout this paper, often without mention.
        
        \begin{thm} \cite[Theorem 0.2 and Corollary 2.4]{Dah_find} \label{theo;finddelta}
        
        There is an algorithm whose input is a finite presentation of a toral relatively hyperbolic group $H$, and whose output is 
        a set of representatives of the parabolic subgroups up to conjugacy,
        (given by a finite collection of bases as free abelian groups), the hyperbolicity constant for
        the coned-off Cayley graph, and the list of orbits of simple loops in
        $\widehat{Cay}(H)$ for any given length.  
     
        \end{thm}

  \begin{remarknum} \label{Rem:FindRH}
  Theorem \ref{theo;finddelta} implies that we do not need to be told {\em how} a group
  is toral relatively hyperbolic, merely that it {\em is}.  Once we know it is, we can apply the above
  algorithm to discover a collection of parabolic subgroups, along with bases, and also
  the relevant `data' of the groups relative hyperbolicity.
  \end{remarknum}
        
        \begin{hypremark}
        In case $\Gamma$ is torsion-free hyperbolic, rather than merely toral relatively hyperbolic, Theorem \ref{theo;finddelta} may be replaced by the above-mentioned algorithm of Papasoglu.
        
        For torsion-free hyperbolic groups, the Cayley graph replaces the coned-off graph.
        \end{hypremark}
        
        Once we have the hyperbolicity constant, we can explicitly
        compute finite state automata for the languages which are used below
        (we say more about this later).

\subsection{Splittings of groups} \label{ss:splittings}
Recall the terminology of graphs of groups (see \cite{Serre} and \cite{RipsSelaJSJ}
 for more details).

\begin{defn}
A {\em splitting} of a group is a graph of groups decomposition.  The splitting is called {\em abelian} if all of the edge groups are abelian.

An {\em elementary} splitting is a graph of groups decomposition for which the underlying graph contains only one edge.

A {\em refinement} of a graph of groups $\Lambda$, at a vertex $v\in \Lambda$ (with vertex group $V$), is a splitting $\Xi$ of  $V$   such that all edge groups adjacent to $V$ in $\Lambda$ are contained in
 vertex groups of $\Xi$, up to conjugacy. The refinement is elementary if  $\Xi$ is elementary.
\end{defn}

\begin{defn}
A splitting is {\em minimal} if there is no proper invariant sub-tree of the associated
Bass-Serre tree.

A splitting is {\em reduced} if it is minimal and, for every valence two vertex with two distinct adjacent edges, the inclusions of the edge groups into the vertex group are proper inclusions.
\end{defn}

We will also need the following defintion, from \cite{RipsSelaJSJ} (see also
\cite{BFAccess} for more information).

\begin{defn} \label{d:unfolded}
A splitting $G = A\ast_{C_1} B_1$ is obtained from the splitting $G = A \ast_C B$ by
{\em folding} if $C$ is a proper subgroup of $C_1$ and $B_1 = C_1 \ast_C B$.

A splitting $G = A_1 \ast_{C_1}$ with stable letter $t$ is obtained from the splitting $G = A\ast_C$ (where the inclusions of $C$ in $A$ are $\alpha, \omega : C \to A$)  by
{\em folding} if
$\alpha(C)$ is a proper subgroup of $C_1$ and $A_1 = A \ast_C (tC_1t^{-1})$.

An {\em unfolding} is the inverse operation of folding.  A splitting is {\em unfolded}
if there is no unfolding (note that it is required in the definition of folding that
$C \ne C_1$).
\end{defn}

We now describe the classes of splittings considered in the sequel.

\begin{defn} \label{EssentialSplitting}
Let $\Gamma$ be a toral relatively hyperbolic group.
A  splitting of $\Gamma$ is called {\em essential}
if it is reduced and if  
\begin{enumerate}
\item all edge groups are abelian; and
\item if $E$ is an edge group and $\gamma^k \in E$
for some $k \neq 0$ then $\gamma \in E$.
\end{enumerate}
\end{defn}

In case $\Gamma$ is a hyperbolic group, an abelian splitting is essential
if and only if all edge groups are trivial or maximal cyclic.  In toral relatively
hyperbolic groups, the edge groups  are in fact direct factors of maximal
abelian subgroups of $\Gamma$. Indeed, if $A$ is the maximal abelian subgroup
containing $E$, the second condition says that $A/E$ has no torsion, thus is
free abelian, and therefore there is a section $s: A/E \to A$, and one has $A= E
\rtimes s(A/E)$, which is a direct product since $A$ is abelian.

\begin{defn} \label{PrimarySplitting}
Let $\Gamma$ be a toral relatively hyperbolic group.  A reduced 
splitting of $\Gamma$ is called {\em primary} if it is essential and if
each noncyclic abelian subgroup is elliptic ({\it i.e.} fixes a point) in the Bass-Serre tree of the splitting.
\end{defn}

\begin{defn} [Dehn twists]
Suppose that $\Gamma = A \ast_C B$ is an abelian splitting
of $\Gamma$.  The {\em Dehn twists} in $c \in C$ is the automorphism
of $\Gamma$ which fixes $A$ element-wise and conjugates $B$ by $c$.

Suppose that $\Gamma = A \ast_C$ is an abelian splitting of $\Gamma$.
The {\em Dehn twist} in $c \in C$ replaces the stable letter $t$ of the HNN
extension by $tc$ and fixes $A$ element-wise.
\end{defn}

\begin{defn} [Generalized Dehn twists]
Suppose that $\Theta$ is an abelian splitting of $\Gamma$, and that
$A$ is an abelian vertex group of $\Theta$.
A {\em generalized Dehn twist with respect to $\Theta$} is an automorphism of $\Gamma$
which fixes each vertex group other than $A$ and each edge
group adjacent to $A$ element-wise, and also fixes $A$ as a set (though not
necessarily element-wise, of course).
\end{defn}

\begin{lem} \label{DehnInf}
Suppose that a toral relatively hyperbolic group $\Gamma$ admits a 
nontrivial elementary primary splitting $\Lambda$ with nontrivial edge group.  
Then some Dehn twist or generalized Dehn twist with respect to $\Lambda$
has infinite order in $\mbox{Out}(\Gamma)$.
\end{lem}
\begin{proof}
Maximal abelian subgroups of $\Gamma$ are malnormal (see \cite[Lemma 2.4]{Groves_RH1} for instance, stated as Lemma \ref{malnormal} below). 
Therefore,  
if $A \ast_C B$ is a nontrivial primary splitting of $\Gamma$ then
at least one of $A$ and $B$ is non-abelian.

If both $A$ and $B$ are non-abelian then the Dehn twist in some $c \in C$ has infinite order in $\mbox{Out}(\Gamma)$.

The HNN case is similar (noting that noncyclic abelian edge groups in
$\Gamma$ are elliptic in primary splittings).

Suppose then that $\Gamma = A \ast_C B$ is a nontrivial primary
splitting and that $A$ is abelian. By Definition \ref{EssentialSplitting}, $C$
is a direct factor of $A$, thus there is a basis  $\mathcal A$ of $A$ extending one of $C$,
and since the splitting is reduced,  there is  $a \in \mathcal A$ so that $a \notin C$.  Suppose that
$c \in C \smallsetminus \{ 1 \}$.  The automorphism of $A$ which sends $a$ to
$ac$, and fixes each other element of $\mathcal A$ (hence of $C$) extends to
a generalized Dehn twist, and has infinite order in $\mbox{Out}(\Gamma)$.
\qed \end{proof}

\subsection{Vertex groups are relatively hyperbolic} \label{ss:vertexRH}
Our approach to solving the isomorphism problem is to construct
the JSJ decomposition (see Section \ref{JSJSection} below), and then
attempt to build an isomorphism between two groups by finding
isomorphisms between vertex groups and `gluing' these together
into an isomorphism between the whole groups.

In order to implement this strategy, it is very important that at each stage we are
working with toral relatively hyperbolic groups, in order that we can
apply the results of Subsection \ref{SS:Algo} above (and the other results
in this paper).  In particular, the 
vertex groups of our splittings should be toral relatively hyperbolic (the 
edge groups are always abelian).  The purpose of this subsection is to
prove that this is indeed the case.

 \begin{thm} \label{vertexgpsRH}
Suppose that $\Gamma$ is a toral relatively hyperbolic group, and that $\Lambda$ is a primary splitting of $\Gamma$. Then, every vertex group of $\Lambda$ is toral relatively hyperbolic, and the parabolic subgroups are the intersections of the parabolic subgroups of $\Gamma$ with the vertex group.  

\end{thm}

\begin{proof}
In \cite[Theorem1.3]{Bowditch_periph}, Bowditch proved that the vertex groups
of a peripheral splitting of a relatively hyperbolic group are themselves
hyperbolic relative to the intersection of the ambient parabolic
subgroups. Here peripheral means that the parabolic subgroups are all
elliptic, and the edge groups are all parabolic. We are not yet in this situation.

Let us denote by $\mathcal P$ the collection of maximal parabolic subgroups of $\Gamma$.
Consider the collection $\mathcal{C}$ of edge groups of $\Lambda$, and their
conjugates, that are
maximal cyclic and non-parabolic in $\Gamma$.  By 
\cite[Lemma 4.4]{DahComb} the group $\Gamma$ is hyperbolic relative to
$\mathcal{P}\cup\mathcal{C}$. By assumption, every group in  $\mathcal{P}$ is
elliptic in the splitting $\Lambda$, thus, every group in
$\mathcal{P}\cup\mathcal{C}$ is elliptic. Moreover, by construction, every
edge group of $\Lambda$ is a
subgroup of some group in $\mathcal{P}\cup\mathcal{C}$. This makes the
splitting  $\Lambda$ peripheral for the structure $\mathcal{P}\cup\mathcal{C}$
(in the sense of Bowditch \cite{Bowditch_periph}). Therefore, \cite[Theorem
1.3]{Bowditch_periph} can be applied, thus ensuring that each vertex group  is hyperbolic relative to their intersections with the
groups in $\mathcal P$, and the groups $\mathcal{C}$. The latter being cyclic,
hence hyperbolic, they can be removed from the list of parabolic subgroups
without affecting the relative hyperbolicity (see Remark \ref{rem;hyprelhyp}).
This proves the result.
\qed \end{proof}

\section{An algorithm that lists monomorphisms} \label{List}

In this section we apply Theorem \ref{theo;decision} to prove one of the key
result we use, 
namely Theorem \ref{theo;algo_comp} (see Theorem \ref{BasicList} and related
comments in the section on the strategy).

We need vocabulary about collections of subgroups.

\begin{defn} \label{def:periph}
        A \emph{primary peripheral structure} on a 
        torsion-free group is a
        family of subgroups with the following properties. The family is
        conjugacy-closed and finite up to conjugacy, and the subgroups are
        either  maximal cyclic subgroups, or  
        (not necessarily maximal) 
        non-cyclic abelian subgroups, each of them marked by the choice
        of a basis for one representative of each conjugacy class.
\end{defn}

 For the purposes of input into an algorithm, a primary peripheral
 structure $\mathcal{P}$ will be given by  basis  
 for representatives of $\mathcal{P}$ up to conjugacy.
 
 Of course, in the case of torsion-free hyperbolic groups, a primary
 peripheral structure is a conjugacy-closed collection of maximal
 cyclic subgroups.
 
\begin{remarknum}
        It is important to note that we distinguish the {\em parabolic} subgroups of a 
        relatively hyperbolic group, from the {\em peripheral} subgroups belonging to some primary peripheral structure. 
        Not all peripheral subgroups will be parabolic, and not all parabolics will be peripheral.

In applications, the subgroups in a peripheral structure of a toral relatively hyperbolic group will be the edge groups of a graph of
groups decomposition (and their conjugates), or equivalently the set of edge stabilizers in an action on a simplicial tree.  

\end{remarknum}

\begin{defn} \label{def:compat}
       Suppose that $\mathcal{P}_1$ and $\mathcal{P}_2$ are
       primary peripheral structures on $H_1$, and $H_2$ respectively.

        We say that a homomorphism $\phi:H_1 \to H_2$ is \emph{compatible} 
        with $\mathcal{P}_1$ and $\mathcal{P}_2$, 
        if  the images of each group in
        $\mathcal{P}_1$ commutes with some group in $\mathcal{P}_2$, 
        with the
        added stipulation that the given generator of each cyclic group in
        $\mathcal{P}_1$ is actually sent to a 
	conjugate of a given generator  of an element of $\mathcal{P}_2$.
         Furthermore, we require that $\phi$ is
        injective on the ball of radius $8$ of $H_1$. 
\end{defn}

The final requirement in Definition \ref{def:compat}
allows us to ensure that $[\phi(u),\phi(v)]\neq 1$
for any $u,v \in \mathcal{B}_{H_1}(2)$ which are so that
$[u,v] \neq 1$.  (Here $\mathcal{B}_{H_1}(2)$ denotes the ball of radius $2$ about the
identity element in $H_1$, with respect to the chosen generating set.)   In particular (though we will need more than this), if
$H_1$ is non-abelian then so is its image in $H_2$.

The main result of the section is:

        \begin{thm} \label{theo;algo_comp}
          There is an algorithm that, 
          given a finitely presented group 
          $H_1$ with a solution to its word problem, 
          and a toral relatively  hyperbolic group  $H_2$, 
          each with primary peripheral structures, terminates if there is a finite subset 
          $A$ of $H_1$ so that there are only 
          finitely many non-conjugate compatible homomorphisms that 
          are injective on $A$.  
          In case it terminates, the algorithm provides a finite list of 
          homomorphisms
          containing one representative of every conjugacy class 
          of monomorphism.
        \end{thm}

For comments, motivation, and main ideas, we advise reading the relevant
paragraphs, about Theorem \ref{BasicList}, in the section ``Strategy''.

As mentioned in Remark \ref{rem;expo}, the proof of this result in the hyperbolic
case is much easier than in the relatively hyperbolic case.  Therefore, in this
section we prove Proposition \ref{prop;Omega} only in the torsion-free hyperbolic case
(the general case is proved in Section \ref{part;proof_gene}).

However, assuming Proposition \ref{prop;Omega} in the general case, we
give a complete proof of Theorem \ref{theo;algo_comp} in the relatively
hyperbolic case in this section.

Assumptions in this section are as follows:  In Subsection \ref{part;not} we are
in the general setting of toral relatively hyperbolic groups.  In Subsection
\ref{para;proof_hyp} we prove Proposition \ref{prop;Omega} in the case
of hyperbolic groups.  In Subsection \ref{ss:proofofmainprop}, we return
to the setting of toral relatively hyperbolic groups, and prove Theorem 
\ref{theo;algo_comp} (assuming Proposition \ref{prop;Omega}).

For the reader interested only in hyperbolic groups, this section contains a
complete proof of Theorem \ref{theo;algo_comp}, which is the only result needed
from this section for Sections \ref{sect:splittings}--\ref{ProofMainThm} (which
solve the isomorphism problem).  Such a reader should keep in mind that
torsion-free hyperbolic groups are, in particular, toral relatively hyperbolic.

    \subsectionind{Notations and Objective}\label{part;not}

      \subsubsection{Basic notations}
        In the following, $H_2$ is a non-elementary toral relatively hyperbolic
        group, and $H_1$ is a finitely presented.

        The application we are aiming at is the case where $H_1$ is also a
        non-elementary toral relatively hyperbolic group, but, for now,  
        all the geometry  takes place in $H_2$. 
        The case where $H_1$ is abelian makes Theorem \ref{theo;algo_comp} a
        matter of some simple linear algebra (which we leave as an exercise for the 
        reader).  Therefore,
        we suppose henceforth that $H_1$ is non-abelian.

        We consider a coned-off Cayley graph $\widehat{Cay} (H_2)$ (which is just $Cay (H_2)$ if $H_2$ is hyperbolic), 
        with distance denoted by $d$.

        The graph $\widehat{Cay}(H_2)$ is $\delta$-hyperbolic.  By Theorem \ref{theo;finddelta} it is possible 
        to algorithmically find $\delta$ from a finite presentation of $H_2$, 
        hence the constants  $L$, $L_1$, $L_2$, $L_1'$ and $L_2'$ and $r$ of the previous section.

        For the group $H_2$, we denote by $F$ the free product, 
        and  $\mathcal{L}_0\subset \mathcal{L}\subset F$,   the
        normalized rational languages over $F$, as defined in
        Paragraph \ref{para;decision}. By the remark above, the language $\mathcal{L}$ (or more precisely the underlying automaton) can be explicitly computed.

         Let $\epsilon$ be as in Proposition \ref{prop;stab_con}, 
         for $(L'_1,L'_2)$-quasi-geodesics 
         without $r$-detours in  $\widehat{Cay}( H_2)$.   
         In the case of hyperbolic groups, this is just a fellow-traveling 
         constant for   $(L'_1,L'_2)$-quasi-geodesics (in either case the constant $\epsilon$ can be explicitly computed).

        Recall the natural projection $\nu:F\to H_2$, and that a preimage of an element of $H_2$ is a representative in $F$.  
        Let us denote by $a_i, i\in I$ the given symmetric generating set of $H_1$.

\begin{defn}      
        For any homomorphism $\phi:H_1 \to H_2$, an \emph{acceptable lift}  of $\phi$ 
        is a choice of representatives in $\mathcal{L} \subset F$,
        of  the elements $\phi(a_i)$ and $\phi(a_i a_j)$, $i,j\in I$.

\end{defn}

\begin{notation}
$\mathcal{B}_{H_1}(2)$ is the ball of radius $2$ in $H_1$. If $\phi : H_1 \to H_2$ is a homomorphism, then
$\tilde{\phi} : \mathcal{B}_{H_1}(2) \to F$ denotes an acceptable lift.
\end{notation}

The following is clear.   
\begin{lem} \label{lem;finlifts}
Each homomorphism $\phi : H_1 \to H_2$ has at least one, but only finitely many,
acceptable lifts, and one can compute the list of them.
\end{lem}
        
      \subsubsection{Short and long morphisms, the property $\Omega$} \label{Long_and_short}
        Let us consider $a,b \in H_1$, in the given generating set, such that
        $a$ and $b$ do not generate an elementary subgroup (or equivalently
        such that $[a,b]\neq 1$).
        We define below (Remark \ref{rem;Omega}) a property 
	$\Omega$ that an acceptable lift of a homomorphism $\phi$ may or may not have. 
        This is based on the following remark (we think of $H_2$ as hyperbolic for this explanation). 
        If $h\in H_2$ is very long, then $h\phi(a)h^{-1}$ has the property
        that paths representing it begin and terminate by segments close to a
        prefix of $h$, except, crucially,  if $h$ commutes  with  $\phi(a)$,
        and in this case, $h\phi(ab)h^{-1}$ has this property.

        If an acceptable lift $\tilde{\phi}$ is so that  $\tilde{\phi}(a)$ and
        $\tilde{\phi}(b)$ (or $\tilde{\phi}(ab)$ and $\tilde{\phi}(b)$, or
        interchanging $a$ and $b$) have almost same initial, and final
        subsegment of large fixed length, we say that $\tilde{\phi}$ does
        \emph{not} satisfy $\Omega$ (see Remark \ref{rem;Omega} for precise definition). 
	We think of $\phi$ as \emph{long} if all
        its acceptable lifts are like that, and as \emph{short} if none are
        like that ({\it i.e. } if all its acceptable lifts satisfy $\Omega$). 
        Unfortunately, this leaves some place for homomorphisms that are neither long nor short, 
	but this is not a serious difficulty.
        
        Indeed, Proposition \ref{prop;Omega} is sufficient for our needs. 
        It states that, in each conjugacy class of homomorphism, only finitely many of them are not long, 
        and at least one of them is short.
        
        There are other simpler properties that guarantee this dichotomy: one could choose to look for 
	morphisms in minimal position 
        in their conjugacy class (that is, minimizing $\max\{ d(1,\phi(a_i)) \}$). 
        But we ultimately want to use it in a solvable system of equations and inequations.
        The property $\Omega$ is interesting for us, because it is defined using normalized rational languages, hence it will enable 
        us to encode in a system of equations, inequations and constraints our search of \emph{short} morphisms.

        \begin{prop} \label{prop;Omega}
          Let $H_1$ be a finitely presented non-abelian group, with $a$ and $b$ two elements such that $[a,b]\neq 1$  
          and $H_2$ be a 
          toral relatively hyperbolic group with primary peripheral structure. 
          Let $S_1 = \{a_i, i\in i \}$ be a finite symmetric  generating set of $H_1$.

          There exists a computable finite subset $\mathcal{QP}$ of $H_2$,
          and, for every $h\in \mathcal{QP}$,  
          a computable normalized rational 
          language $\mathcal{L}_h$, such that:

        \begin{itemize}
          \item Any compatible homomorphism $\phi:H_1 \to H_2$ has a conjugate $\psi$
            for which every acceptable lift $\tilde{\psi}: \mathcal{B}_{H_1}(2) \to F$ 
            satisfies the three conditions, for all $h\in \mathcal{QP}$ 
            \begin{enumerate}
              \item[] $\Omega(i)$ \hskip .2cm  either $\tilde{\psi}(a)$ or $\tilde{\psi}(b)$ is outside $\mathcal{L}_h$,
              \item[] $\Omega(ii)$ \hskip .2cm if $\tilde{\psi}(b)$ is in $\mathcal{L}_h$
                then either $\tilde{\psi}(ab)$ or 
                $\tilde{\psi}(a^{-1}b)$ is not,
              \item[] $\Omega(iii)$ \hskip .1cm if $\tilde{\psi}(a)$ is in $\mathcal{L}_h$ then
              either $\tilde{\psi}(ba)$ or $\tilde{\psi}(b^{-1}a)$ is not. 
            \end{enumerate}

          \item Any compatible homomorphism $\phi : H_1 \to H_2$ has 
            only finitely many  conjugates   $\psi$ so that $\psi$ has an
            acceptable lift $\tilde{\psi} : \mathcal{B}_{H_1}(2) \to F$  
            satisfying  $\Omega(i) \wedge \Omega(ii) \wedge \Omega(iii)$ for all $h\in \mathcal{QP}$. 
        \end{itemize}  

        \end{prop}

        \begin{remarknum}\label{rem;Omega}
          The property $\Omega = ( \forall h\in \mathcal{QP}, \, \Omega(i)\wedge \Omega(ii)\wedge \Omega(iii))$
          is clearly a boolean combination of properties of  membership to
          the $\mathcal{L}_h$ or their complements for the representatives
          $\tilde{\psi}(a),\tilde{\psi}(ab),\tilde{\psi}(a^{-1}b),\tilde{\psi}(b),\tilde{\psi}(ba)$
          and, $\tilde{\psi}(b^{-1}a) $. 

          Hence, it is a boolean combination of \emph{normalized rational constraints} on
          these representatives.
        \end{remarknum}

        We think of $\tilde{\psi}$ satisfying $\Omega$ as a indication of $\psi$ being "short" (as justified by the study to come).
        
        For convenience, let us precise that 
        $\tilde{\psi}$ does not satisfy $\Omega$ if there is $h\in
        \mathcal{QP}$ so that either (i) both $\tilde{\psi}(a)$ and
        $\tilde{\psi}(b)$ are in $\mathcal{L}_h$; (ii) all three of
        $\tilde{\psi}(b)$,  $\tilde{\psi}(ab)$ and $\tilde{\psi}(a^{-1} b)$ are
        in ${\mathcal L}_h$; or (iii) all three of $\tilde{\psi}(a)$, $\tilde{\psi}(ba)$ and
        $\tilde{\psi}(b^{-1}a)$ are in ${\mathcal L}_h$.

        In the next section, we concentrate on the case of hyperbolic groups. We briefly comment its content.

        In Paragraph \ref{para;QP}, 
        we explain  {\em quasi-prefixes}. 
        In the hyperbolic case, 
        $\mathcal{QP} = \{h\in H_2, \,  d(1,h) = 8\epsilon + 80\delta \} $.  For each such element $h$, 
        the language $\mathcal{L}_h$, is the set of elements $\tilde{g}$ of $\mathcal{L}\subset F $ labeling a path passing within 
        a certain neighborhood of $h$ and of $\nu(\tilde{g})h$ 
        (we interpret this as $h$ being a quasi-prefix and a quasi-suffix of $\tilde{g}$). We then   prove
        that these languages are normalized rational.

        Paragraph \ref{part;finiteness} serves to prove the second 
        point of the proposition:  
        in each conjugacy class of compatible homomorphism, 
                    only finitely many homomorphisms have an acceptable 
                    lift satisfying $\Omega$ 
                    (this is Corollary \ref{coro;fini}).  
        The reason is that if $h$ is very large, and $h_0\in [1,h]\cap
        \mathcal{QP}$, then  $h\phi(a)h^{-1}$  and $h\phi(b)h^{-1}$  usually 
        have $h_0$ as quasi-prefix and quasi-suffix (and hence fail to satisfy
        $\Omega$). Unfortunately this might not be the case if $h^{-1}$ is close to
        the centralizer of $ \phi(a)$.  We show in Corollary \ref{coro;fini} that if this happens,  then $h^{-1}$ is far
        from the centralizer of $\phi(b)$,  of $\phi(ab)$ and of $\phi(a^{-1}b)$, and this makes 
        every acceptable lift of $\phi$ fail to satisfy $\Omega$.

On the other hand, Paragraph  \ref{part;exists} serves to 
        prove the first point of the proposition, embodied as 
        Corollary \ref{coro;exists}. 
        The tactic here is simple: if a monomorphism has an acceptable lift
        contradicting $\Omega$,  it does not minimizes  $\max\{ d(1,\phi(a)),d(1,\phi(b))  \}$. 

        In fact, the same tactics are used for the relative case, but some extra
        difficulties appear.  We discuss these extra difficulties later, in Section 
        \ref{part;proof_gene}.

      \subsection{Proof of Proposition \ref{prop;Omega} 
        in the case of hyperbolic groups} \label{para;proof_hyp}

          In this subsection,  $H_2$ is a torsion-free
          hyperbolic group, with a word metric $d$.

          \subsubsection{Quasi-prefixes, $\mathcal{QP}$ and $\mathcal{L}_h$
            for $h\in H_2$.} \label{para;QP} 
        Recall that $(L'_1,L'_2)$-quasi-geodesics are $\epsilon$-close to a geodesic.
        We choose constants
        $\eta = \epsilon +
        10\delta$,  $\rho = 8\eta$ and $\eta' = \eta + \epsilon$. 
        We will emphasize later why we choose these values.

\begin{defn} \label{def:QP}
        Let $\mathcal{QP}= \{ h \in H_2, d(1,h) =\rho \}$. 
\end{defn}
Clearly, $\mathcal{QP}$ is a finite set, which can be computed from a solution
to the word problem and a knowledge of $\delta$ (both of which can be obtained
from a finite presentation using the  algorithm of Papasoglu mentioned in Subsection 
\ref{ss:uniform} above).

The following definition introduces the key idea in the proof of Proposition \ref{prop;Omega}.
\begin{defn} \label{def:quasi-prefix}
        Let $h$ be in  $\mathcal{QP}\subset H_2$. 
        Given  
        $\tilde{g} \in \mathcal{L}  \subset 
        F$, we say that $h$ is a {\em quasi-prefix} of $\tilde{g}$ if  
        $h$ is at distance at most $\eta$ from the path labeled by    
        $\tilde{g}$ in  $Cay (H_2)$.  
\end{defn}        

\begin{defn} \label{def:Lh}
Given $h \in \mathcal{QP}$, 
the language $\mathcal{L}_{h} \subset
         \mathcal{L}$ is the set of 
         elements $\tilde{g}$ of $\mathcal{L}$ whose normal form
        is of length at least $2\rho L'_1 + L'_2$ and so that $\tilde{g}$ and 
        $\tilde{g}^{-1}$ both have $h$ as a quasi-prefix. 
\end{defn}

Given these definitions, the condition $\Omega$ can be compared to the ``forbidden states'' in Sela's solution \cite{SelaIso}. 
However, as we have said, our use of rational constraints allows for considerable
streamlining in the solution of the isomorphism problem in later sections.

        \begin{lem} \label{lem;Lh_normalized}
          For all $h\in \mathcal{QP}$, the language $\mathcal{L}_h \subset F$ 
          is normalized rational, and computable.
        \end{lem}
\begin{proof}
           There are finitely many elements at distance $\eta$ from $h$, and
           for each of them, say $h'$, the  elements in $F$  
           that label  $(L'_1,L'_2)$-quasigeodesics in $Cay (H_2)$ from $1$
           to $h'$  all have length at most $L'_1 \rho +
           L'_2$.  
           It is therefore possible to compute the finite 
           list $\mathfrak{L}_h$ of all the normal forms of these  elements. 
           The language $\mathcal{L}'_h$ of the elements of $F$ that have a
           normal form with both a prefix and a suffix in $\mathfrak{L}_h$, is
           easily seen to be a normalized rational language. 
           But  $\mathcal{L}_h$ is the subset of $\mathcal{L} \cap \mathcal{L}'_h$ of words of length at least $2\rho L_1' + L_2'$.  
 \qed \end{proof}

    \subsubsection{A finiteness result}\label{part;finiteness}
In the next lemma, we will use the fact that $\eta = \epsilon + 10\delta$ (in fact we use $\geq$). 
Given an element $g\in H_2$, we denote by $Cent(g)$ its centralizer in $H_2$.

         \begin{lem}\label{lem;prefix} 
            For all $g\in H_2$, 
           there is a constant $K_{pre}(g)$ such that if the coset $hCent(g)$
           is at word distance at least  $K_{pre}(g)$ from $1$, 
           and if $h_0\in \mathcal{QP}$ is on a geodesic $[1,h]$,  
           then any representative 
           of $hgh^{-1}$ in $\mathcal{L}\subset F$ is in 
           $\mathcal{L}_{h_0}$.  
         \end{lem}

  \begin{proof} 
        If $hgh^{-1}=h'gh'^{-1}$, then $h'^{-1}h \in Cent(g)$ and so
        $hCent(g) = h'Cent(g)$.  Therefore, given $K(g)$, there is a constant 
        $K_{pre}(g)$ such that if $d(1,hCent(g)) \ge K_{pre}(g)$
       then $d(1,hgh^{-1}) > K(g)$. 
       
       Let $K(g) = 2\rho + 30\delta + d(1,g)$. Let $h_0$ be as in the statement. It suffices to prove that, if $d(1,hgh^{-1}) >K(g)$, then $h_0$ is at distance at most $\eta-\epsilon = 10\delta$ from any segment $[1,hgh^{-1}]$.
 We prove the contrapositive:
           if  $h_0$ is at distance at least $10\delta$ 
           from a segment $[1,hgh^{-1}]$,   
        then   $d(1,hgh^{-1})\leq  K(g)$
        (and similarly if $hgh^{-1} h_0$ is at distance at least 
        $10\delta$ from  $[1,hgh^{-1}]$).

           By hyperbolicity in the quadrilateral
           $(1,h,hg,hgh^{-1})$, the vertex $h_0$ is at distance at most 
           $10\delta$ from either $[h,hg]$ or $[hg,hgh^{-1}]$.

           First assume that there is $v\in [h,hg]$ at 
           distance $10\delta$
           from $h_0$. Then we bound the distances, by triangular inequalities:
           \[d(1,h) \leq
           d(1,h_0)+d(h_0,v) + d(v,h)\] and also 
           \[d(1,hgh^{-1}) \leq d(1,h_0) +
           d(h_0,v) + d(v,hg)+d(hg,hgh^{-1}).\]
           Using $d(1,h_0) = \rho$ and  $d(hg,hgh^{-1}) =d(1,h)$, 
           one obtains
\begin{eqnarray*}
d(1,hgh^{-1}) & \leq& 2\rho +20\delta + d(v,hg)+d(v,h)\\
& = & 2\rho  +20\delta +d(1,g) \\
& <& K(g).
\end{eqnarray*}

           Suppose now that there is $w \in [hg,hgh^{-1}]$ 
           at distance $10\delta$ from $h_0$. One has the bounds, again by
           triangular inequalities 
\begin{eqnarray*}
d(1,hgh^{-1}) &\leq& d(1,h_0) + d(h_0,w)+ d(w,hgh^{-1})\\
& \leq & \rho + 10\delta + d(w,hgh^{-1}),
\end{eqnarray*}
and also 
\[d(1,h) \leq d(1,h_0) + d(h_0,w)+d(w,hg) +d(hg,h),\]
which implies
\[d(1,h) \leq \rho + 10\delta +
           d(w,hg) + d(1,g). \]
           Since 
\begin{eqnarray*}
d(1,h) &=&d(hgh^{-1},hg)\\
&=&d(w,hgh^{-1}) + d(w,hg),
\end{eqnarray*}
one deduces 
\[d(w,hgh^{-1})\leq \rho + 10\delta  + d(1,g).\] 
Together with the first bound obtained, this gives 
\[d(1,hgh^{-1}) \leq  2\rho + 20\delta +  d(1,g) < K(g),\]
as required. 
\qed \end{proof}

         \begin{cor} \label{coro;fini}
           Given a compatible homomorphism $\phi: H_1 \to H_2$, 
           only finitely many conjugates of $\phi$ have an acceptable lift 
           satisfying $\Omega$.   
         \end{cor}
\begin{proof}
           Let $K_{pre}>\max_{x\in \mathcal{B}_{H_1}(2) \smallsetminus \{1\}} \{ K_{pre}(\phi(x))\}$. For $x\in H_1$ we write $N_x=\{h\in H_2, d(1,hCent\phi(x)) \leq K_{pre}\}$. This means that $h\in N_x$ if and only if $h^{-1}$ is in the $K_{pre}$-neighborhood of $Cent(\phi(x))$. Note that $N_x=N_{x^{-1}}$.

           Let us remark now that, by the compatibility of $\phi$,  $\langle\phi(a)\rangle$ and $\langle\phi(b)\rangle$ are infinite cyclic subgroups, not in the same elementary subgroup.  Therefore, the $K_{pre}$-neighborhood of their centralizers have finite intersection, and $N_a\cap N_b $ is finite. Similarly, since $[a,ab]\neq 1$, by compatibility of $\phi$, $\langle\phi(a)\rangle$ and  $\langle\phi(ab)\rangle$ are not in the same elementary subgroup, and for the same reason as above, $N_a\cap N_{ab} $ is finite (and similarly $N_a\cap N_{a^{-1}b}$).

In order to obtain a contradiction, assume that for different $h(n)$, the
           homomorphism  $h(n)\phi h(n)^{-1}$ (which we write
           $\phi^{h(n)}$ for readability)  has an acceptable lift satisfying $\Omega$ ({\it
           i.e. } is ``short"). By Lemma  \ref{lem;prefix} (and adopting the
           notation $h(n)_0$ from it),   for all $n$, $h(n)\in N_a \cup N_b$, 
           since otherwise, $\phi^{h(n)}(a),\phi^{h(n)}(b)$ have their
           representatives in $\mathcal{L}$ 
           in $\mathcal{L}_{h(n)_0}$, falsifying $\Omega(i)$.
           
We can assume that for all $n$, $h(n)\in N_a$.  Also, since $N_a\cap N_b$ is
finite,  we may assume that for all $n$, $h(n)\notin N_b$. Thus, all representatives 
in $\mathcal{L}$ of $\phi^{h(n)}(b^{\pm 1})$ are in $\mathcal{L}_{h(n)_0}$.  
Since $\Omega(ii)$ is satisfied for some acceptable lift, at least one representative 
in $\mathcal{L}$ of  $\phi^{h(n)}(ab)$ or of  $\phi^{h(n)}(a^{-1}b)$ is not in
$\mathcal{L}_{h(n)_0}$.
By Lemma   \ref{lem;prefix}, this implies that $h(n)\in N_{ab} \cup N_{a^{-1}b}$ for all
$n$. But we noted that 
$N_a\cap N_{ab}$ and $N_a\cap N_{a^{-1}b}$ are finite, thus contradicting the 
assumption that the $h(n)$ are all distinct.
\qed \end{proof}

     \subsubsection{An existence result}\label{part;exists}
          It is obvious that,  in every
          conjugacy class of homomorphism from $H_1$ to $H_2$, there is 
          a (possibly nonunique) homomorphism minimizing $\max\{d(1,\phi(a)),d(1,\phi(b))\}$. 
	We aim to prove that such a 
          morphism cannot contradict $\Omega$. We will find a conjugation decreasing this quantity for 
          every morphism contradicting $\Omega$.

     Let $\eta'=\eta+\epsilon$.

        \begin{lem} \label{lem;pour_les_geod1}
          Suppose that $h\in \mathcal{QP}$, that $\tilde{g} \in \mathcal{L}_h$, and that
          $g$ is the image of $\tilde{g}$ in $H_2$. Let $[1,g]$ be a geodesic 
          in $Cay (H_2)$. 

           Then there is a vertex  $w\in [1,g ]$ such that
          $d(w,h) \leq \eta'$.

        \end{lem}

\begin{proof}
          The path defined by $\tilde{g}$ in $Cay (H_2)$ contains a vertex 
           $v$ such that $d(h,v)\leq \eta$. 
          Since this
          path is an $(L'_1,L'_2)$-quasi-geodesic with end-points $1$ and $g$, 
          it is contained in the
          $\epsilon$-neighborhood of $[1,g]$. Thus there is $w\in [1,g]$ at
          distance at most $\eta$ from $v$. One obtains $d(w,h)\leq
          \eta+\epsilon  = \eta'$. 
\qed \end{proof}

For the next result, we use the fact that $\rho > 2\eta'$.

        \begin{lem} \label{lem;reduc_reguliere}
          Let $h\in\mathcal{QP}$. 
          If $g\in H_2$ is such that $d(h,[1,g])<\rho/2$ and 
          $d(gh,[1,g])<\rho/2$, then  $d(1,h^{-1}gh)< d(1,g)$.
          In particular, this is true if a representative $\tilde{g}$ of $g$ 
          is in 
          $\mathcal{L}_h \subset F$.

        \end{lem}

\begin{proof}
          If $w$ and $w'$ are vertices in $[1,g]$ at distance less than
          $\rho/2$ from $h$ and $gh$ respectively, one computes  $d(h,gh)
          < \rho + d(w,w') = \rho +d(1,g) -d(1,w) - d(w',g)$.  Since
          $d(1,w)\geq d(1,h) -\rho/2 = \rho/2  $, and similarly for
          $w'$, one gets $d(h,gh)  < d(1,g)$.  
          
Now assume that $\tilde{g} \in \mathcal{L}_h$.  
          By Lemma \ref{lem;pour_les_geod1}, there is a vertex $w\in [1,g]$
          with $d(w,h)\leq \eta'$, and similarly, there is $w'\in [1,g]$
          such that $d(w',gh)\leq \eta'$.  By choice of $\rho$, we have $ \eta' > \rho/2$, so 
the second assertion follows from the first.
\qed \end{proof}

For the next result, we use the fact that
$\rho  > 4\eta' + \delta = 8\epsilon + 41\delta$.

        \begin{lem}\label{lem;4cas}
       Let $\phi:H_1 \to H_2$ be a homomorphism with an acceptable lift $\tilde{\phi}:
        \mathcal{B}_{H_1}(2) \to \mathcal{L}$ not satisfying $\Omega$, and  
        let $h$ be as provided by the fact that $\Omega$ does not hold. 
   
           Then   $d(1,h^{-1} \phi  (a) h) < \max\{d((1,\phi (a)),
           d((1,\phi (b))\}  $, and similarly for $h^{-1} \phi
          (b) h$. 
          
          In particular, $\phi$ does not minimize the quantity $\max\{d((1,\psi (a)), d((1,\psi (b))\}$ over its 
	conjugacy class.

        \end{lem}

    \begin{proof}  
          If $\tilde{\phi}(a) \in \mathcal{L}_h$, by  Lemma
          \ref{lem;reduc_reguliere}, $d(1,h^{-1}\phi(a) h) <d(1,\phi (a))$.
          We assume now that  $\tilde{\phi}(a) \notin \mathcal{L}_h$. Since
          $\Omega$ is not satisfied by $\tilde{\phi}$, and by definition of
          $h$, necessarily $\Omega(ii)$ is false, and one has that
          $\tilde{\phi}(b)$, $\tilde{\phi}(ab)$ 
          and $\tilde{\phi}(a^{-1}b)$  are in $\mathcal{L}_h$.

          Most of the discussion will hold in the triangle 
          $(1,\phi(a), \phi(ab))$, and for one case, in $(1,\phi(a^{-1}), \phi(a^{-1}b))$. 
          The discussion will essentially 
          hold on the possible combinatorial configurations for the approximating trees for the vertices 
	$1,\phi(a),\phi(ab),h,\phi(a)h$.

         We choose geodesic segments for the sides of the triangle $(1,\phi(a), \phi(ab))$. 
         By Lemma \ref{lem;pour_les_geod1} there are vertices 
          $v$ and $v_a$, such that  
          $v\in [1,\phi(ab)]$, $d(v,h)\leq \eta'$,
          and 
          $v_a\in [\phi(a),\phi(ab)]$,   $d(v_a,\phi(a)h) \leq \eta'$.

          By hyperbolicity, there is $v' \in [1,\phi(a)] \cup [\phi(a),\phi(ab)]$ with   $d(v,v')\leq \delta$.   We distinguish four cases (illustrated in Figure 
          \ref{fig;4cas_bis}). 
          The first dichotomy
          (cases $(\alpha)$ and $(\beta)$) 
          concerns whether $v'\in  [\phi(a),\phi(ab)]$ or 
          $v'\in  [1,\phi(ab)]$.

          In case $(\alpha)$, both $v'$ and $v_a$ are in
          $[\phi(a),\phi(ab)]$, and we denote $(\alpha_1)$ 
          the case where they
          appear on the segment in order  $(\phi(a), v',v_a, \phi(ab))$, 
          and
          $(\alpha_2)$ when they appear in the order  
          $(\phi(a),v_a,v', \phi(ab))$.

          In case $(\beta)$,  $v'$ is on $[1,\phi(ab)]$, and we make a
          dichotomy on the position of $v_a$. Either there is $v_a'$ in  $
          [1,\phi(a)]$ within a distance at most $2\eta'$ of $v_a$, 
          (case  $(\beta_2)$), or there is not (case $(\beta_1)$).

              \begin{figure}[hbt]
                      \begin{center}
              \includegraphics[width=5in]{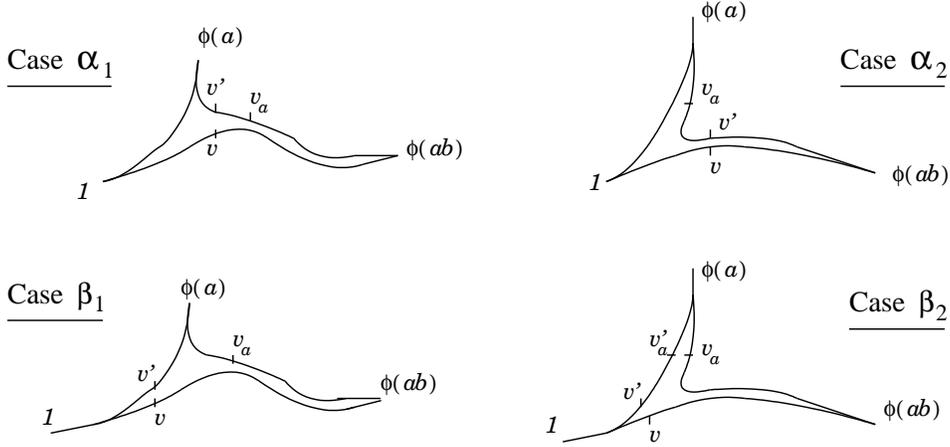}
              \caption{The four cases of Lemma \ref{lem;4cas} 
              (and the second set of four cases in Lemma \ref{lem;4cas_rh})}
              \label{fig;4cas_bis}
            \end{center}
          \end{figure}

          We now treat case $(\alpha)$. One has 
          \begin{eqnarray*}
          d(h,\phi(a)h) &\leq & d(h,v) + d(v,v')+ d(v',v_a) +d(v_a,\phi(a)h) \\
          & \leq & 2\eta' +  \delta + d(v',v_a).
          \end{eqnarray*} 
        We want to prove that $d(h,\phi(a)h) <d(\phi(a),\phi(ab))=d(1,\phi(b))$, 
        and for that it is sufficient to prove that
        \[d(v',v_a) < d(\phi(a),\phi(ab))-2\eta'-\delta.\]

          In case $(\alpha_1)$, on the segment  
          $[\phi(a),\phi(ab)]$  we
          have  
          $d(v',v_a)=  d(\phi(a),\phi(ab)) -  d(\phi(a), v')) -
          d(\phi(ab), v_a)$. Since  $(\phi(b))^{-1}$ has also $h$ as a quasi-prefix,  and by the length requirement on elements of $\mathcal{L}_h$ in Definition \ref{def:Lh}, we have
        $d(v_a,\phi(ab)) > 
        \rho-\eta-\epsilon \geq 2\eta'+\delta$, and we get the required conclusion.

          In case $(\alpha_2)$,  
          $d(v_a,v')=  d(\phi(a),\phi(ab)) -  d(\phi(a), v_a))
          - d(\phi(ab), v' )$. Since   
          $d(\phi(a),v_a) \geq d(\phi(a),\phi(a)h)-\eta'= \rho-\eta'$, one has 
        $d(\phi(a),v_a)> 2\eta'+\delta$, 
        hence the result.

          In  case $(\beta_1)$,  since
          $v_a$ is at distance greater than $2\eta'$ 
          from $[1,\phi(a)]$, the distance from $\phi(a)h$ to     $[1,\phi(a)]$
          is  greater than $\eta'$,      
           and by Lemma \ref{lem;pour_les_geod1}, no representative in 
          $\mathcal{L}$ of $\phi(a^{-1})$  has 
          $h$ as quasi prefix. Therefore, 
          $h$ lies at distance at least
          $\eta'$ from $[1,\phi(a^{-1})]$.         
             Thus, $\phi(a^{-1})$ falls into case $(\alpha)$, in the triangle $(1, \phi(a^{-1}), \phi(a^{-1}b))$
          and in this case we get that  
          $d(1,h^{-1}\phi (a) h )< d(1,\phi(b))$.

          Finally, in case $(\beta_2)$ both $v'$ and $v_a'$ are 
          on $[1,\phi(a)]$. Since $\rho/2 > 2\eta'$,  Lemma \ref{lem;reduc_reguliere}
          implies that $d(1, h^{-1} \phi (a) h) < d(1,\phi (a))$. 
          This finishes the proof in all the cases.
\qed \end{proof}

The next result now follows immediately.
        
        \begin{cor}\label{coro;exists}
          Any acceptable lift of a  compatible homomorphism which  minimizes $\max\{d(1,\phi(a)), d(1,\phi(b))\}$
          satisfies $\Omega$.
        \end{cor}
        
        Finally, we complete the proof of Proposition \ref{prop;Omega} in
        the hyperbolic case.
        
\begin{proof}[of Proposition \ref{prop;Omega}, hyperbolic case]
        As defined in Paragraph \ref{para;QP}, the set $\mathcal{QP}$ is
        computable, and for any given $h$, the language $\mathcal{L}_h$ is explicit,
        computable, 
        and normalized rational
        by Lemma \ref{lem;Lh_normalized}.  
        The first point is then ensured by Corollary \ref{coro;exists}, and
        the second point is ensured by Corollary \ref{coro;fini}. 
\qed \end{proof}

    \subsection{An algorithm for listing compatible monomorphisms}
    \label{ss:proofofmainprop}
      In this part, we return to the generality of toral relatively hyperbolic
      groups, but it will (hopefully) not be confusing for the reader interested only in
      hyperbolic groups, since very little specificity is to be considered. We
      consider Proposition \ref{prop;Omega} established in this setting 
      (see the proof in  the relative case in Section \ref{part;proof_gene}).

      Thus, $H_2$ is a non-elementary toral relatively 
      hyperbolic group, and $H_1$ is 
      finitely presented and non-abelian.  We assume that we are given a
     finite presentation, and a solution to the word problem in $H_1$. Let $d_0$ be
      the word metric in $H_1$ for this presentation.

          \subsubsection{Coding: the sentences $\Sigma_R$, and  $\Lambda_{R,\ell}$} \label{para;code}
        Assume we are given $R\geq 8$ and $\ell$ a finite list of compatible morphisms $H_1\to H_2$.  These homomorphisms are encoded as elements of $H_2$ which
        represent the images of the (fixed) generating set of $H_1$.
      
       For any given $R\geq 8$, let $B_R$ be the set of non-trivial 
elements of length at most $R$ in the word metric $d_0$. 
      
      Moreover, let $x_{1}^{(1)}, x_2^{(1)}, \dots, x_{c_1}^{(1)}$ be given generators for representatives of the maximal cyclic subgroups of $\mathcal{P}_1$, the primary peripheral structure of $H_1$. 
	Also, let  
	$\beta_1^{(1)}, \dots, \beta_{a_1}^{(1)}$ be the given basis for representatives of the 
	other abelian subgroups of $\mathcal{P}_1$. 
      Let now $S^0_R = B_R \cup (\bigcup_j \beta_j^{(1)} ) \cup \{x_{1}^{(1)},\dots, x_{c_1}^{(1)} \}$. 
      
      Since we are given a solution to the word problem in $H_1$, we can explicitly compute $B_R$, hence $S^0_R$, 
      and a presentation of $H_1$ for the generating set $S^0_R$. With a simple triangulation argument, and possibly by enlarging 
      $S^0_R$ into a larger computable finite set $S_R$, we can assume that  all relations in this presentation are of length $2$ or $3$.

        Let $\{u(s),  s\in S_R \}$ a set of unknowns that we will use to define systems of equations 
        (formally speaking, this is just a set, another copy of $S_R$). Also,
       let $\underline{g}=\{g_1, \dots, g_{c_1+a_1}\}$ be another set of
        unknowns.

        We now describe the set of parameters of the system of equations we
        want to write. Let $x_{1}^{(2)}, x_2^{(2)}, \dots, x_{c_2}^{(2)}$ be
        given generators for representatives of the maximal cyclic subgroups
        of $\mathcal{P}_2$, the primary peripheral structure of $H_2$. Let
        also  $\beta_1^{(2)}, \dots, \beta_{a_1}^{(2)}$ be the given bases for
        representatives of the other abelian subgroups of $\mathcal{P}_2$.

        Let us now define the system $\Sigma$ of equations in $H_2$, whose
        unknowns are the elements of $S_R$ and parameters the elements  of
        $H_2$ above.

  \[
        \begin{array}{ll}
        \Sigma_R(S_R , \underline{g})        
          = & 
   \left( \displaystyle\bigwedge_{s_i s_j =1} u(s_i)u(s_j) = 1 \right) 
                 \wedge    \left( \displaystyle\bigwedge_{s_i s_j s_k=1}
                             u(s_i)u(s_j)u(s_k) =1  \right)  \\
&      \wedge   
\left( \displaystyle\bigwedge_{i=1}^{c_1}  \displaystyle\bigvee_{j=1}^{c_2} 
 (u(x^{(1)}_i))^{g_i}  =x^{(2)}_j     
 \right)  \\
&  \wedge 
\left( \displaystyle\bigwedge_{i=1}^{a_1}   
  \displaystyle\bigvee_{j=1}^{a_2}   
                \,    \displaystyle\bigwedge_{x\in \beta^{(1)}_i, y\in
                  \beta^{(2)}_j}  [ (u(x))^{  g_{c_1+i}} , \, y ]=1  \right)

       \end{array}
 \]
(where for readability we wrote $(u(x^{(1)}_i))^{g_i}$ for $g_i u(x^{(1)}_i) g_i^{-1}$).

Let $\tilde{\ell}$ be the finite and computable (Lemma \ref{lem;finlifts})
 set of acceptable 
lifts of morphisms in the list $\ell$. 
We now define the constraints on some representatives $\widetilde{u(s)}$ 
in $F$ of the unknowns $u(s)$.

\[ \begin{array}{ll}
 \Lambda_{R,\ell}(\widetilde{S_R})  
 = & \left(  \displaystyle\bigwedge_{\psi \in \tilde{\ell}}
 \displaystyle\bigvee_{s\in S_R} \widetilde{u(s)} \in \mathcal{L} \setminus \{\psi
 (s)\}  \right)  
  \wedge \left(   \displaystyle\bigwedge_{s \in B_R } 
          \widetilde{u(s)}   \in \mathcal{L} \setminus \mathcal{L}_0   \right)
 \\

\end{array}
\]

\begin{lem}
    If there is a compatible homomorphism $\phi$
    from $H_1$ to $H_2$ injective on the ball of radius $R$ of $H_1$ (for the given word metric $d_0$), not in
    $\ell$, then there is a solution to the system $\Sigma_R(S_R,\underline{g})$ such that for all $s\in S_R$, $u(s)=\phi(s)$. Conversely, if   there is a solution of the system $\Sigma_R(S_R,\underline{g})$ with constraints 
    $\Lambda_{R, \ell}( \widetilde{S_R} )$, then $s \mapsto u(s)$ can be extended in  a
    compatible homomorphism
    from $H_1$ to $H_2$ injective on the ball of radius $R$ (for $d_0$) of $H_1$,  not in
    $\ell$. 
\end{lem}
\begin{proof}
For any homomorphism $\phi$, the first two
         blocks of $\Sigma_R(S_R, \underline{g})$ are satisfied by the elements $u(s)= \phi(s)$ for each $s\in S_R$. 
          The third and fourth blocks, are satisfied whenever the homomorphism $\phi$ is compatible. The first block of  
        $\Lambda_{R,\ell}(\underline{\tilde{u}}, \underline{\tilde{v}})$ is satisfied by any acceptable lift of a homomorphism which is not
        in $\ell$, and the second block is satisfied by any acceptable lift of
        homomorphism injective on the $R$-ball of $H_1$. 
        Thus, one has the first part of the lemma.

        Conversely, if the unknowns satisfy $\Sigma(S_R,\underline{g})$, then every relation of
        length $2$ and $3$ (hence every relation of our triangular presentation of $H_1$)
        is satisfied in $H_2$ by the unknowns. Therefore, the map $s \to u(s)$
         extends to a homomorphism.
  
        By the fourth block of $\Sigma$, this homomorphism sends 
        the peripheral cyclic groups of $H_1$ to conjugates of peripheral
        cyclic groups of $H_2$, and by the fifth block,  the images of non
        cyclic peripheral subgroups  commute with  some  
        peripheral non-cyclic subgroups of $H_2$. 
        This makes the homomorphism compatible.

        The constraints in $\Lambda$,  
        are respectively the fact that $\phi
        \notin \ell$, and the injectivity on the $R$-ball of $H_1$. 
\qed \end{proof}

\begin{cor} \label{coro;decision}
  Let $H_2$ be a toral relatively hyperbolic group, and $H_1$ a finitely
  presented group, with a solution to its word problem. Let $a,b \in H_1$ be
  non-commuting elements in the given generating set.
  
  There is an explicit algorithm whose input is $R, \ell$, that always
  terminates, and whose output is
  ``yes'' if there is a compatible homomorphism  $\phi$, not in $\ell$, injective on the
  $R$-ball,  so that every acceptable lift $\tilde{\phi}$ of $\phi$
  satisfies 
  $\Omega$  for the elements $a,b$, and outputs ``yes" only if there
  exists such a homomorphism $\phi$, which has an acceptable lift  $\tilde{\phi}$
        satisfying  $\Omega$ for the same elements.

\end{cor}

\begin{proof}
        First, as noticed in the previous paragraphs,  
        we can compute $S_R$, a triangular presentation associated, and the systems above.
        Recall (Remark \ref{rem;Omega}) that $\Omega$ is a boolean combination
        of rational constraints on representatives of certain $u(s), s\in
        S_R$. By elementary boolean operations, the sentence 
        $\Sigma_R(S_R,\underline{g}) \wedge \Lambda_{R,\ell}(\widetilde{S_R})
        \wedge \Omega$ is equivalent to a disjunction of finitely 
        many systems of equations, inequations, and constraints, say
        $\Sigma'(R,\ell,1) \vee \dots \vee \Sigma'(R,\ell,m)$ (for ease of
        notation, we drop the variables for the remainder of the proof).
        
        For each of them, $\Sigma'(R,\ell,k)$, there is an explicit 
        algorithm $\mathcal{A}(k)$, given by Theorem
        \ref{theo;decision}, that always terminates and says ``yes'' 
        if a homomorphism has all its acceptable lifts satisfying 
        $\Sigma'(R,\ell,k)$, hence  
        $\Sigma'(R,\ell,1) \vee \dots \vee \Sigma'(R,\ell,m)$, that is $   
        \Sigma_R \wedge  \Lambda_{R,\ell} \wedge \Omega$.

        On the other hand, if, for some $k$,  $\mathcal{A}(k)$ says ``yes'', 
        then there is a homomorphism with an acceptable lift satisfying 
        $\Sigma'(R,\ell,k)$, hence 
        $\Sigma_R \wedge  \Lambda_{R,\ell} \wedge \Omega$. 

        An algorithm satisfying the requirements of the corollary is then as follows:
        run all the
        algorithms $\mathcal{A}(k)$, $k=1,\dots, m$.  The algorithms will
        eventually terminate. We define the final answer to be ``yes'' if and only if one 
        of the $\mathcal{A}(k)$  says ``yes''.
\qed \end{proof}

     \subsubsectionind{The main result}
\begin{proof}[of Theorem \ref{theo;algo_comp}]        
We now describe an algorithm that fulfills the conditions of 
Theorem \ref{theo;algo_comp}.  \footnote{Recall that we are excluding the case
that $H_1$ is abelian, since this case may be easily proved using linear
algebra.}

        Start with the empty list $\ell=\emptyset$, and $R=8$.    

First, run the algorithm from Corollary \ref{coro;decision}.  If this
algorithm outputs the answer ``yes", then search for a solution to
the equations $\Omega \wedge 
\Sigma_R(S_R,\underline{g}) \wedge \Lambda_{R,\ell} 
(\widetilde{S_R})$. 
Since there is a solution, we will eventually find one. 
 Such a solution $(S_R,\underline{g})$ corresponds to a
compatible homomorphism $\psi_R : H_1 \to H_2$ defined by $\psi_R(s) = u(s)$, that is injective on the $R$-ball of $H_1$ and has at least one acceptable lift satisfying $\Omega$.
Add $\psi_R$ to the list $\ell$, replace $R$ by $R+1$, and start again
the algorithm from Corollary \ref{coro;decision}.

If the algorithm from Corollary \ref{coro;decision} outputs ``no",
then we stop and the output is the list $\ell$.

        We need to check that, if, for some $R_0$,  
        there are finitely many conjugacy classes of compatible 
        homomorphisms injective on the ball of 
        radius $R_0$ of $H_1$, then the algorithm from Corollary \ref{coro;decision} 
	eventually outputs ``no" at some point, 
        and that, when it does, 
        the final output list contains at least one representative of every conjugacy class of 
        compatible monomorphism.

In order to obtain a contradiction, assume that at every step of the recursion above 
(that is
for every $R$),  the algorithm  outputs  ``yes" (and therefore, finds a
suitable $\psi_R$). There exist arbitrarily many compatible morphisms $\psi_R,
R\in \mathbb{N}$, all different,  with $\psi_R$ injective on the ball of
radius $R$ in $H_1$, and each has an acceptable lift satisfying
$\Omega$. By the second point of  Proposition  \ref{prop;Omega}, 
there must be infinitely many conjugacy classes of morphisms among the $\psi_R$. Hence, there is no finite subset $A$ of $H_1$ on which only finitely many non-conjugate compatible morphisms $H_1\to H_2$ are injective, and this proves the first assertion.

        Now assume that the algorithm from Corollary \ref{coro;decision} outputs
        ``no" on the input $R,\ell$. This means that for every homomorphism satisfying $\Sigma_R \wedge  \Lambda_{R,\ell}$, 
        there is an 
        acceptable 
        lift that fails to satisfy $\Omega$.  
        Consider a compatible monomorphism, we want to show that one of its conjugates is in $\ell$.
        By the first point of Proposition \ref{prop;Omega},   
        it has a conjugate such that every acceptable lift satisfies $\Omega$. This means that this conjugate is not solution of 
        $\Sigma_R \wedge \Lambda_{R,\ell}$. 
        Since it is a compatible injective homomorphism, 
        the only possibility for it failing to satisfy $\Sigma_R \wedge \Lambda_{R,\ell}$ is that it is already in the list $\ell$.  This is exactly what we require, and the proof is complete.
\qed \end{proof}

\section{Splittings} \label{sect:splittings}

In this section we investigate what happens when the algorithm
from Theorem \ref{theo;algo_comp} does not terminate. 

\begin{hypremark} \label{rem;BPR} Our approach can be exemplified briefly by the case of torsion free hyperbolic groups.
Recall Bestvina-Paulin's theorem (and variations on it) for a torsion free
hyperbolic group $H_1$ with an essential peripheral structure. 
If $H_1$  has infinitely many non-conjugate compatible automorphisms,  (and
even weaker: if for some hyperbolic groups $H_2$ with essential peripheral
structure, and for any finite subset $A\subset H_1$, there are  infinitely many
non-conjugate compatible homomorphisms $H_1\to H_2$ which are injective on $A$)   then $H_1$ admits a 
so-called compatible essential small action on a $\mathbb{R}$-tree.   
By deep results of Rips' theory of actions on $\mathbb{R}$-trees, this implies
that $H_1$ has a compatible essential splitting. 
Conversely, if $H_1$
has a compatible essential splitting, it has infinitely many Dehn twists.  By 
pre-composing with these Dehn twists we find infinitely many non-conjugate
compatible automorphisms.

Summing up our discussion, if the algorithm of Theorem \ref{theo;algo_comp}    does
not terminate for $H_1$ (and arbitrary hyperbolic group  $H_2$), then  $H_1$ has a compatible
essential splitting, and conversely,  if $H_1$ has a compatible essential
splitting, the algorithm of Theorem \ref{theo;algo_comp} applied to  $H_1=H_2$  does
not terminate. 
\end{hypremark}

In the next paragraphs, we generalize this to toral relatively hyperbolic
groups and primary splittings (Theorem \ref{H1Splits}). Most of the
important technical tools were developed elsewhere, 
see for instance  \cite{Groves_RH1}.

We then apply this analysis to find a maximal primary splitting for every
toral relatively hyperbolic group
(see Theorem \ref{MaxSplit}).  In the next section, we will turn this
into a primary JSJ decomposition of $\Gamma$.

\subsection{Limit $\R$-trees for toral relatively hyperbolic groups}
\label{section:Rtrees}
In this subsection we recall a construction of an action on an $\mathbb{R}$-tree, from 
\cite{Groves_RH1} (see also \cite{GrovesCAT(0)1}).

\begin{thm} \label{LimitTree}
\cite[Theorem 6.4]{Groves_RH1}
Suppose that $G$ is a finitely generated group, and that $\Gamma$
is a toral relatively hyperbolic group.  Suppose that 
$\{ h_i : G \to \Gamma \}$ is a sequence of pairwise non-conjugate
homomorphisms.  There is a subsequence $\{ f_i \}$ of
$\{ h_i \}$, an $\R$-tree $T$ and an isometric $G$-action on 
$T$ with no global fixed point which satisfies the following properties:
Let $K$ be the kernel of the $G$-action on $T$, and
let $L = G/K$.
\begin{enumerate}
\item \label{SegStabAb}
Stabilizers in $L$ of non-degenerate 
segments in $T$ are abelian;
\item \label{Tline}
If $T$ is isometric to a real line then for all but finitely many
$n$ the group $f_n(G)$ is free abelian;
\item \label{tripod}
If $g \in G$ stabilizes a tripod in $T$ pointwise then $g \in \text{ker}(f_n)$ for all but finitely many $n$;
\item \label{Stable}
Let $[y_1,y_2] \subset [y_3,y_4]$ be a pair of non-degenerate
segments of $T$ and assume that $\text{Stab}_L[y_3,y_4] \neq 
\emptyset$.  Then
\[      \text{Stab}_L[y_1,y_2] = \text{Stab}_L[y_3,y_4].        \]
In particular, the action of $L$ on $T$ is stable.
\item \label{seqStable}
If $g$ is not in the kernel of the $G$-action on $T$ then for
all but finitely many $n$ we have $g \not\in \text{ker}(f_n)$; and
\item \label{tf}
$L$ is torsion-free.
\end{enumerate}
\end{thm}

We comment briefly on the construction of the $\R$-tree $T$.
To the group $\Gamma$ is associated a space
$X$ (defined in \cite{Groves_RH1}):  Fix a finite generating set $S$ of $\Gamma$.  The space $X$ is constructed from the Cayley graph of $\Gamma$ with respect to $S$ by
`partially coning' and then adding Euclidean flats for each coset
of each parabolic (see \cite{Groves_RH1} for more details). The group
$\Gamma$ acts properly and cocompactly by isometries on $X$, and
thus to
a homomorphism $h : G \to \Gamma$ is associated a $G$-action on
$X$.  If we choose a basepoint $x \in X$, then we define
\[      \| h \| = \min_{\gamma \in \Gamma} \max_{g \in S} 
d_X(\gamma .x , h (g) \gamma . x).      \]

Given the sequence $\{ h_i \}$, a limiting action of $G$ is extracted 
on an asymptotic cone $X_\omega$ of $X$, which is formed using the scaling factors
$\left\{ \frac{1}{\| h_i \|} \right\}_{i=1}^{\infty}$ (see \cite{DS} for more details on asymptotic
cones in this context).  The asymptotic cone $X_\omega$ is a `tree-graded' space
(see \cite{DS}).  In the context of toral relatively hyperbolic groups
(rather than the more general setting of \cite{DS}), a $G$-action
on an $\R$-tree was defined in \cite{Groves_RH1}, which is the 
$T$ from Theorem \ref{LimitTree}.

\begin{cor} \label{MaxCycCor}
Suppose that $g \in L$ is so that for some $k$ the element
$g^k$ stabilizes some non-degenerate segment $[a,b]$ in $T$.
Then $g$ stabilizes $[a,b]$ (pointwise).
\end{cor}

\begin{proof}
The argument from the proof of Theorem 4.4(2) on page 1353 of
\cite{GrovesCAT(0)1} shows that there is no element $g \in L$ so
that $g^2$ stabilizes a segment in $T$ and $g$ acts by inversion
on this segment.  Thus we may assume that $k > 2$.

If $g^k$ fixes a non-degenerate segment $[a,b]$ (not fixed by $g^{k-1}$) then $g$ must fix a point
when acting on $T$.  Consider the subtree $T' \subset T$ consisting 
of $[a,b]$, together with an interval in $T$ from $[a,b]$ to 
$\text{Fix}(g)$.  Let $T'' = \langle g \rangle T'$.  Then $T''$ is fixed
pointwise by $g^k$.  

Since $k \ge 3$, we may assume that $T''$ contains a tripod. The element
$g^k$ must be in the kernel of the $G$-action, by (\ref{tripod}) and
(\ref{seqStable}) of Theorem \ref{LimitTree}.  In this case,
since $L$ is torsion-free, by Theorem \ref{LimitTree}.(\ref{tf}),
$g$ acts trivially on $T$, and so certainly fixes $[a,b]$.
\qed \end{proof}

We also need one other property of this limiting $\R$-tree,
which will be crucial in our ability to decide if, given a splitting
$\Lambda$ of a toral relatively hyperbolic group,
there is a splitting of a vertex group $V$ of $\Lambda$ which refines
the splitting $\Lambda$.

\begin{lem} \label{BoundEll}
Suppose that $G$ and $\Gamma$ are as in Theorem \ref{LimitTree},
that $\{ h_i : G \to \Gamma \}$ is a sequence of pairwise non-conjugate
homomorphisms, with associated limiting $\R$-tree $T$.

Suppose that $g \in G$ is such that there exists $D \ge 0$ so that
for all $i$ there exists $x_i \in X$ so that $d_X(h_i(g) . x_i , x_i)
\le D$.  Then $g$ fixes a point in the limiting action of $G$ on $T$.
\end{lem}
\begin{proof}
Suppose that the word length of $g$ with respect to $S$ is $n$.
Suppose that we have applied a conjugation to $h_i$ so that
$\| h_i \| = \max_{s \in S}d_X(h_i(s).x, x)$ for the chosen basepoint
$x$ of $X$.  Then $d_X(h_i(g).x,x) \le n \| h_i \|$.

By \cite[Theorem 4.12]{Groves_RH1}, the space $X$ has 
{\em relatively thin triangles}, which means that
there is some $\nu > 0$ so
that for any geodesic triangle $\Delta$ there is a flat $E \in X$ so that 
any side of $\Delta$ is contained in the $\nu$-neighborhood of the
union of $E$ and the other two sides.

Consider the quadrilateral with vertices $x, h_i(g).x, x_i$ and
$h_i(g) . x_i$ (and sides $[h_i(g) . x , h_i(g) . x_i ] = h_i(g) [x,x_i]$).
If all four vertices are contained in a flat $E$, then because the
stabilizer in $\Gamma$ of $E$ acts by translations, we have
$d_X(h_i(g).x,x) = d_X(h_i(g).x_i,x_i) \le D$.  In this case, let
$y_i = x$.

Otherwise, it is not hard to see, by checking the possible combinatorial configurations of the quadrilateral,  
	that there is a point $y_i$ within
$2n \| h_i \|$ of $x$ so that $d_X(h_i(g).y_i,y_i) \le D + 2\nu$.

In the limit, there is a point $y_\omega$ (a limit of a subsequence
of $\{ y_i \}$) so that $d_{X_\omega}(x_\omega,y_\omega) \le
2n$ and $g$ fixes $y_\omega$.  Note that it is required
that the points $y_i$ don't move away from $x$ too quickly, because
$X_\omega$ is the path component containing $x_\omega$, the
limit of the constant sequence $\{ x \}$.
\qed \end{proof}

Recall also the following:
\begin{lem} \label{AbMal}
\cite[Lemma 6.7]{Groves_RH1}
Each abelian subgroup of $L$ is contained in a unique maximal abelian 
subgroup. Maximal abelian subgroups of $L$ are malnormal.  
\end{lem}

\begin{thm} \label{SplitEss}
Suppose that $\Gamma$ is a toral relatively hyperbolic group, that
$H$ is a finitely presented group, and that $\{ \phi_i : H \to \Gamma \}$ is a sequence of
pairwise non-conjugate homomorphisms that converges into an
action of $H$ on a tree $T$.  Suppose furthermore that $H$ is
non-abelian, and that the action of $H$ on $T$ is faithful.

Then $H$ admits nontrivial essential splitting.
\end{thm}
\begin{proof}
We have found a faithful stable action of $H$ on an $\R$-tree
$T$.   By \cite[Theorem 9.5]{BFAccess}, $H$ admits a
nontrivial splitting
over a group of the form $E$-by-cyclic, where $E$ fixes a 
non-degenerate segment of $T$.

By Theorem \ref{LimitTree}.(\ref{SegStabAb}), stabilizers in $H$
of non-degenerate segments in $T$ are abelian.

Lemma \ref{AbMal} implies that maximal abelian subgroups of $H$ are
malnormal, which means that any abelian-by-cyclic subgroup of 
$H$ is abelian.  

It remains to check Condition (2) of Definition \ref{EssentialSplitting}, namely
that any edge group $E$ of the splitting is a direct factor of the (unique) maximal
cyclic group containing $E$.

The $\R$-tree
$T$ was found from the tree-graded metric space $\mathcal C_\infty$,
which in turn was found by taking a limit of homomorphisms
$\{ f_i : G \to \Gamma \}$.  

There are now two cases to consider: either $E$ stabilizes a nontrivial segment
in $T$ or else it arises from a surface or a toral piece in the band complex found
by the Rips Machine (the Rips Machine is explained in \cite{BFAccess}).
 
In case an edge group $E$ stabilizes a nontrivial segment in $T$,
the fact that Condition (2) is satisfied follows from Corollary \ref{MaxCycCor}.

Otherwise, $E$ arises either from a surface or a toral piece in
the band complex found by the Rips Machine.  In case of the surface, the
splitting corresponds to cutting the surface along a simple closed curve, which
represents a primitive element of the fundamental group of the surface.  
Let $T_0 \subset T$ be the subtree corresponding to the surface piece.  For a surface piece, if an element of $H$ fixes a point in $T_0$ then it 
fixes all of $T_0$ pointwise.  However, $T_0$ is not a line and tripod in $H$
stabilizers are trivial.  Therefore, in this case $E$ is cyclic.  We claim that $E$
is maximal cyclic in $H$.  If not, then there is an element $h \in H$ so that
$h \not\in E$ but $h^k \in E$.  Then $h$ permutes the $T_0, h. T_0, \ldots
h^{k-1} . T_0$, which are also disjoint.  Therefore $h$ fixes the (finite) tree consisting
of the shortest paths between the sub-trees $h^i . T_0$ and $h^j . T_0$ for
$i,j \in \{ 0, \ldots , k-1 \}$.  This implies that $h$ fixes a point $x \in T \smallsetminus
T_0$.  Clearly $h^k$ also fixes $x$.  But $h^k$ leaves invariant and acts hyperbolically on $T_0$.  It is impossible for an isometry of an $\R$-tree to leave invariant and
act hyperbolically on some subtree and also fix any point in the tree (since the 
minimal invariant subtree is unique).  This is a contradiction which shows that
$E$ is maximal cyclic. 

Now consider the case when $E$ corresponds to a toral piece, with
associated subtree $T_0$.  In this case $T_0$ is a line.  In this case,
$E$ is an extension of a group $E_0$ which fixes $T_0$ pointwise by
a free abelian group $E_1$ acting freely on $T_0$ by translations.  Corollary
\ref{MaxCycCor} shows that $E_0$ is a direct factor in the maximal abelian
subgroup containing it.  The above argument from the surface case shows
that $E_1$ is also a direct factor of the maximal abelian subgroup containing
it.  Since abelian subgroups of $H$ are malnormal, the extension of $E_0$
by $E_1$ is abelian, and is also a direct factor in the maximal abelian subgroup
containing it, as required.
\qed \end{proof}

%Given Lemma \ref{AbMal}, the proof of the following is identical
%to that of \cite[Lemma 2.1]{SelaDio1}

Let us recall also the following fact, that follows from  Lemma \ref{AbMal}.

\begin{lem} \label{MakeAbEll}\cite[Lemma 2.10]{Groves_MR}
Let $L = G/K$ be as in Theorem \ref{LimitTree},
let $M$ be a maximal abelian subgroup
of $L$, and let $A$ be any abelian subgroup of $L$.  Then
\begin{enumerate}
\item If $L = U \ast_A V$ then $M$ can be conjugated into
either $U$ or $V$.
\item If $L = U \ast_A$ then either (i) $M$ can be conjugated into
$U$; or (ii) $M$ can be conjugated to $M'$ and $L = U \ast_A M'$.
\end{enumerate}
\end{lem}

Using Lemma \ref{MakeAbEll}, if we have an abelian 
splitting $L = U \ast_A$ where $A$ is a subgroup of a 
non-elliptic maximal abelian subgroup $M$, we can convert
it into the amalgamated free product $L  = U \ast_A M$.

This will allow us to deduce the existence of compatible primary 
splittings in the sequel.

\subsection{Splittings} \label{SplitSection}
In this subsection, we investigate what happens when the algorithm
from Theorem \ref{theo;algo_comp} does not terminate.

It will be easy to deal with abelian groups, and we will not
have to apply the algorithm from Theorem \ref{theo;algo_comp} to
them.  Thus we will assume for this subsection that the groups
under consideration are non-abelian.

\begin{defn}
Suppose that $G$ is a group with a peripheral abelian structure.  A graph
of groups decomposition $\Lambda$ of $G$ is called {\em compatible} if for each
subgroup $A$ in the peripheral structure there is some vertex group $V$ in $\Lambda$
so that $A$ is conjugate into $V$.
\end{defn}

\begin{thm} \label{H1Splits}
Let $H_1$ and $H_2$ be non-abelian toral relatively hyperbolic groups
with abelian peripheral structures.  If the algorithm
from Theorem \ref{theo;algo_comp} does not terminate
then $H_1$ admits a nontrivial compatible primary splitting.
\end{thm}
\begin{proof}
If the algorithm from Theorem \ref{theo;algo_comp} does not
terminate, then there exists a sequence of compatible 
homomorphisms $\{ \phi_i : H_1 \to H_2 \}$ so that:
\begin{enumerate}
\item The $\phi_i$ are pairwise non-conjugate; and
\item For each $i$, the homomorphism $\phi_i$ is injective
on the ball of radius $i$ in $H_1$.
\end{enumerate}

The techniques of \cite{Groves_RH1} are designed to handle
exactly this situation.  By Theorem \ref{LimitTree}
we can extract an action of $H_1$ on an $\R$-tree $T$ from 
the sequence $\{ \phi_i \}$.  We will use this action to prove that there
exists a compatible splitting of $H_1$.

Denote by $\{ f_i \}$ the subsequence of $\{ \phi_i \}$ as in the 
statement of Theorem \ref{LimitTree}.

Since $H_1$ and $H_2$ are non-abelian, and each $\phi_i$ is injective
on the ball of radius $8$ about $1$ in $H_1$, the image $\phi_i(H_1)$
is non-abelian.  Therefore, by Theorem \ref{LimitTree}.(\ref{Tline}),
$T$ is not isometric to a real line.  Thus there is a tripod in $T$.

Suppose that $g \in H_1$ is in the kernel of the $H_1$-action on $T$.
Then by Theorem \ref{LimitTree}.(\ref{tripod}), $g \in \text{ker}(f_i)$
for all but finitely many $i$.  However, $\phi_i$ is injective on the ball
of radius $i$ about $1$ in $H_1$.  Therefore, if $g \neq 1$ then
$g$ is only in finitely many of the $\text{ker}(f_i)$ (since $\{ f_i \}$ is
a subsequence of $\{ \phi_i \}$).  This shows that
the kernel of the $H_1$-action on $T$ is trivial.

Theorem \ref{SplitEss} now implies that $H_1$ admits a 
nontrivial essential splitting, $\Lambda$.  

It needn't be the case that $\Lambda$ is a primary
splitting.  However, we claim that $H_1$ admits a nontrivial
compatible primary splitting.  

In order to simplify the situation, collapse all edge groups
of $\Lambda$ except one.  
Denote the resulting elementary
splitting of $H_1$ by $\Xi$. 

If there is a noncyclic abelian edge group that is not elliptic in
$\Xi$, then we apply Lemma \ref{MakeAbEll} to obtain a splitting $\Xi'$
where all noncyclic abelian subgroups are elliptic.  In particular,
all noncyclic groups in the peripheral structure of $H_1$ are
elliptic in $\Xi'$.  Note that $\Xi'$ is a primary splitting.

Consider a cyclic group $C = \langle c \rangle$ in the peripheral structure of $H_1$,
and let $E_1 = \langle e_1 \rangle , \ldots , E_j = \langle e_j \rangle$ 
be representatives of conjugacy classes of cyclic edge groups in
the peripheral structure of $H_2$.  

For all $i$, there exists $k$ so that $f_i(c)$ is conjugate to either $e_k$ or
$e_k^{-1}$.  In any case, it is clear that the translation length
of $c$ acting on the space associated to $H_2$ is bounded, which
implies (by Lemma \ref{BoundEll}) that $c$ fixes a point in the limiting
tree $T$.

The vertex groups of the splitting extracted from the action of $H_1$ 
on $T$ correspond to orbits of branching points with nontrivial
stabilizer (see \cite[Theorem 3.1, p.545]{SelaAcyl}).  Passing to the splitting $\Xi$ can only increase the 
vertex groups, so $c$ is conjugate into
a vertex group of $\Xi$, and so also into a vertex group of $\Xi'$.  Thus
the splitting $\Xi'$ is compatible, as required.
\qed \end{proof}

 The converse of  Theorem
\ref{H1Splits} might not be true for any group $H_2$, but it holds for $H_2= H_1$, as follows.

\begin{thm} \label{NoTermiffSplit}
Let $H_1$ be a non-abelian toral relatively hyperbolic group,
with abelian peripheral structure.  The algorithm from Theorem
\ref{theo;algo_comp} applied to compatible homomorphisms from
$H_1$ to itself terminates if and only if $H_1$ does not admit
a nontrivial compatible primary splitting.
\end{thm}
\begin{proof}
If the algorithm does not terminate then Theorem \ref{H1Splits} implies
that $H_1$ admits a nontrivial compatible primary abelian splitting.

Conversely, if $H_1$ admits a nontrivial compatible essential
primary splitting then by Lemma \ref{DehnInf} there is a Dehn twist or generalized
Dehn twist of infinite order in $\mbox{Out}(H_1)$, which
means that there are infinitely many conjugacy classes of compatible
monomorphisms from $H_1$ to itself, so by Theorem 
\ref{theo;algo_comp} the algorithm does not terminate.
\qed \end{proof}

\subsection{When $H_1$ and $H_2$ admit no compatible primary splittings}
We now solve a special case of the isomorphism problem.  This
case is already of substantial interest.  In particular, in the special
case that the groups are hyperbolic, and the peripheral structures
are empty, this result is the main result of \cite{SelaIso} (which
is \cite[Theorem 7.3, p.256]{SelaIso}).  This result will also be a key
ingredient in the final proof of Theorem \ref{Toral}. Also, this is sufficient for our proof of Theorem \ref{HomeoHyp} in Section \ref{HypSection}.

\begin{thm} \label{IsoNoSplit}
Let $H_1$ and $H_2$ be toral relatively hyperbolic groups
with abelian peripheral structures, and suppose that $H_1$
and $H_2$ do not  admit any nontrivial compatible primary
splittings (including compatible free product decompositions).  
Then it is decidable whether or not there is a 
compatible isomorphism between $H_1$ and $H_2$.

Moreover, in case there is a compatible isomorphism between
$H_1$ and $H_2$, the algorithm will provide a list consisting
of compatible isomorphisms which contains a representative of each
conjugacy class of compatible isomorphism between $H_1$
and $H_2$.
\end{thm}
\begin{proof}
By Theorem \ref{H1Splits}, the algorithm from Theorem 
\ref{theo;algo_comp} terminates when applied to compatible
homomorphisms from $H_1$ to $H_2$.  Therefore, we can
algorithmically find a finite list ${\mathcal C}_1$ of homomorphisms from $H_1$
to $H_2$ which contains a representative of every conjugacy class
of compatible monomorphism from $H_1$ to $H_2$ (and in particular 
if there is a compatible isomorphism, a conjugate of it will be in this list).

By applying the algorithm to compatible homomorphisms from
$H_2$ to $H_1$ (Theorem \ref{H1Splits} once again ensures
that this algorithm terminates), we find a finite list $\mathcal{C}_2$ of compatible homomorphisms
from $H_2$ to $H_1$ which contains a representative of each
conjugacy class of compatible monomorphism from $H_2$ to $H_1$.

Now, there is a compatible isomorphism from $H_1$ to $H_2$ if and
only if there is a compatible homomorphism $\phi_1 : H_1 \to H_2$
from $\mathcal{C}_1$ and $\phi_2 : H_2 \to H_1$ from $\mathcal{C}_2$
so that $\phi_2 \circ \phi_1$ is an inner automorphism of $H_1$ and
$\phi_1 \circ \phi_2$ is an inner automorphism of $H_2$.  (If there
is a compatible isomorphism $\phi_1$, choose $\phi_2$ to be the compatible
homomorphism in the conjugacy class of $\phi_1^{-1}$.)

We can certainly decide whether or not a homomorphism from
$H_1$ to itself is an inner automorphism, since this is simply a
matter of solving a finite system of equations, which we can do
by \cite{Dah_eq}.  

We therefore consider each of the possible compositions
of homomorphisms from $\mathcal{C}_1$ and $\mathcal{C}_2$ 
(and then from $\mathcal{C}_2$ and $\mathcal{C}_1$) in turn.
If we do not find a compatible isomorphism in this manner, then  
there does not exist a compatible isomorphism.  
In case there is a compatible isomorphism, it is straightforward
to reduce the list $\mathcal{C}_1$
 so that it contains only compatible isomorphisms.
\qed \end{proof}

\subsection{Tietze transformations}
We now investigate the situation when there are
nontrivial compatible primary splittings of $H_1$ and $H_2$.

Recall the concept of {\em Tietze transformations}.
Let $\mathcal Q = \langle \mathcal A \mid \mathcal R \rangle$ be a presentation.  Let $F$ be the free group on $\mathcal A$, and let $N$ be the normal closure of $\mathcal R$ in $F$.  Let $r \in N$.  Then $\langle \mathcal A \mid \mathcal R \cup \{ r \} \rangle$ defines a group isomorphic to the group defined by $\mathcal Q$.
Similarly, let $x$ be a letter not in $\mathcal A$, and let $w \in F$.  Then $\langle \mathcal A \cup \{ x \} \mid \mathcal R \cup \{ x^{-1}w \} \rangle$ defines a group isomorphic to the group defined by $\mathcal Q$. 
A {\em Tietze transformation} is the passage from a given (finite) presentation to another in either of the above two ways (in either direction).

The following theorem of Tietze is central to this part of the proof:

\begin{thm} \label{Tietze}
[Tietze; see \cite{LS}, Theorem II.2.1]
Two finite presentations define isomorphic groups if and only if it is possible to pass from one to another by a finite sequence of Tietze transformations.
\end{thm}

Given a finite group presentation $\langle \mathcal{A} \mid \mathcal{R} 
\rangle$, we can systematically enumerate
all finite presentations which define the same group as
$\langle \mathcal{A} \mid \mathcal{R} \rangle$, using Tietze 
transformations.  This gives us a list of all finite presentations
defining our given group.

\begin{defn} \label{def:exhibit}
Suppose that a finitely presented group $G$ admits a
decomposition as $G = A \ast_C B$, where $A$ and $B$ are
finitely presented and $C$ is finitely generated.  A finite 
presentation for $G$ is said to {\em exhibit}
the splitting $G = A \ast_C B$ if it is of the form
\[      \big\langle \mathcal{A}_1, \mathcal{A}_2, \mathcal{A}_3 \mid
\mathcal{R}_1, \mathcal{R}_2, \{ c = i_1(c) = i_2(c) \mid c \in 
\mathcal{A}_3 \} \big\rangle    ,\]
where $\langle \mathcal{A}_1 \mid \mathcal{R}_1 \rangle$ is a
presentation of $A$, $\langle \mathcal{A}_2 \mid \mathcal{R}_2 \rangle$
is a presentation of $B$, $C$ is generated by $\mathcal{A}_3$ and
the monomorphisms $i_1 : C \to A$ and $i_2 : C \to B$ are those
that define the amalgamated free product.

There is an entirely analogous definition of a presentation which
exhibits an HNN extension.
\end{defn}

The following result is clear.

\begin{lem} \label{lem:exhibit}
Suppose that a finitely presented group $G$ admits a nontrivial
splitting $G = A\ast_C B$ or $G = A \ast_C$, where $A$ and $B$ are finitely presented and $C$ is finitely generated.  Then there is a 
finite presentation of $G$ which exhibits this splitting.
\end{lem}

Using Tietze transformations, we can systematically enumerate
all finite presentations presentations of a finitely presented group $G$.
We will eventually
find a presentation exhibiting a splitting (except that we may not be able to check
whether the maps $i_1$ and $i_2$ are monomorphisms).   
However, we have to be
able recognize what kind of splitting we have found. 
In particular, for an arbitrary finitely presented group, we have no
way of knowing whether the splitting exhibited is nontrivial or not.

However, in the case of toral relatively hyperbolic groups, we have
the following result which surmounts these difficulties.

\begin{thm} \label{t:FindSplit}
Suppose that $G$ is a toral relatively hyperbolic group.  There is
an algorithm whose input is a finite presentation for $G$ along
with finite generating sets for representatives of conjugacy
classes of subgroups in an abelian peripheral structure for $G$ and
which terminates if and only if there is a nontrivial compatible primary
splitting of $G$.  In case it terminates, it terminates with a presentation
which exhibits a nontrivial compatible primary splitting, along with a proof that the
splitting is primary, compatible and nontrivial.
\end{thm}

\begin{proof}
We enumerate finite presentations of $G$ using Tietze transformations
searching for presentations which exhibit splittings.  In order to
find a presentation which exhibits a primary splitting if one exists, in parallel 
we continue to enumerate presentations, whilst examining those 
presentations which exhibit splittings to decide the following:
\begin{enumerate}
\item \label{S1} whether the splitting has abelian edge groups;
\item \label{S2} which of the edge groups are contained in parabolic subgroups;
\item \label{S3} which edge groups are cyclic;
\item \label{S4} whether the cyclic edge groups are maximal cyclic;
\item \label{S4.5} whether the noncyclic edge groups have
elements with roots not in the edge group;
\item \label{S5} whether the splitting is nontrivial;
\item \label{S6} whether noncyclic abelian subgroups of $G$ are elliptic; and
\item \label{S7} whether the peripheral subgroups of $G$ are elliptic (and thus
if the splitting is compatible).
\end{enumerate}
Note that in case $G$ is torsion-free hyperbolic, we do not require
the statements (\ref{S2}), (\ref{S3}) or (\ref{S6}).
It also is worth noting that all abelian subgroups of toral relatively hyperbolic groups
are finitely generated. Furthermore, the vertex groups of a primary 
splitting are toral relatively  hyperbolic, by Theorem \ref{vertexgpsRH},
so all nontrivial primary splittings have (finite) presentations which 
exhibit them.

We now describe a procedure which will terminate if the splitting
is primary.  Note that by running this in parallel with a search of
further presentations which exhibit splittings (and investigating
these in parallel also) if there is a nontrivial primary splitting we
will find one, together with a proof that the splitting is nontrivial
and primary. We only proceed to consider an item in the above list
if all of the previous items have been determined to hold.

A solution to the word problem allows us to decide (\ref{S1}). 
A simple system of equations can be used to decide (\ref{S2}), using
conjugacy and the malnormality of the parabolic subgroups.
Having determined which edge groups are contained in parabolics,
those that are not must be cyclic.  It is straightforward to determine whether
or not a finitely generated subgroup of a parabolic is cyclic, and thus we can decide
(\ref{S3}).  
Once we know whether a cyclic edge group is parabolic or not, we can decide
if the edge group is maximal as follows:  If the cyclic edge group is parabolic,
being maximal in the parabolic is equivalent to being maximal in $G$, and this
is straightforward to decide.  If the cyclic edge group $E$ is not parabolic, a 
theorem of Osin \cite[Theorem 1.16(3)]{Osin} states that it is possible to decide whether 
or not there is a nontrivial root of a generator of $E$. Thus we  can decide (\ref{S4}).

If an edge group $E$ is noncyclic, it is conjugate into a unique
parabolic subgroup, $P$ say. The group $P$ can be found
by a simple enumeration
process, and a solution to the word problem.  There are roots of some
$\gamma \in E$ which lie in $G \smallsetminus E$ if and only
if there are roots of $\gamma$ which lie in $P \smallsetminus E$
(since $P$ is malnormal).  This is the case if and only if $P/E$ is
not free abelian.  Since $P$ and $E$ are free abelian groups, and the generators
of $E$ are given as words in the generators of $P$, to decide if $P/E$ is free abelian is equivalent to decide if $E$ is a direct factor of $P$: if it is a direct factor we eventually see it on a certain basis of $P$, and if not, we will eventually find an element outside $E$ that has a power in $E$.
Therefore, we can decide  (\ref{S4.5}).

Once we know that the cyclic edge groups are maximal, it is straightforward
to decide whether or not the inclusions of edge groups into vertex groups are proper.
Thus we can decide (\ref{S5}).
For (\ref{S6}) and (\ref{S7}), suppose that $A$ is a parabolic or 
peripheral subgroup of $G$, and suppose that $\mathcal B$ is a basis
for $A$ (which we can always find).  We may search conjugates of 
the vertex groups of $G$ for one which contains the elements of $\mathcal B$.  This procedure will terminate if (\ref{S6}) and (\ref{S7})
hold, which is all that is required of this algorithm by the
above discussion.
\qed \end{proof}

The algorithm from Theorem \ref{t:FindSplit} above is the
complementary algorithm to that from Theorem \ref{theo;algo_comp}.

\begin{defn} \label{data}
By the {\em data} of a graph of groups we mean a description of the underlying
graph, presentations for the vertex groups and generators for the edge groups, a list of the
Bass-Serre generators, and relations corresponding to the usual
description of the fundamental group of a graph of groups.

If the graph if finite, vertex groups are finitely presented and the edge groups finitely generated,
then the data of a graph of groups is clearly finite.  When we say an algorithm {\em outputs}
a graph of groups decomposition of a finitely presented group $G$, we mean that it outputs the (finite) data of a
graph of groups, as well as an explicit
isomorphism from the fundamental group of the graph of groups to $G$ (with the given finite
presentation).
\end{defn}

\begin{remarknum}
We will always find our splittings by enumerating Tietze transformations.  By keeping a record
of the Tietze transformations performed, we can always exhibit an explicit isomorphism to
the presentation we started with.
\end{remarknum}

\begin{thm} \label{t:DecideSplit}
There exists an algorithm whose input is a finite presentation for $H$, a freely
indecomposable toral relatively hyperbolic group, along with a primary peripheral
structure for $H$, and which outputs ``yes'' or ``no''
depending on whether or not $H_1$ admits a nontrivial compatible primary splitting.

In case the algorithm terminates with ``yes'', it also outputs a finite presentation which
exhibits a nontrivial compatible primary splitting, along with a proof that the
splitting is primary, compatible and nontrivial.
\end{thm}
\begin{proof}
Run in parallel the algorithm from Theorem \ref{theo;algo_comp} applied to 
homomorphisms from $H_1$ to itself, and the algorithm from Theorem \ref{t:FindSplit}.
If $H_1$ admits a nontrivial compatible primary splitting, Theorem \ref{t:FindSplit}
assures us that we will find it.

However, if $H_1$ does not admit a nontrivial compatible primary splitting then
the algorithm from Theorem \ref{theo;algo_comp} terminates, by Theorem
\ref{NoTermiffSplit}.  If the algorithm from Theorem \ref{theo;algo_comp} does terminate
then Theorem \ref{NoTermiffSplit} also implies that $H_1$ does not admit a 
nontrivial primary splitting.
\qed \end{proof}

\begin{thm} \label{MaxSplit}
There exists
an algorithm whose input is a finite presentation for $H_1$, a freely indecomposable toral relatively hyperbolic group, and
whose output is a primary splitting of $H_1$, whose
vertex groups (equipped with the peripheral structure coming from
the adjacent edge groups) do not admit a nontrivial compatible
primary splitting.
\end{thm}
\begin{proof}
Begin with a trivial (empty) peripheral structure on $H_1$.  With Theorem \ref{t:DecideSplit}
 we can decide whether or not $H_1$ admits a nontrivial primary splitting.
If it does not, then the algorithm outputs the finite presentation for $H_1$, as a trivial
primary splitting.

Suppose that $H_1$ does admit a nontrivial primary splitting, that we
have found by the algorithm from Theorem \ref{t:FindSplit}. 
Give each of the non-abelian vertex groups of this splitting the peripheral
 structure coming from the adjacent edge groups.  Note that each of the vertex groups
is toral relatively hyperbolic, by Theorem \ref{vertexgpsRH}.  Apply the algorithm
from Theorem \ref{t:DecideSplit} in turn to each of the non-abelian vertex groups with
their attendant peripheral structures.

Proceeding in this manner, we continue to refine the splittings of $H_1$.
Generalized accessibility (see \cite[Main Theorem, p.451]{BFAccess}), and
the fact the $H_1$ is finitely presented, assures us that this process
will eventually terminate, and when it does we have found a splitting
as in the statement of the theorem.
\qed \end{proof}

The splitting that we have found in Theorem \ref{MaxSplit} above is not
suited to our purposes, because it is not canonical.  The way around this is
to turn this splitting into a {\em JSJ decomposition}, which will be canonical
enough for our needs.

\section{Finding the JSJ} \label{JSJSection}

The purpose of this section is to find the correct JSJ decomposition
for freely indecomposable toral relatively hyperbolic groups, and 
prove that there is an algorithm which finds the JSJ.   This JSJ decomposition
is the {\em primary JSJ decomposition}.

\begin{remarknum}
In case $\Gamma$ is torsion-free hyperbolic, the primary JSJ
decomposition is just the essential JSJ decomposition, whose
existence and uniqueness (up to slidings, etc.) was proved by Sela 
\cite[Theorem 1.8, p.569]{SelaGAFA}.

In case $\Gamma$ is toral relatively hyperbolic, we prove Theorem
\ref{JSJTheorem} in Section \ref{JSJ-RH} at the end of this paper.
The primary JSJ decomposition of a 
toral relatively hyperbolic group is closely
related to the abelian JSJ decomposition found for limit groups by
Sela \cite[Section 2]{SelaDio1}, although we need to be a little
more careful about the kinds of splittings we allow.
This is exactly
captured in the notion of a {\em primary} splitting, from Definition 
\ref{PrimarySplitting}.  
\end{remarknum}

The following theorem states the existence of a {\em primary JSJ
decomposition} for freely indecomposable relatively hyperbolic groups.
The reader who is unfamiliar with JSJ decompositions of groups
(and in particular with the terminology in Theorem \ref{JSJTheorem})
may wish to consult Section \ref{JSJ-RH}.  Also, recall the definition
of an {\em unfolded} splitting from Definition \ref{d:unfolded}, and see Theorem
\ref{t:unfolded} for a justification that unfolded splittings exist.

\begin{thm} \label{JSJTheorem}
[cf. Theorem 2.7, \cite{SelaDio1}; see also Theorem 7.1, \cite{RipsSelaJSJ}]
Suppose $\Gamma$ is a freely indecomposable toral relatively hyperbolic group.  There exists
a reduced unfolded splitting $\Lambda$ of $\Gamma$ with abelian edge groups, which we call a {\em primary JSJ decomposition} of $\Gamma$, satisfying the following:
\begin{enumerate}
\item Every canonical socket of a CEMQ subgroup of $\Gamma$ is conjugate to a vertex
group in $\Lambda$.  Every QH subgroup of $\Gamma$ can be 
conjugated into one of the CEMQ
subgroups of $\Lambda$.  Every vertex group in $\Lambda$
which is not a socket subgroup
of $\Gamma$ is elliptic in any primary splitting of $\Gamma$;
\item A one edge primary splitting $\Gamma = D \ast_A E$ or $\Gamma = D \ast_A$ which is 
hyperbolic in another primary splitting is obtained from the primary JSJ decomposition of $\Gamma$ 
	by cutting a surface  corresponding to a CEMQ subgroup of $\Gamma$ along an essential s.c.c;
\item Let $\Theta$ be a one edge primary splitting $\Gamma = D \ast_A E$ or $\Gamma = D \ast_A$,
which is elliptic with respect to any other one edge primary splitting of $\Gamma$.  Then $\Theta$ is
obtained from $\Lambda$ by a sequence of collapsings, foldings and
conjugations;
\item \label{JSJunique}  If JSJ$_1$ is another primary JSJ 
decomposition of  $\Gamma$ then JSJ$_1$ is obtained from $\Lambda$ 
by a sequence of slidings, conjugations and modifying boundary
monomorphisms by conjugations.
\end{enumerate}
\end{thm}

The remainder of this section is devoted to proving the
following result.

\medskip

{\noindent \it Theorem {\bf \ref{FindJSJ}.}}---
{\em
There is an algorithm
which takes a finite presentation for a freely indecomposable
toral relatively hyperbolic group, $\Gamma$ say, as input
and outputs a graph of groups which is a primary JSJ decomposition 
for $\Gamma$.
}

\medskip

The approach to proving Theorem \ref{FindJSJ} has two components.
The first is to find a maximal primary splitting of $\Gamma$.
That this can be done is the content of Theorem \ref{MaxSplit} above.

Having found this maximal splitting, we collapse some parts of this
graph of groups in order to find the JSJ decomposition.  (More precisely,
we collapse some parts to find the CEMQ pieces of the JSJ.)

Let $\Gamma$ be a freely indecomposable toral relatively hyperbolic
group and let $\Lambda_{Max}$ be the maximal splitting of $\Gamma$
obtained from Theorem \ref{MaxSplit}.  The splitting 
$\Lambda_{Max}$ is obtained from a primary JSJ decomposition
of $\Gamma$ by
cutting the CEMQ pieces along  maximal collections of disjoint, non-parallel
essential s.c.c. (this follows from the proof of Lemma \ref{HaveQH}
and Proposition \ref{HaveaJSJ} below).  We now decide which of the vertex groups arise in this manner,
and glue them back together in order to obtain a primary JSJ 
decomposition of $\Gamma$.

 In the perspective of the socket groups of the JSJ decomposition, 
 a {\em socket subgroup} of a group $G$ is a subgroup of $G$ which is
 obtained from a $QH$-subgroup of $G$ by adding all the roots of the
 punctures.  As an abstract group, a socket subgroup is just a free group.  
 However, socket subgroups always come equipped with
peripheral structures, corresponding to the cyclic subgroups
generated by the added roots of the punctures.

  A particular class of  groups can arise after this splitting process: they are 
  $QH$-subgroups with 
 some (non-essential) amalgamation over $\mathbb{Z}$  on the boundaries, and   
which have the following presentations as QH subgroups: 
 \[ \langle p_1,p_2,p_3 \mid p_1 (p_2)^{n_2} (p_3)^{n_3}  \rangle, \; n_2,n_3 \in \mathbb{N}  \]
 (where $p_1$, $p_2^{n_2}$ and $p_3^{n_3}$ correspond to the punctures of 
 a thrice-punctured sphere).
 As an abstract group, such a group $P$ is  free of rank $2$, with $\{ p_2, p_3 \}$
 a basis.  The group $P$ comes equipped with a peripheral structure consisting of the 
(conjugacy classes of) the cyclic subgroups generated by $p_1$, by $p_2$ and
by $p_3$.  It is straightforward to see that such a $P$ is freely indecomposable
relative to the peripheral structure and admits no nontrivial compatible primary
splitting.
 Call such a $P$ a {\em basic socket group} (equipped with the above peripheral
 structure).

\begin{prop}\label{prop;find_sockets} For each of the vertex groups $V$ in $\Lambda_{Max}$, it is decidable 
  whether or not
 $V$ admits a compatible isomorphism to any of the  basic socket groups, and in such case, determines such an isomorphism.  Moreover, the  basic socket group to which $V$ can be compatibly isomorphic is uniquely defined.   
 \end{prop}

\begin{proof}
The group $V$ can only be compatibly isomorphic
to a basic socket subgroup if it is (abstractly) isomorphic to a free group
of rank $2$, which we can decide as follows:
By Proposition \ref{vertexgpsRH}, $V$ is a toral relatively hyperbolic
group.  Therefore, it is decidable whether or not $V$ is a hyperbolic
group.  In case $V$ is hyperbolic we may compute its Grushko decomposition, by Gerasimov's algorithm (see also \cite{DahG_freeprod}). 

In case $V$ is isomorphic to a free group of rank $2$, it is possible to compute an explicit isomorphism, and thus to find a basis $(e_1,e_2)$ of the group, and express each peripheral subgroup  in this basis. The problem is then to solve: given a peripheral subgroup generated by $c$, does there exists a basis $a,b$ of the free group, in which $a^nb^m = c$ for some $n$, and $m$ in $\mathbb{N}$. Although this is not a sentence of the first order theory, the problem can be solved as follows.  

  The element $c$ is equal to $a^nb^m$, with $a$ and $b$ primitive elements, if and only if  $c=\alpha \beta$ for some $\alpha$ and $\beta$ such that   $[\alpha,a] =1 = [\beta, b]$. Indeed, if $[\alpha,a] =1 $, then $a$ and $\alpha$ have a common root, but $a$ is primitive. 

 Moreover, by a classical characterization (see Proposition 5.1 in \cite{LS}),
 $(a,b)$ is a basis if and only if there exists $g$ such that $[a,b] =
 ([e_1,e_2]^{\pm 1})^g$.  Thus, the problem   is equivalent to decide whether
 there exists $g, \alpha, \beta, a, b$ such that  $[a,b] = g ([e_1,e_2]^{\pm
 1}) g^{-1}$,  $c=\alpha \beta$ and $[\alpha,a] =1 = [\beta, b]$. This can be decided by Makanin's algorithm (see \cite{Makanin}). 

It remains to see that the  basic socket group to which $V$ can be compatibly
isomorphic is uniquely defined.  Each basic socket $V$ (with its peripheral
structure) has at least one peripheral
subgroup $\langle c\rangle$, and possibly one or two other, up to
conjugacy. The quotient of $V$ by (the normal closure of) any of these peripheral
subgroups is either infinite cyclic, or a generalized Baumslag-Solitar group
$\langle a,b | a^n=b^{-m} \rangle $ (when quotienting by $c$, and in this
case, $n$ and $m$ are the orders of the roots added). If neither $n$
nor $m$ are $1$, then the generalized Baumslag-Solitar group is rigid (see
\cite{Levitt_rigid}), meaning that $n$ and $m$ are characteristic of the group, which makes the structure of
basic socket unique. If now $n= 1$, this means that  $\langle a\rangle$ is in
the peripheral structure (otherwise, the basic socket has a free compatible
splitting), thus preserved by any compatible isomorphism. Then  $m$ is
characterized to be the order of the finite quotient of $V$ by $a$ and $c$.
\qed \end{proof}

Apply the algorithm from Proposition \ref{prop;find_sockets} to each of the
vertex groups of $\Lambda_{Max}$ to determine which vertex groups (with
the peripheral structure coming from adjacent edge groups) are
basic socket groups.

We can now decide whether there are adjacent vertex groups in $\Lambda_{Max}$
which are compatibly isomorphic to basic socket groups, and are joined by an edge
group corresponding to boundary components of the underlying surface (without
roots attached).  
 If there are such vertex groups, collapse
the edge joining them, to get a larger socket subgroup.  Proceed in this manner to obtain a
splitting $\Lambda'_{Max}$ which contains
socket  subgroups which are attached only to vertex groups of $\Lambda_{Max}$ which are
not compatibly isomorphic to any socket group (basic or otherwise).

We claim that the splitting $\Lambda'_{Max}$ is a primary
JSJ decomposition of $\Gamma$ (see Proposition \ref{HaveaJSJ}
below).  Call the vertex groups of $\Lambda'_{Max}$ which have been created
from more than one basic socket vertex group by gluing along boundary components
as described above the {\em socket} vertex groups of $\Lambda'_{Max}$.  The
socket vertex groups contain canonical maximal $QH$-subgroups.

Define
a splitting $\Lambda_{QH}$ of $\Gamma$ as follows: first, refine $\Lambda'_{Max}$ by splitting each socket vertex group along the boundary curves of the associated surface group (this is not an essential splitting !); second, collapse all the
edges in $\Lambda'_{Max}$ which are  not  adjacent to 
 any of the obtained  surface  vertex groups.  The {\em $QH$-subgroups} of $\Lambda_{QH}$
 are those obtained from socket subgroups of $\Lambda'_{Max}$ (which in particular are built 
 from more than one basic socket subgroup).
  
  See \cite[Theorem 4.21, p.87]{RipsSelaJSJ}, or 
Theorem \ref{EQD} below, for information about
the essential quadratic decomposition of $\Gamma$.  The key
step to proving that $\Lambda'_{Max}$ is a primary JSJ decomposition
of $\Gamma$ is the following result.  The reader unfamiliar with the essential
quadratic decomposition of finitely generated groups from \cite{RipsSelaJSJ} is
advised to consult Section \ref{JSJ-RH} at this point.

\begin{lem} \label{HaveQH}
$\Lambda_{QH}$ is an essential quadratic decomposition of $\Gamma$
(see Theorem \ref{EQD}).
\end{lem}

\begin{proof}
Let $\Theta$ be an essential quadratic decomposition of $\Gamma$.
The edge groups of $\Theta$ correspond to edge groups adjacent
to CEMQ subgroups, which are in particular maximal essential
QH-subgroups.  Suppose that an edge group $E$ of $\Theta$ is not elliptic
in $\Lambda'_{Max}$.  Let $\Theta'$ be the one-edge splitting arising from
$\Theta$ by collapsing all edges of $\Theta$ except the one associated to $E$.
It is easy to see that there is a one-edge splitting $\Pi$ arising from collapsing
all edges of $\Lambda'_{Max}$ except one so that $E$ is hyperbolic in $\Pi$.
Since $\Pi$ and $\Theta'$ are primary one-edge splittings, Theorem \ref{NoEllHyp}
implies that $\Pi$ and $\Theta'$ form a hyperbolic-hyperbolic pair of splittings, so
since all noncyclic abelian subgroups are elliptic in a primary splitting, $\Pi$ is
a cyclic splitting.  Theorem \ref{EQD}.(ii) now implies that $\Theta'$ is obtained
from $\Theta$ by cutting a surface corresponding to a CEMQ along an
essential s.c.c.  This is clearly not the case, so we have shown  that
all edge groups of $\Theta$ are elliptic in $\Lambda'_{Max}$.

Let $\Sigma$ be a $QH$-subgroup of $\Lambda_{QH}$.  First note that $\Sigma$ is a 
QH-subgroup, so by Theorem \ref{JSJTheorem}.(1) there is a CEMQ vertex
group $\Sigma'$ of $\Theta$ so that $\Sigma$ is conjugate into
$\Sigma'$.  Let $S$ be the surface corresponding to $\Sigma$ and
$S'$ the surface corresponding to $\Sigma'$.  We suppose (by conjugating) 
that $\Sigma \le \Sigma'$.  By Theorem \ref{EQD}, we can consider $S$ to be a sub-surface of $S'$

We claim that, in fact, $\Sigma = \Sigma'$.  Suppose not.  Then
some boundary components of $S$ are boundary components of $S'$ 
and some boundary components of $S$ are essential simple closed
curves on $S'$.  Let $S_0$ be a component of $S' \smallsetminus S$,
and let $\Sigma_0$ be the subgroup of $\Sigma'$ corresponding
to $S_0$.

Suppose that all boundary components of $S_0$ lie on $S$.  In
this case it is clear that $\Sigma$ could be made larger in 
$\Lambda'_{Max}$, contradicting the construction of $\Lambda'_{Max}$.

Therefore, some boundary component $c$ of $S_0$ lies in 
$S' \smallsetminus S$.  Let
$E$ be the edge group of $\Theta$ adjacent to $\Sigma'$ corresponding to
$c$.  We have
already noted that $E$ is elliptic in $\Lambda'_{Max}$, which means by
Theorem \ref{NoEllHyp} that all edge groups of $\Lambda'_{Max}$ are
elliptic in the one-edge splitting $\Theta_0$ corresponding to $E$.
Suppose that $E$ is not conjugate into an edge group of
$\Lambda'_{Max}$, but is conjugate into a vertex group $V$.  Then $E$
induces a compatible primary refinement of $V$.  If $V$ is not a socket 
vertex subgroup of $\Lambda'_{Max}$ then this contradicts the maximality of
$\Lambda'_{Max}$ (since non-socket vertex groups of $\Lambda'_{Max}$
admit no compatible primary splittings).  
However, we claim that $V$ cannot be a socket vertex
subgroup of $\Lambda'_{Max}$.  Indeed, the only compatible primary refinements
of socket vertex subgroups arise from cutting the corresponding surface along
an essential s.c.c.  However,
$E$ is elliptic in any primary splitting, since it corresponds to a boundary
component of the surface of a CEMQ subgroup.  If $E$ corresponded to
cutting along an essential s.c.c. of the surface of $V$, then since socket vertex
groups were constructed from at least two basic socket subgroups, there would
be some splitting in which $E$ were hyperbolic.

We have proved that all edge groups of $\Theta$ adjacent to $\Sigma'$ 
are in fact conjugate into edge groups of $\Lambda'_{Max}$.
The subgroup $\Sigma_0$ (along with the roots attached to form socket subgroups
in the construction of $\Lambda'_{Max}$) corresponds to a subgraph of $\Lambda'_{Max}$
which is a basic socket or socket subgroup, and it is glued to $\Sigma$ along
a subgroup corresponding to two boundary components.  Therefore,
$\Sigma_0$ is nontrivial, we could make $\Sigma$ larger,
in contradiction to the construction of $\Lambda'_{Max}$.
Thus, $\Sigma = \Sigma'$.

We now claim that any CEMQ subgroup of $\Theta$ must be conjugate into
a socket subgroup of $\Lambda'_{Max}$, which proves that $\Lambda_{QH}$ must be an
essential quadratic decomposition of $\Gamma$, as required.

Suppose that $\Sigma$ is a CEMQ subgroup of $\Theta$.  Then by Theorem \ref{NoEllHyp}
and Theorem \ref{EQD}.(ii), each edge group of $\Theta$ is elliptic in $\Lambda'_{Max}$.
Suppose that $\Sigma$ is not conjugate into a vertex group of $\Lambda'_{Max}$. Then
$\Sigma$ admits a nontrivial primary splitting induced by $\Lambda'_{Max}$.  Since all edge
groups of $\Sigma$ are conjugate into edge groups of $\Lambda'_{Max}$, the vertex groups
of this splitting of $\Sigma$ (along with the roots attached to form sockets in $\Lambda'_{Max}$) 
correspond to subgraphs of $\Lambda'_{Max}$, are are either basic socket or socket subgroups.
Therefore, these basic (socket) subgroups can be merged to form a bigger socket subgroup,
in contradiction to the construction of $\Lambda'_{Max}$.  We have proved that any
CEMQ subgroup of $\Theta$ is conjugate into a socket subgroup of $\Lambda'_{Max}$,
which completes the proof.
\qed \end{proof}

\begin{prop} \label{HaveaJSJ}
The splitting $\Lambda'_{Max}$ is a primary JSJ decomposition of $\Gamma$.
\end{prop}

\begin{proof}
Let us check that $\Lambda'_{Max}$ satisfies the three first points of Theorem \ref{JSJTheorem}. 

First, every canonical socket subgroup $H$ of $\Gamma$ contains a $CEMQ$ subgroup, 
and all the roots of its puncture elements.  Using the first point of Theorem \ref{EQD}, and the proof of Lemma 
\ref{HaveQH} above, we see that $H$ is conjugate to a vertex group of $\Lambda'_{Max}$. 
Given a vertex group of  $\Lambda'_{Max}$, which is not a socket group, assume 
that it is not elliptic in a primary splitting of $\Gamma$. Since it admits no 
compatible primary splitting, one of its adjacent edge groups is hyperbolic in this splitting. 
 But, by Definition \ref{PrimarySplitting},   
 this edge group is cyclic, and by Theorem   \ref{EQD}(ii)
 it  must occur as a simple closed curve of a CEMQ group.  Since (by the proof of Lemma 
\ref{HaveQH} above) all of the edge groups adjacent to a CEMQ subgroup of $\Gamma$ are
elliptic in $\Lambda'_{Max}$, this contradicts the maximality of the socket groups  
of  $\Lambda'_{Max}$.  This proves the first point, and similarly, one gets the second point.

Let $\Theta$ be a one-edge primary splitting elliptic in any other primary splitting.  if the edge
group of $\Theta$ is cyclic, Theorem
\ref{EQD}.(iii) implies that the edge group of $\Theta$ is conjugate into a vertex group $V$ of 
$\Lambda_{QH}$.  However, $\Theta$ now induces a (possibly trivial) compatible primary
splitting of $V$.  Since $\Theta$ is elliptic with respect to any other primary splitting,
the edge group of $\Theta$ does not correspond to an essential s.c.c. on a surface corresponding
to a $QH$-subgroup of $\Lambda_{QH}$.  This, and the maximality of the splitting $\Lambda'_{Max}$
shows that the edge group of $\Theta$ is conjugate to an edge group of $\Lambda'_{Max}$.
Therefore $\Theta$ can be obtained (by collapses, conjugations and slidings) from
$\Lambda'_{Max}$. 
  If now the edge group of $\Theta$ is noncyclic abelian,  it is elliptic in $\Lambda_{QH}$. It is therefore conjugated to a splitting of a vertex group of $\Lambda_{QH}$, and the same conclusion holds.  

This ensures the third point.
\qed \end{proof}

Since we can algorithmically find the splitting $\Lambda'_{Max}$,  Proposition 
\ref{HaveaJSJ} implies Theorem \ref{FindJSJ}.

\begin{remarknum}
Since we have an explicit isomorphism between the
canonical socket subgroups of the JSJ decomposition $\Lambda'_{Max}$ 
and a socket group given by a `standard' presentation (see Definition 
\ref{QH}), we may assume that the canonical socket subgroups of 
$\Lambda'_{Max}$ are given with standard presentations.
\end{remarknum}

\section{Proof of the main result} \label{ProofMainThm}

The main result of this section is the following.

\begin{thm} \label{JSJIsoPb}
There exists an algorithm whose input is a pair of finite
presentations, each defining a freely indecomposable toral relatively hyperbolic group, and whose output is ``yes" or ``no" depending upon whether or not
the groups defined by presentations are isomorphic.
\end{thm}

Together with the main result of \cite{DahG_freeprod} (which computes Grushko decompositions), Theorem \ref{JSJIsoPb} 
completes the proof of Theorem \ref{Toral} (our solution to the isomorphism problem for toral relatively hyperbolic groups).

Suppose that $\Gamma_1$ and $\Gamma_2$ are toral relatively hyperbolic groups, given by finite
group presentations $\langle X_1 \mid \mathcal{R}_1 \rangle$ and $\langle X_1 \mid 
\mathcal{R}_2 \rangle$.  

By Theorem \ref{FindJSJ} it is possible to algorithmically find JSJ decompositions
$\Xi_1$ of $\Gamma_1$ and $\Xi_2$ of $\Gamma_2$.  The vertex groups of $\Xi_1$
and $\Xi_2$ are themselves toral relatively hyperbolic, by Theorem \ref{vertexgpsRH}.
They are either socket groups, abelian groups, or else we
call them {\em rigid}.   Define peripheral
structures on each of them  by taking adjacent edge groups and their conjugates.  
The key property of rigid vertex groups is that they do not
admit any nontrivial compatible primary splittings.

By the construction of the JSJ decomposition, 
we know whether or not a given vertex group is a socket group or an abelian group (in any case, 
being abelian is easily recognised, and being rigid is recognised by Theorem \ref{t:DecideSplit}).

\begin{remarknum} \label{r:BoundedMovesonJSJ}
Since noncyclic abelian subgroups are elliptic in all primary splittings, 
and maximal abelian subgroups are malnormal, there is no
nonempty sequence of sliding moves on a single edge (without 
backtracking) which  brings an edge back to where it began.  Therefore,
there
are only finitely many slidings that can be performed on this graph
in order that it remains a JSJ decomposition.  It is clear that all of these
moves can be performed algorithmically. Thus using these sliding moves, we can, from one JSJ decomposition $\Xi_2$ of $\Gamma_2$,  
compute effectively  
 the list $\Xi_2^1, \ldots,
 \Xi_2^k$ of all $JSJ$ decompositions of $\Gamma_2$, up to conjugation and modifying 
boundary morphisms  by conjugation.
\end{remarknum}

\begin{defn} \label{consistent}
Suppose that $\Lambda$ is a graph of groups decomposition
of a group $G$, and $\Lambda'$ a graph of groups decomposition
of $G'$, and that $\mu$ is a graph isomorphism between the graph
underlying $\Lambda$ and the graph underlying $\Lambda'$.  Let
$e$ be an edge in the graph underlying $\Lambda$, and let
$v$ be the initial vertex of $e$.  Denote by $\iota_e$ the 
injection of $G_e$ into $G_v$, given by the graph of groups, 
and by $\iota_{\mu(e)}$ the given
inclusion of $G_{\mu(e)}$ into $G_{\mu(v)}$.  Suppose that
$\pi_1$ is an isomorphism between $G_e$ and $G_{\mu(e)}$.

An isomorphism $\pi_2 : G_v \to G_{\mu(v)}$ is called {\em
consistent with $\pi_1$} if, on $G_e$, one has, up to conjugacy in $G_{\mu(v)}$,   
\[      \pi_2 \circ \iota_e  = \iota_{\mu(e)} \circ \pi_1.     \]
\end{defn}

Let us define a \emph{complete item} to be a finite collection  of the following form:
\begin{enumerate}
\item A  primary JSJ decomposition $\Xi_2^i$ of 
$\Gamma_2$;
\item A certain   graph isomorphisms $\mu$ from
the underlying graph of $\Xi_1$ to the underlying graph of $\Xi_2^i$,
along with:
\begin{enumerate}
\item for each edge $e$ in the graph of $\Xi_1$, an isomorphism
between the group associated to $e$ and the group associated to
$\mu(e)$; and
\item for each vertex $v$ in the graph of $\Xi_1$, an isomorphism
between the group associated to $v$ and the group associated
to $\mu(v)$ which is consistent with the above isomorphisms of
edge groups.
\end{enumerate}
\end{enumerate}

Note that one can easily make the list of all graph isomorphism between  the underlying graph of $\Xi_1$ and those of  $\Xi_2^i$, $i\leq k$.

We now explain how to find, if there exists one, a consistent isomorphism between two vertex groups. First there can only be a 
compatible isomorphism between vertex groups of the same type (rigid, abelian, or socket), and we know the type of each vertex.

Given, for two rigid vertex groups,  By Theorem \ref{IsoNoSplit}, it is possible
 to find a list containing a representative of each conjugacy class of compatible isomorphism.
  Moreover, given an  isomorphism compatible with adjacent edge groups, it is easy to decide whether the isomorphism is consistent in the sense of our definition (one need to check, for each edge group,  whether images by the isomorphism of its given generators  are simultaneously conjugate to the given generators in the target group). 

Given two abelian groups, with  Tietze enumeration, one can find basis for each of them, and check whether there is  
a consistent isomorphism, using elementary linear algebra.    

Given a pair of socket groups, we can find, using   Tietze enumeration, their standard presentations (see Definition \ref{QH}), and check whether the surfaces and the given roots of the boundary components are the same, up to conjugation. Note that a boundary element is conjugated to its inverse if and only if the surface is non-orientable.

Therefore,  using  $\Xi_2$, the conclusion of Remark \ref{r:BoundedMovesonJSJ}  and the discussion above, one can effectively decide whether a complete item exists, and find one if it exists.

By uniqueness of the JSJ decompositions (Theorem \ref{JSJTheorem},  point \ref{JSJunique}),  
if there is no complete item, then $\Gamma_1$ 
and $\Gamma_2$ are not isomorphic.

\begin{prop}
If there exists   a complete item  then $\Gamma_1$ and 
$\Gamma_2$ are isomorphic.
\end{prop}

\begin{proof}
Consider a complete item: a  pair of primary JSJ decompositions of $\Gamma_1$ and $\Gamma_2$, a graph isomorphism $\mu$, and a collection of 
isomorphisms between the edge groups and vertex groups which
are consistent in the sense of Definition \ref{consistent}.  For ease
of notation, we continue to refer to the primary JSJ decomposition of 
$\Gamma_2$ as $\Xi_2$, rather than $\Xi_2^i$.

Fix a vertex group $V_1$ of $\Xi_1$.  The graph isomorphism $\mu$
takes the vertex of $V_1$ to the vertex of $V_1'$, a vertex
group of $\Xi_2$, and we are given a compatible isomorphism 
 $\pi$ from $V_1$ to $V_1'$.   Having fixed our initial vertex, and an
isomorphism between $V_1$ and $V_2$, we have used up our
freedom of conjugating the graphs of groups. 

Choose an oriented edge $e$ adjacent to $V_1$, $\iota_e$ the boundary morphism of its group in $V_1$,  
	and let $V_2$ be the 
vertex group associated to the vertex at the other end of 
$e$ (possibly the vertices at either end of $e$ are the same,
in which case $V_1 = V_2$, and $\iota_e$ is related to the orientation of $e$).  Denote by $\pi_e$ the given isomorphism between the group of $e$ and that of $\mu(e)$.
 We have consistent isomorphisms between the edge groups and between the vertex groups. 
Hence, we can conjugate,  in $V_1'$,  
the  boundary morphism  $\iota_{\mu(e)}$ so that  it coincide with $\pi\circ \iota_e\circ  \pi_e^{-1}$.  

If the two vertices of $e$ are distinct, using a suitable conjugate of a consistent isomorphism between $V_2$ and $V_2'$,   
we can extend this into an isomorphism between the
subgroups of $\Gamma_1$ and $\Gamma_2$ corresponding
to the one edge subgraphs.  
If the two vertices are the same, we can also form an isomorphism between the groups of the  one edge subgraphs, by also conjugating the second boundary morphism of $\mu(e)$ in $V_1'$ so that it coincide with the image by $\pi$ of the second boundary morphism of $e$ (this is possible since $\pi$ is consistent).

 In both cases, all that is required to extend the consistent isomorphisms, is to modify the boundary
morphisms by conjugation in the graph of group $\Xi_2$. The obtained isomorphism, between the groups of the  one-edge subgraphs, is clearly consistent (in the sense of Definition \ref{consistent}) with respect to  the adjacent edge groups of the  vertices.

Proceeding in this manner, we can construct an isomorphism
between $\Gamma_1$ and $\Gamma_2$ by adding one
edge at a time and modifying boundary morphisms by
conjugation.
\qed \end{proof}

We have thus shown that there exists a complete item if and only if there is an isomorphism between  $\Gamma_1$ and $\Gamma_2$. 
Since, as we discussed,
 we can effectively decide the existence of a complete item, we have proved 
Theorem \ref{JSJIsoPb}.

This finally finishes our solution of the isomorphism
problem for toral relatively hyperbolic groups, modulo the
results in Sections \ref{part;proof_gene} and \ref{JSJ-RH}.

In particular, we have completed the solution to the isomorphism
problem for torsion-free hyperbolic groups.

\section{Hyperbolic manifolds} \label{HypSection}

We now turn to the homeomorphism problem for finite
volume hyperbolic $n$-manifolds, for $n \ge 3$.  Mostow-Prasad
Rigidity implies that this is equivalent to the 
isomorphism problem for their fundamental groups.  

Therefore, the main purpose of this section is to prove
the following theorem, from which Theorem \ref{HomeoHyp}
follows.

\begin{thm} \label{th:manifold}
The isomorphism problem is solvable for the class of 
fundamental groups of finite-volume hyperbolic $n$-manifolds,
for $n \ge 3$.
\end{thm}

\begin{remarknum}
It is not immediately obvious what it means to be given
a finite-volume hyperbolic manifold as input to an algorithm.
There are a number of possible answers to this question.
We could want the triangulation of a compact core.  Or
possibly a discrete subgroup of $SO(n,1)$, given by
a collection of generating matrices, say.\footnote{This approach
may run into problems of real (or complex) arithmetic.}
 However, Mostow-Prasad Rigidity implies that the fundamental group
determines the manifold, and hence it is enough to be given
a (finite) presentation for the fundamental group.  This is the
point of view we take.
\end{remarknum}

Let $\Gamma_1$ and $\Gamma_2$ be fundamental groups
of finite-volume hyperbolic manifolds (of dimension at least $3$),
$M_1$ and $M_2$, say.

We apply the following analysis to $\Gamma_1$.  

First, the cusp groups of $\Gamma_1$ are all finitely generated
virtually abelian.  By a result of Hummel,  \cite{Hummel},  there is a finite index
subgroup of $\Gamma_1$ which has abelian cusp groups.

We can enumerate the homomorphisms from $\Gamma_1$ to 
finite groups, and find a presentation for the kernel of each such
map.  One way to find a presentation for such a kernel $K$ is to build the presentation
$2$-complex of $\Gamma_1$ and then build the finite cover corresponding
to $K$.  Contracting a maximal tree in this cover gives a presentation $2$-complex
for $K$, from which a presentation is obvious.

  By then running the one-sided algorithm which recognizes
toral relatively hyperbolic groups, in parallel with searching further
maps to finite groups, we will eventually find a finite index subgroup
$H_1$ of $\Gamma_1$ so that $H_1$ is toral relatively hyperbolic.

Similarly, find a finite index subgroup $H_2$ of $\Gamma_2$ which
is toral relatively hyperbolic.  

Let $N_1$ be the intersection of all subgroups of $\Gamma_1$ of
index at most $d = [\Gamma_1:H_1].[\Gamma_2:H_2]$, and
$N_2$ the intersection of all subgroups of $\Gamma_2$ of index
at most $d$.  Both $N_1$ and $N_2$ are toral relatively hyperbolic,
are normal in $\Gamma_1$ and $\Gamma_2$ respectively, and
finite presentations for them can be effectively found.

\begin{remarknum}
Note that we can find the relatively hyperbolic structure for $N_1$ and
$N_2$, and in particular a basis for the cusp subgroups.  Therefore,
we can decide the dimensions of the manifolds $M_1$ and $M_2$.
\end{remarknum}

Clearly, if $\Gamma_1$ and $\Gamma_2$ are isomorphic then
$N_1$ and $N_2$ are isomorphic, as are $G_1 = \Gamma_1/N_1$
and $G_2 = \Gamma_2/N_2$.  

By Theorem \ref{Toral}, we can effectively decide whether or not
$N_1$ and $N_2$ are isomorphic.  Also, since $G_1$ and $G_2$
are finite, we can effectively decide whether or not they
are isomorphic.  Therefore, supposing that both pairs
are isomorphic, we need to decide whether the appropriate extensions
(namely $\Gamma_1$ and $\Gamma_2$) are isomorphic.

We want to classify the extensions of $N$ by $G$, where $N$
is the toral relatively hyperbolic group of finite index in $\Gamma$
(and $N$ is characteristic), and $G = \Gamma / N$.  The extension
gives a homomorphism $G \to \mbox{Out}(N)$ via the action of
$\Gamma$ on $N$ by conjugation.  Recall the following result:

\begin{prop} \cite[Corollary IV.6.8, p. 106]{Brown}
\label{extensions}
If $N$ has a trivial center then there is exactly one extension of
$G$ by $N$ (up to equivalence) corresponding to any 
homomorphism $G \to \mbox{Out}(N)$.
\end{prop}

We apply Proposition \ref{extensions} to the pairs $(N_1, G_1)$ and
$(N_2,G_2)$ to decide whether they determine equivalent extensions
(equivalent extensions, in particular, correspond to isomorphic groups).

\begin{remarknum} The fundamental group of a finite volume hyperbolic manifold is never 
virtually abelian, and therefore,  
the center of $N_1$ is trivial. 
 Hence the hypotheses of Proposition \ref{extensions}
are satisfied.
\end{remarknum}

Being fundamental groups of finite-volume hyperbolic
manifolds of dimension at least $3$, neither $N_1$ nor $N_2$ admits
a nontrivial splitting over an abelian group (see, for example, 
\cite[Theorem 1.6(i)]{Bel}).  
Therefore $\mbox{Out}(N_1)$ is finite, and the algorithm from Theorem  \ref{theo;algo_comp} must terminate when applied to homomorphisms from $N_1$ to itself, and we can find
representatives for the conjugacy classes of isomorphisms from
 $N_1$ to itself.  In other words, we can
effectively calculate $\mbox{Out}(N_1)$, and also the homomorphism
from $G_1 \to \mbox{Out}(N_1)$ determined by the conjugation action
of $\Gamma_1$ on $N_1$.  

Similar considerations obviously apply for $N_2$, $\mbox{Out}(N_2)$,
and the homomorphism $G_2 \to \mbox{Out}(N_2)$.

Proposition \ref{extensions} now implies that $\Gamma_1$ and $\Gamma_2$ are isomorphic if and only if
\begin{enumerate}
\item $N_1$ and $N_2$ are isomorphic;
\item $G_1$ and $G_2$ are isomorphic; and
\item After composing with an automorphism of $G_2$ and an 
automorphism of $\mbox{Out}(N_2)$, the homomorphisms 
$G_1 \to \mbox{Out}(N_1)$ and $G_2 \to \mbox{Out}(N_2)$ determined
by conjugation of $\Gamma_1$ (respectively $\Gamma_2$) on
$N_1$ (resp. $N_2$) are the same.
\end{enumerate}
We can effectively decide each of these issues and thus we have
proved Theorem \ref{th:manifold}.  Mostow-Prasad Rigidity now immediately implies Theorem \ref{HomeoHyp}.

\section{Proof of Proposition \ref{prop;Omega} in the case of toral relatively
  hyperbolic groups}\label{part;proof_gene}

     We give here a proof of the technical Proposition \ref{prop;Omega} in the
     generality of toral relatively hyperbolic groups. We advise reading the
     paragraph \ref{Long_and_short}.

     We use the notation of
    Section \ref{part;not}, which we briefly recall now, for convenience.  Thus $H_2$ is a non elementary toral relatively hyperbolic group, generated by $S$, with maximal parabolic subgroups the free abelian groups $G_1, \dots, G_p$ and their conjugates. For each we are given a basis $S_i$ (with an order on each, so that a lexicographical order is defined on words). The graph $\widehat{Cay}(H_2)$ is the coned-off Cayley graph, its vertices of finite valence are identified with the elements of $H_2$, thus $H_2 \subset Cay (H_2) \subset \widehat{Cay}(H_2)$, and its distance is denoted by $d$. 

This distance is not to be confused with the word distance on $H_2$ for the generating set $S$, which we denote by $dist$.  

The free product $F=  F_S * G_1*\dots * G_p$ maps onto $H_2$. Recall that, in this setting, the normal form of an element of $F$ is the reduced word in the alphabet $S^{\pm 1}\cup S_1^{\pm 1}\cup \dots \cup S_p^{\pm 1}$ such that every subword in   $(S_i^{\pm 1})^*$ is minimal for the lexicographical order. The normal form of an element in $F$ labels a path in   $\widehat{Cay}(H_2)$ as explained in Remark \ref{rem;labels}.  To avoid confusion, let us precise that it also labels a path in $Cay(H_2)$, and both paths are related: we get the path in $\widehat{Cay}(H_2)$ by replacing every maximal subpath in a coset of a parabolic subgroup by the pair of edges passing through the infinite valence vertex associated to this coset. Unless otherwise precised, the labeled paths we will use  are in $\widehat{Cay}(H_2)$ and defined according to  Remark \ref{rem;labels}

Note that the length of the labeled path in  $\widehat{Cay}(H_2)$  is different from the length of the normal form seen as a word (in general the latter is longer): a maximal subword in  $(S_i^{\pm 1})^*$ labels a single pair of edges adjacent to an infinite valence vertex.  This difference is similar to that of $d$ and $dist$:  if one add the restriction that the labeled path has no angle larger than an explicit constant $K$, then the length of the normal form is at most $(K'+1)$ times that of the path, where $K'$ is the maximal length of an element of $G_i$ labeling a pair of edges of angle at most $K$.

 The set $\mathcal{L} \subset F$ is the normalized rational language consisting of elements labeling local quasi-geodesics without detours in $\widehat{Cay}(H_2)$ (see Definition \ref{def;lang_L} for details and constants).

     The structure of the proof is similar to the one for hyperbolic
     groups. First (in Paragraph  \ref{para;QP_rh}) we define the finite set
     $\mathcal{QP}$, and the languages
     $\mathcal{L}_h$. 
     
     Defining  $\mathcal{QP}$ to be a sphere for $d$, as in the hyperbolic
     case,     would be unwise, since
     this is an infinite set. Defining it to be a sphere for $dist$ would miss
     the flavor of hyperbolic geometry, and would make arguments more
     complicated. Thus we define it in terms of the distance $d$, and of
     angles in the coned-off graph. 

     Also,  for $h\in \mathcal{QP}$, the language that is really analogous to
     $\mathcal{L}_h$ in the hyperbolic case is here called
     $\mathcal{L}_{1,h}$. But, perhaps surprisingly, this does not suffice in
     ensuring that for ``most'' elements  $h$, the element $hgh^{-1}$ has a quasi-prefix that
     is a prefix of $h$ (as it was the case for hyperbolic groups, see Lemma
     \ref{lem;prefix}).

     If $h$ is long for $dist$, but $hgh^{-1}$ is short for $d$, a phenomenon can
     occur, the ``cascade effect'', 
     which I. Bumagin identified
     in her study of the conjugacy problem for relatively hyperbolic
     groups in \cite{Bumagin}. This justifies the introduction of an additional language
     $\mathcal{L}_{2,h}$ to define $\mathcal{L}_h=\mathcal{L}_{1,h}\cup
     \mathcal{L}_{2,h}$.

     We prove that each $\mathcal{L}_h$ is normalized rational.
     The proof of the regularity of
     the languages is a little more complicated than in the hyperbolic case, due to the presence of the
     new sets $\mathcal{L}_{2,h}$.

        In Subsection \ref{part;finiteness_rh}, we prove the second point 
        of Proposition \ref{prop;Omega}: 
        in all conjugacy class of compatible homomorphism, 
        only finitely many homomorphisms have an acceptable lift 
        satisfying $\Omega$. 
        As we already mention, the analogue of Lemma \ref{lem;prefix}, 
        which is Lemma \ref{lem;prefix_rh} requires a study of the ``cascade
        effect'', and is rather involved.

        Finally, Subsection  \ref{part;exists_rh} serves to prove the first point
        of Proposition \ref{prop;Omega}.  
        We find a  
        quantity $Q$  associated to every homomorphism, and we show that 
        in any conjugacy class 
        there is a homomorphism  minimizing  
        $Q$.  On the other hand, if an acceptable lift of a homomorphism
        $\phi$ contradicts $\Omega$, 
        then there is come conjugate $\phi^h$ with $Q(\phi^h) < Q(\phi)$.  
        The argument is similar to the one in the
        hyperbolic case, but there are four more cases to study.

        Let us remark that, in the hyperbolic case,  $Q (\phi)$ was simply the
        maximum between two distances. Here again, we need to take angles into
        account in a sensible way to define $Q$.

      \subsection{The set $\mathcal{QP}$, and the languages 
        $\mathcal{L}_h$, for $h\in \mathcal{QP}$} \label{para;QP_rh}
	With notations as above, let $\mathcal{L}_p \subset \mathcal{L}$ be the set of the elements of $\mathcal{L}$
        that map on an element of $(G_1 \cup \dots \cup G_p )\subset H_2$  under $F \to H_2$ 
	(for convenience $\mathcal{L}_p(i)$ is the subset of $\mathcal{L}_p$ that maps on elements of $G_i$).

        \begin{lem}
          The language $\mathcal{L}_p$ is normalized rational.
        \end{lem}

       \begin{proof}
		Let $\tilde{g} \in \mathcal{L}_p(i)$, its image in $H_2$ is in  
		 $G_i$. 
		We can decompose the normal form of $\tilde{g}$ as the concatenation of  ShortLex words in $(S_i^{\pm 1})^*$, and 
		of normal forms of elements in $\mathcal{L}$, with
          	no prefix nor suffix in $(S_i^{\pm 1})^*$, and that define elements of $G_i$ in $H_2$. 
		Conversely, if an element  $\tilde{g}\in \mathcal{L}$ has such a decomposition, then it defines 
		clearly an element of $G_i$, hence it is in  $\mathcal{L}_p$

	The elements of $\mathcal{L}$ of this decomposition  define paths in
          $\widehat{Cay}(H_2)$ (labeled by their normal form) 
	between two points that are at distance $2$, by definition of the decomposition. Hence, these
          paths are of length at most $L'_1$.   We can bound the maximal angle of such a path in terms of $L'_1$ and $r$. Indeed, if the path makes an angle greater than $L'_1+3$ at a vertex $v$, then, the loop consisting of the path together with the pair of edges between its end points has to pass twice at that vertex, because it is of length $L'_1+2$ (a simple loop of this length contradicts the definition of angle). 
In this case, the absence of $r$-detours gives the desired bound.  
	Therefore, the length  of the normal forms 
		of these elements  
	is bounded by some  computable constant. Hence  there are only
          finitely many such elements, and  they are computable with a
          solution of the membership problem for $G_i$ (here this is equivalent to commuting with an element of the basis of $G_i$, hence can be solved by solving a certain equation).

          This shows that membership of $\mathcal{L}_p(i)$  is characterized on the normal forms by being the concatenation
          of several words belonging in two regular languages. 

		Therefore, the
          language of the normal forms of the elements of  $\mathcal{L}_p$ is
          the finite union (over $i=1, \dots, p$) 
		of regular languages.  Therefore, it is a regular language, and this makes  $\mathcal{L}_p$ normalized rational.
\qed \end{proof}

        The constant $\epsilon$ is chosen as in Proposition \ref{prop;stab_con}, for 
        $(L'_1,L'_2)$-quasi-geodesics without $r$-detours.

        We choose constants  $\eta=3\epsilon + 300\delta $,  
        $\rho=
        20(\epsilon+\eta+20\delta)$,  
        and  $\Theta
        \geq 4(\eta +\epsilon + 100\delta) + 200\delta $, 
        and sufficiently large
        such that no angle at a vertex of finite valence in $\widehat{Cay}(H_2)$
        is greater than $\Theta$.

\begin{defn} \label{def:QP:RH}
        Let $\mathcal{QP}$ be the subset of $H_2$ consisting of the elements $h$
        with $d(1,h) \leq \rho$, such that 
        either (i)  $d(1,h) \geq \rho -1$ and $\MaxAng[1,h ] < \Theta$; or 
(ii) for some choice of geodesic $[1,h]$, if $v$ is the vertex adjacent to $h$ in $[1,h]$ 
then $\Ang_v[1,h] \in [\Theta,\Theta+1]$ and if $[1,v] \subset [1,h]$ is the sub-geodesic
from $1$ to $v$ then $\MaxAng[1,v]< \Theta$.
\end{defn}

        The set $\mathcal{QP}$ is contained in a certain cone of 
        $\widehat{Cay}(H_2)$, and
        therefore is a finite (computable) set. Recall that $H_2$ has been identified with the set of vertices of finite valence of 
  	$\widehat{Cay}(H_2)$.

\begin{defn} \label{def:quasi-prefix:RH}
         Let $h$ be an element of  $\mathcal{QP}\subset H_2$. 
         Given an element 
         $\tilde{g} \in (\mathcal{L}\setminus \mathcal{L}_p) \subset 
        F$, we say that $h$ is a {\em quasi-prefix} of $\tilde{g}$ 
         if a path of  ${Cay}(H_2)$ labeled by 
         the letters of a shortest word defining 
         $\tilde{g}$   contains
         some edge $e$ with $h \in \Cone_{\eta,\eta}(e)$  
	(we emphasize that we are considering the path in $Cay H_2 \subset \widehat{Cay}H_2$, 
	thus without infinite valence vertex, labeled by some shortest word, and not necessarily minimal for the 
	lexicographical order).

	We say that $h$ is a {\em quasi-suffix}  of $\tilde{g}$  if it is a quasi-prefix of  $\tilde{g}^{-1}$, and if it is both a quasi-prefix and a quasi-suffix, we say that they are disjoint if the edge  from the definition of quasi-prefix, lies on the path before the edge from the definition of quasi-suffix.

\end{defn}        

\begin{defn} \label{def:Lh:RH}
         Given $h \in \mathcal{QP}$, 
         we define the language $\mathcal{L}_{1,h} \subset
         (\mathcal{L}\setminus \mathcal{L}_p)$ to consist of 
         those elements $\tilde{g}$ of $\mathcal{L}$       so that $h$ is a quasi-prefix and a quasi-suffix of $\tilde{g}$, and these quasi-prefix and quasi-suffix are disjoint. 

\end{defn}

         \begin{lem}
           For all $h$, the language $\mathcal{L}_{1,h} \subset F$ 
           is normalized rational.
        \end{lem}

\begin{proof}
          One needs to show that the language $L_{h}$ of 
          the normal forms of the
          elements of $\mathcal{L}_{1,h}$, is regular.

          Let $\mathcal{E}$ be the set of edges $e$ in $Cay(H_2)$ 
          (between elements  
          of $H_2$), such that the cone of radius and angle $\eta$ centered at 
          $e$, in  $\widehat{Cay} (H_2)$,   contains $h$.

          Let $\mathfrak{F}$ be the set of normal forms  of
          the elements of $\mathcal{L} \subset F$ labeling  paths in
          $\widehat{Cay}(H_2)$  
          whose last edge is in  $\mathcal{E}$.  
          Note that, by Proposition  \ref{prop;stab_con}, every such normal form labels a quasi-geodesic path in a 
	  conical (hence finite)  neighborhood of a geodesic. Therefore,  $\mathfrak{F}$ is a finite set, and
          it is computable with a  solution to the word problem in $H_2$.

	  For $f \in  \mathfrak{F}$, consider the decomposition $f=rp$ where $p$ is the maximal suffix of $f$ in some $(S_i^{\pm 1})^*$, for $i\geq 1$ (that is, in some abelian factor of the free product). The word $p$ might be empty, in this case, any index $i$ (that is,  any abelian factor) is suitable for the argument. 
Let us say that $p \subset p' \in (S_i^{\pm 1})^*$ if $p'$ is reduced, minimal in lexicographical order,  and contains all the letters of $p$ with multiplicity.   

	Let $L_{h}(pre)[f] = \{ rp', \,  p\subset p' \}$.  This set is the set of normal forms that equal some element of $F$ labeling in $Cay(H_2)$ a path starting by the path labeled by $f$ (note that they may not label themselves this path: a permutation of the letters of $p'$, in order to make $p$ appear, may be necessary).

	  Given $f$, the set $L_{h}(pre)[f]$ is clearly recognized by a finite state automaton (it is finite). Let $W_f$ be the set of normal forms without prefix in $(S_i^{\pm 1})^*$. By definition, to have $h$ as a quasi-prefix is equivalent to have normal form in $\bigcup_{f  \in  \mathfrak{F}}  L_{h}(pre)[f] \cdot W_f$ where $\cdot$ is the concatenation of word-languages.

      	Similarly, we define $L_{h}(suf)[f] =\{ (p'')^{-1}r^{-1}, \,  p\subset p'' \} $. We claim that, if $L_h$ is the set of normal forms of elements of $\mathcal{L}_{1,h}$, then       $L_{h} = (\bigcup_{f  \in  \mathfrak{F}}   L_{h}(pre)[f] \cdot W_f )  \cdot   \bigcup_{f  \in  \mathfrak{F}}   L_{h}(suf)[f]$ and is reduced and ShortLex.   The inclusion $\subset $ is direct from the definition of $\mathcal{L}_{1,h}$. Conversely, if $w= r_1p' w' (p'')^{-1} (r_2)^{-1}$  with the above property, and moreover is reduced and shortest for the lexicographical order, then  it has $h$ as quasi-prefix, by definition of  $\mathfrak{F}$, and    $w^{-1} = r_2p'' w'^{-1} p'^{-1} r_1^{-1}$ satisfy $w'^{-1} \in  W_{f'}$  for the $f'$ associated to $r_2p''$ (otherwise it would not be reduced and ShortLex). Hence $h$ is also a quasi-suffix, and clearly, quasi prefix and quasi suffix are disjoint.

   This makes $L_h$ a regular language, hence the result.
\qed \end{proof}

\begin{defn}
        We now define, for all $h\in \mathcal{QP}$, 
        the language $\mathcal{L}_{2,h}$  as follows. 

	If every geodesic segment
        $[1,h]$ has no angle greater than $\Theta$, 
        or is of length $2$, then $\mathcal{L}_{2,h} =
        \emptyset$.

        Now assume that some 
	geodesic segment $[1,h]$ is of length $\geq 3$ and that, for $v\in [1,h]$, 
	$\Ang_v[1,h] \geq \Theta$ (thus by definition of  $\mathcal{QP}$, $d(v,h)=1$).
        The stabilizer $Stab(v)$  is then $hG_ih^{-1}$ for some $i$ (and for the given $h$).

	 We say
        that  $\tilde{g} \in \mathcal{L}$ is in $\mathcal{L}_{2,h}$ if and only if
	its normal form is a product $w_1w_2w_3$ with $w_1$ and $w_3^{-1}$ defining the same element $h' \in H_2$ with the property that $d(h',v) =1$ and  
	$\Ang_v([v,1][v,h']) \leq \epsilon$ for some segments,	 and with $w_2$ the normal form of an element of   $\mathcal{L}_p(i)$, defining $h''$  (thus at distance one from $h^{-1}v$)  such that  $\Ang_v([h^{-1}v,1][h^{-1}v,h'']) \geq 100\delta+ 2\epsilon$.

\end{defn}

        \begin{lem}
          For all $h$, the language $\mathcal{L}_{2,h} \subset F$ 
          is normalized rational.
        \end{lem}
\begin{proof}

          Let $h$ be such that   $\mathcal{L}_{2,h} \neq \emptyset$, and $v$, the vertex as in the definition. 
	There are only finitely many $h'\in H_2$ adjacent to $v$ such that  
        $\Ang_v([v,1][v,h']) \leq \epsilon $, 
         and, given $h'$,      
          there are only finitely many elements $w$ in 
          $\mathcal{L}$ representing them. Let $W(h')$ be this set, which is finite, and computable.

	On the other hand, only finitely many elements of  $\mathcal{L}_p(i)$ define an element $h''\in H_2$ such that  $\Ang_v([h^{-1}v,1][h^{-1}v,h'']) \leq 100\delta+ 2\epsilon$ (this follows from the conical stability of elements of $\mathcal{L}$, Proposition \ref{prop;stab_con}). In  $\mathcal{L}_p(i)$, the complement of this finite set is still normalized rational. Let $W'$ the set of normal forms of such elements.

          The set of  normal forms of elements of $\mathcal{L}_{2,h}$ is then the  elements of the concatenation of the sets  $W$, $W'$ and $W^{-1}$, that are reduced and ShortLex,  and therefore this set is a regular language. This makes  $\mathcal{L}_{2,h}$ normalized rational.

\qed \end{proof}

        We define, for all $h\in \mathcal{QP}$, the normalized rational
        language $\mathcal{L}_h= \mathcal{L}_{1,h} \cup  \mathcal{L}_{2,h}$.

        We will need two sufficient conditions to be in $\mathcal{L}_h$.   

\begin{defn}	A \emph{relative geodesic} between two vertices of $Cay(H_2)$ is a path in $Cay(H_2)$ such that each maximal subpath in a  
	coset of a parabolic subgroup is a geodesic in this coset, and, if each such subpath is replaced by the two edges of $\widehat{Cay}(H_2)$ passing through the infinite 
	valence vertex of this coset, one gets a geodesic of  $\widehat{Cay}(H_2)$. 
\end{defn}

Recall that if  $\Ang_v[1,h]\geq \Theta$, 
	the vertex $v$ is of infinite valence, hence is the vertex of some coset of a parabolic group.

\begin{figure}[hbt]
            \begin{center}
              \includegraphics[width=5in]{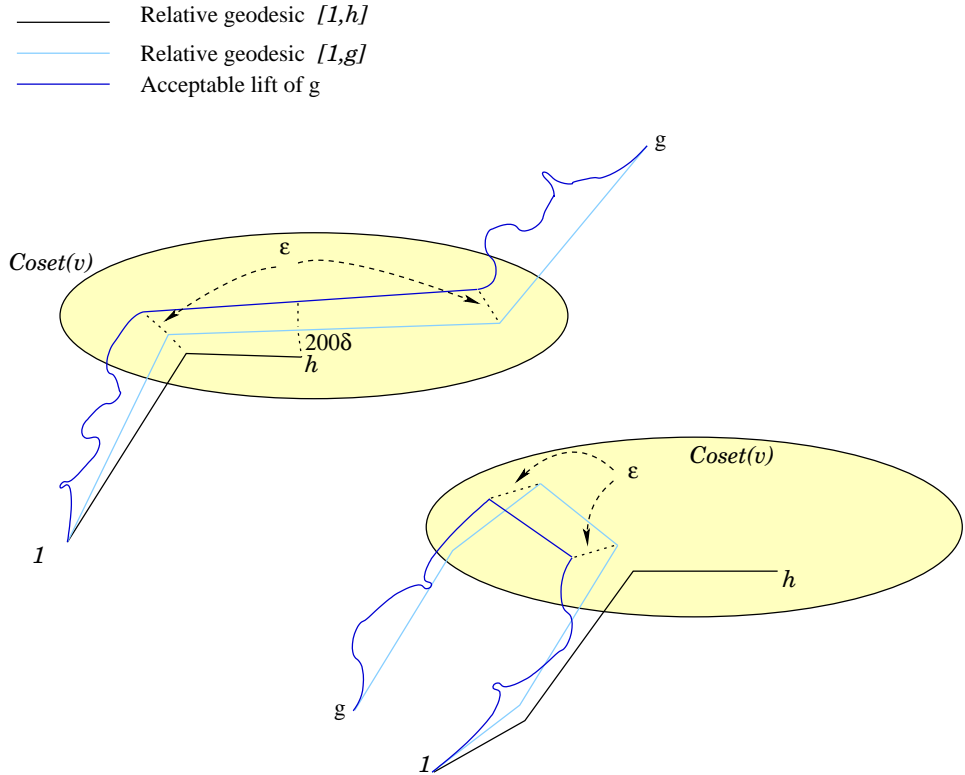}
              \caption{Sufficient conditions to be in, respectively,
               $\mathcal{L}_{1,h}$, and $\mathcal{L}_{2,h}$.}
                 \label{fig;44}
            \end{center}
          \end{figure}

        \begin{lem}\label{lem;suff_Lh}
         Let $h \in \mathcal{QP}$, and assume that there is a 
	vertex $v\in [1,h]$ adjacent to $h$   such that  
	$\Ang_v[1,h]\geq \Theta$. 
   
        Assume that $v\in [1,g]$, and that some $g' \in Coset (v)$ is on a relative geodesic from $1$ to $g$  (thus $g'$ is a vertex adjacent to $v$ in $\widehat{Cay}(H_2)$)  with  
        $\Ang_v([v,g'],[v,h]) \leq 200\delta$. 
        
        Then for all $\tilde{g}$ representing $g$ in $\mathcal{L}$, 
        $h$ is a quasi-prefix of $\tilde{g}$.

        \end{lem}
\begin{proof}
        One has $\Ang_v[1,g] \geq \Theta - 200\delta > 4\epsilon$. 
        By property of conical stability (Proposition \ref{prop;stab_con}, last point),  
        the element $\tilde{g}$ 
        labels a path in 
        $\widehat{Cay}(H_2)$ that contains the vertex  $v$, and a a shortest word defining $\tilde{g} \in F$
   contains a      subword $p$ that defines a path in the coset $Coset(v)$
        starting and ending at distance at most $\epsilon$ from the entering point
        and the exiting  point of the relative geodesic 
		from $1$ to $g$ in
        that coset. 
        Therefore, some shortest word defining $\tilde{g}$ labels
        a path in $Cay (H_2)$ that contains a vertex $g''$ at distance at most
        $3\epsilon $                              
        from  $g'$ in the coset $Coset(v)$. Therefore, it is at
        distance (in the Cayley graph) at most $3\epsilon +200\delta$ 
        from $h$. This ensures $\Ang_v[h, g'']\leq \eta$, and therefore,
        $h$ is a quasi-prefix of $\tilde{g}$. For an illustration, see Figure
        \ref{fig;44},  left side.
\qed \end{proof}

        \begin{lem} \label{lem;suff_prime}
        Let $h\in \mathcal{QP}$ such that  $d(1,h)\geq 3$, and  with $v\in[1,h]$ adjacent to $h$, 
        at which $\Ang_v[1,h]\geq \Theta$.

        Let $g$ be in $H_2$. If $v\in [1,g]$, if $g$ is in the stabilizer  $Stab(v)$ of $v$, and if
        $\Ang_v [1,g] > \Theta$,   
        then any representative in $F$ of $g$ is in $ \mathcal{L}_{2,h}$. 
        \end{lem}     
\begin{proof}        
        It follows from the property of conical stability (Proposition \ref{prop;stab_con}, last point) 
	that any representative in
        $\mathcal{L}$ of $g$ can be written as the concatenation of three
        elements of $\mathcal{L}$, $w_1, w_2, w_3$, such that $w_1$ and
        $(w_3)^{-1}$ each represent an element adjacent to $v$, respectively
        $h_1, h_3$  such that $\Ang_v[h_i,v][v,1] \leq \epsilon$. Therefore,
        $w_2$ is an element of $\mathcal{L}_p$ that produces an angle at least
        $\Theta - 2 \epsilon > 100\delta +2\epsilon$ at $v$. Hence, $w_1w_2w_3 \in
        \mathcal{L}_{2,h}$. For an illustration, see Figure
        \ref{fig;44},  right side. 
        \qed \end{proof}

    % FOR DANIEL:  
  
 % FOR FRANCOIS:
%\begin{figure}[hbt]
%            \begin{center}
%              \includegraphics{44.eps}
%              \caption{Sufficient conditions to be in, respectively,
%              $\mathcal{L}_{1,h}$, and $\mathcal{L}_{2,h}$.}
%                \label{fig;44}
%           \end{center}
%         \end{figure}

    \subsectionind{A finiteness result}\label{part;finiteness_rh}

       \subsubsection{A remark of I.~Bumagin:  the cascade effect}      
In her solution to the conjugacy problem in relatively hyperbolic groups, 
Inna Bumagin points out an interesting phenomenon, 
that she calls the `cascade effect'. This phenomenon can happen in the
situations we have to study, and we need a slight variation 
of the control she gets. The idea presented in this paragraph is already in \cite[\S 5.1]{Bumagin}.

Let $H$ be a relatively hyperbolic group, and $g, h$  elements of $H$. We
consider a quadrilateral of vertices $1,h,hg,hgh^{-1}$ in the coned off Cayley graph 
(we make the choice of a geodesic segment for each side). 
We say that the \emph{cascade effect} occurs if there is a vertex $w_0$ of
infinite valence, on $[1,hgh^{-1}]\cap [hgh^{-1}, hg]$, and an integer $n>0$ such that
the vertex $w_i = hg^{-i}h^{-1} w_0$ is in  $[1,h]\cap [hgh^{-1}, hg]$ for all $i=1, \dots, n$. In
such case, we choose 
$n$ maximal for this property.
 If the cascade effect occurs, since $w_n \in [hgh^{-1}, hg]$ and $n$ is maximal, the 
 vertex $w_{n+1} = hg^{-1}h^{-1} w_n$ is  in $[1,h]$,
but not on $ [hgh^{-1}, hg]$.

Understanding the cascade effect gives  control over certain angles.

\begin{lem} \label{lem;cascade}
        For all $g$, there exists a constant $Casc(g)$, which may be chosen to be greater than $\MaxAng[1,g]$,  
        such that, for all
        $h$,  if the cascade effect occurs in $(1,h,hg,hgh^{-1})$,  
        for a vertex $w_0 \in [1,hgh^{-1}]\cap [hgh^{-1}, hg] $  
        such that $hg^{-1}h^{-1} w_0 \notin [1,hgh^{-1}]$, 
        then $\Ang_{w_0}[hgh^{-1},hg] \leq Casc(g)$.
\end{lem}
 
The lemma can be compared with Bumagin's  \cite[Lemma 5.8]{Bumagin}. We first
define some notation.

        The segment $[1,h]$ contains all the vertices $w_1, \dots, w_{n+1}$. 
	It enters $Coset(w_i)$ in $x_i \in H_2$ and exits it in $y_i\in H_2$. 
	Let $p_i$ be the parabolic element  $p_i=x_i^{-1} y_i$ (in some $G_j$).  
	Let $h_i=y_{i-1}^{-1} x_{i}$ for $i\neq 1, n+2$, $h_1= x_1$, and $h_{n+2}= y_{n+1}^{-1}h$.    

        One has $h=h_1 p_1 h_2 p_2 \dots h_{n+1} p_{n+1} h_{n+2}$.

\begin{figure}[hbt]
                \begin{center}
              \includegraphics[width=5in]{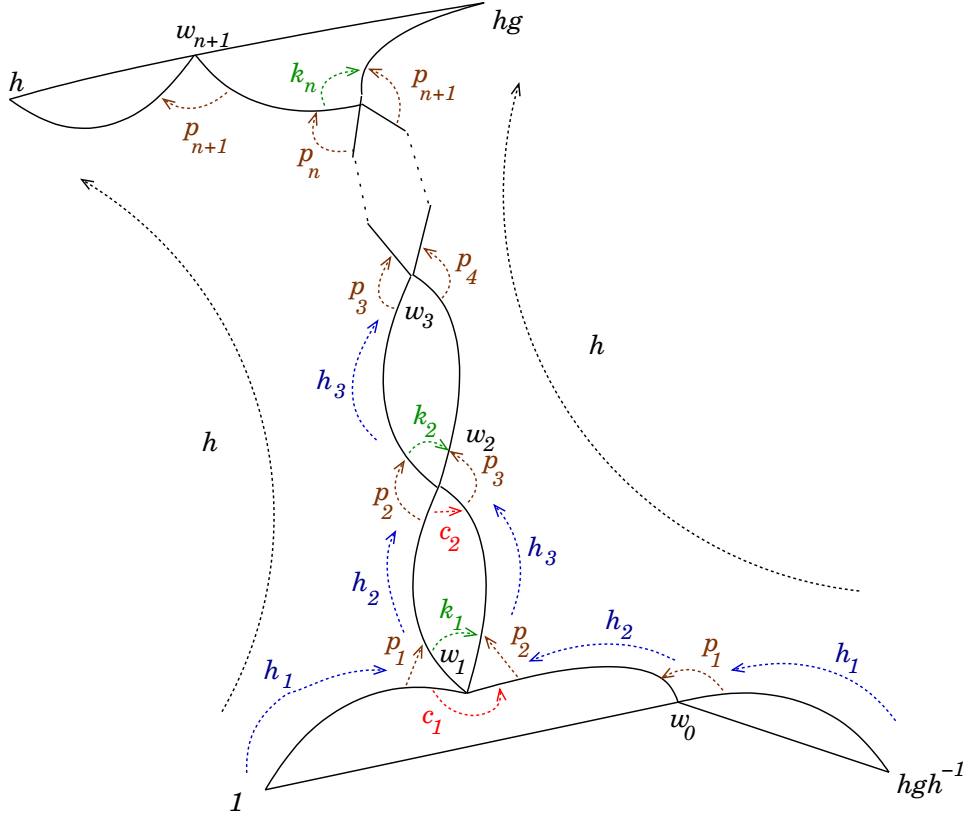}
              \caption{The notations in the cascade.}
              \label{fig;cascade}
            \end{center}
          \end{figure}

        The segment $[hgh^{-1},hg] = hgh^{-1} [1,h]$ contains $w_0$, 
        and all the
        vertices $w_1, \dots, w_n$ (but not $w_{n+1}$). 
	The parabolic elements that are difference between exiting point and  entering point 
	of $[hgh^{-1},hg]$ in $Coset(w_i)$ are the elements $p_{i+1}$ 
        (note the shift in the index). For all $i$ between $1$ and $n$, we
        note respectively $c_i$ and $k_i$ the parabolic elements that join
        respectively the
        entering points, and the exiting points, of $[1,h]$ and $[hgh^{-1}
        hg]$ in $Coset(w_i)$. These notations are illustrated in 
        Figure \ref{fig;cascade}.

        The aim of the lemma is to bound the size of $p_1$ in terms of constants
        computable from the path $[h,hg]$.

\begin{proof}[of Lemma \ref{lem;cascade}]
  First we ``collapse'' all the intermediate levels of
        the cascade. More precisely, 
        for all $i\leq n$, we defined $k_i=p_i^{-1} c_i
        p_{i+1}$. Hence, one has 
\[  k_1k_2\dots k_n = p_1^{-1} c_1c_2\dots p_{n+1}.  \]

        This allows us to express $p_1$ as $p_1=  (k_1\dots k_n)^{-1} (c_1\dots
        c_n) p_{n+1}$.

        Second, we translate all the levels of the cascade between 
        $w_1$ and $w_2$.  More precisely,    
        for all $i\leq n-1$, one has $k_i=h_{i+1} c_{i+1}
        h_{i+2}^{-1}$, and  
       \[  k_1k_2\dots k_{n-1} = h_2 c_2c_3\dots c_n h_{n+1}^{-1}.   \]

        Consider then the four vertices  $(y_1), (y_1 h_2), (y_1 h_2c_2c_3\dots c_n),  (y_1 k_1\dots k_{n-1})$. 
	The vertices   $y_1$ and $ (y_1 k_1\dots k_{n-1})$ are adjacent to $w_1$, and the two others are adjacent to $w_2$.
	A geodesic $[y_1, y_1 h_2  ]$ does not contain $w_1$ nor $w_2$, by definition of $h_2$.

 On the other hand,  $[(y_1 k_1\dots k_{n-1}),y_1 k_1\dots k_{n-1} h_{n+1}]  = \gamma[y_n, y_{n}h_{n+1}]$ for $\gamma$ such that $\gamma w_n = w_1$ and $\gamma w_{n+1} = w_2$ (because $y_1 k_1\dots k_{n-1} h_{n+1}= y_1 h_2c_2c_3\dots c_n$). By definition of $h_{n+1}$,  neither $w_n$ nor $w_{n+1}$ is on $[y_n, y_{n+1}h_{n+1}]$, hence  the segment $[(y_1 k_1\dots k_{n-1}),y_1 k_1\dots k_{n-1} h_{n+1}]$   does not contain $w_1$ nor $w_2$. Therefore, $\Ang_{w_1}[y_1, y_1 k_1\dots k_{n-1}] \leq 50\delta$, since one gets a path from one to the other avoiding $w_1$ by thinness of the bigone $(w_1,w_2)$. Similarly, one gets     $\Ang_{w_2}[y_1 h_2, (y_1 h_2c_2c_3\dots c_n)] \leq 50\delta$. We therefore get a computable bound (depending only on $\delta$ and $\widehat{Cay}(H_2)$) on the length of shortest words defining  $k_1\dots k_{n-1}$ $c_2c_3\dots c_n$ in the basis of the parabolic subgroup.

        Thus, to bound $p_1$, 
        it suffices to get a bound on the word length of $k_n$, $c_1$
        and $p_{n+1}$. One easily gets that the angle produced by $k_n$ is
        bounded by $\MaxAng[h,hg] + 50\delta$, and thus $k_n$ is bounded by a
        constant depending only on $g$. Since $w_{n+1} \notin[hgh^{-1}, hg]$,
        one also can bound $\Ang_{w_{n+1}} \leq \MaxAng[h,hg] +50\delta$,
        hence a bound on $p_{n+1}$ in terms of $g$ only.

        Finally, to bound $c_1$, we use the hypothesis of the lemma: 
        $w_1 \notin [1,hgh^{-1}]$, and therefore, the
        angle produced by $c_1$ is at most $50\delta$. This gives a universal 
        bound on the word length of $c_1$. 
        \qed \end{proof}

        % FOR DANIEL

% FOR FRANCOIS
%    \begin{figure}[hbt]
%                \begin{center}
%              \includegraphics{cascade.eps}
%              \caption{The notations in the cascade.}
%              \label{fig;cascade}
%            \end{center}
%          \end{figure}

       \subsubsection{Production of quasi-prefixes} 
When dealing with hyperbolic groups in Section \ref{List}, it
was sufficient to take {\em any} quasi-prefix of an element.  However,
in the setting of toral relatively hyperbolic groups, we need to be a little
more careful.  In particular, if $h$ has a long subsegment in a single
parabolic coset, we need to choose the quasi-prefix so that, in some sense, it
testifies that there is this longer subsegment. This is achieved with the following definition.

\begin{defn} \label{def:prefix}
                 Given an element $h\in H_2$, 
                 we say that $h_0$ is a \emph{prefix} of $h$ if there
                is a vertex $v$ in a geodesic  $[1,h]$ of $\widehat{Cay}(H_2)$, 
                and elements $h'$ and $h''$ neighbors of $v$ in $[1,h]$, 
                such that, $h_0=h'p$, 
                with $p\in Stab(h'^{-1} v)$,  
so that                if $v$ is of infinite valence then:
                in $Stab(h'^{-1}v)$ (with the word metric given by a chosen
                 basis) seen as the lattice of integers in
                 $\mathbb{R}^n$, $p$ is a point of the lattice at distance at most $2$ from
                 the segment $[1,h'^{-1}h''] \subset \mathbb{R}^n$. 
                 We say, in this case, that $p$ is on a \emph{diagonal} 
                 from $1$ to
                 $h'^{-1}h''$ 
                 in the abelian group  $Stab(h'^{-1}v)$.
\end{defn}
Note that if the vertex $v$ in Definition \ref{def:prefix} above is {\em
not} infinite valence, then $Stab(h'^{-1} v) = \{ 1\}$ and $h_0 = h'$.

         \begin{lem}\label{lem;prefix_rh}             
           For all $g\in H_2$, 
           there is a constant $K_{pre}(g)$ such that if the coset $hCent(g)$
           is at word distance at least 
           $K_{pre}(g)$ from $1$, 
           and if $h_0 \in \mathcal{QP}$ is a prefix of $h$, 
           then any representative of 
           $hgh^{-1}$ in $\mathcal{L}\subset F$ is in 
           $\mathcal{L}_{h_0}$.  
         \end{lem}
   
         The proof of this lemma will require some technicality: we will
         distinguish between 5 main steps, though for hyperbolic groups, only
         the first one is necessary.

        Let $K(g)= 2\rho+d(1,g)+200\delta$, 
        and $K'(g)=3(3+\Theta)(Casc(g) + 100\delta)$.

\begin{proof}[of Lemma \ref{lem;prefix_rh}]
         If $hgh^{-1}=h'gh'^{-1}$, then $h'^{-1}h \in Cent(g)$ and so
         $h'Cent(g) = hCent(g)$. 
        Thus, given   $K(g)$ and $K'(g)$, 
                there exists a constant $K_{pre}(g)$, such that  
                if $dist(1,hCent(g))\geq K_{pre}(g)$, 
                then  either $d(1,hgh^{-1})>K(g)$  
           or there is $v$ such that 
           $\Ang_{v} [1,hgh^{-1}] > K'(g)$.

           Moreover, since the parabolic subgroups of 
           $H_2$ are abelian and malnormal, there
           is at most one conjugate of $g$ in each of them. We choose
           $K_{pre}(g)$ large enough to ensure that $hgh^{-1}$ 
           is not in the
           stabilizer of a vertex adjacent to $1$ in $\widehat{Cay}(H_2)$ 
           (there
           are only finitely many).

\vskip .5cm

           {\sc Step 1: }  If $d(1,hgh^{-1})> K(g) $,
           then  any representative in 
           $\mathcal{L}$ of $hgh^{-1}$ is in $\mathcal{L}_{1,h_0}$.

\vskip .5cm
            
           This is similar to the hyperbolic case. 
           We  prove the contrapositive: that 
           if  $h_0$ is not a quasi-prefix of some acceptable lift of 
           $hgh^{-1}$ (in fact, it is enough to assume that $h_0
           \notin ConN_{100\delta,100\delta}([1,hgh^{-1}])$),  then
           $d(1,hgh^{-1})\leq  K(g)$ (and similarly for  $hgh^{-1} h_0$).

                 By hyperbolicity in the quadrilateral
           $(1,h,hg,hgh^{-1})$, the vertex $h_0$ is in the conical
           neighborhood 
           $ConN_{100\delta, 100\delta}( [h,hg] \cup  [hg,hgh^{-1}])$. The
           computation is then identical to the one of Lemma
           \ref{lem;prefix}. We reproduce it.

            First assume that there is an edge $e\subset [h,hg]$ for which 
            $h_0 \in \Cone_{100\delta,100\delta}(e)$, 
            and $v\in e$ at distance at most
            $100\delta$  from $h_0$. 
            Then we bound the
            distances $d(1,h) \leq
           d(1,h_0)+d(h_0,v) + d(v,h)$ and also 
           $d(1,hgh^{-1}) \leq d(1,h_0) +
           d(h_0,v) + d(v,hg)+d(hg,hgh^{-1})$.    
           Using $d(1,h_0)  \leq \rho$ and  $d(hg,hgh^{-1}) =d(1,h)$, 
           one obtains 
           $d(1,hgh^{-1}) \leq 2\rho +200\delta + d(v,hg)+d(v,h) = 2\rho
           +200\delta +d(1,g) \leq K(g)$.

              Similarly, if   there is an edge $e\subset [hg,hgh^{-1}]$ with  
              $h_0 \in \Cone_{100\delta,100\delta}(e)$. Let $v\in e$  
              be at distance at most
            $100\delta$  from $h_0$. 
            One has the bound 
           $d(1,hgh^{-1}) \leq
           d(1,h_0) + d(h_0,v)+ d(v,hgh^{-1}) \leq \rho
           + 100\delta + d(v,hgh^{-1})$, 
           and also $d(1,h) \leq d(1,h_0) +
           d(h_0,v)+d(v,hg) +d(hg,h)$, giving 
           $d(1,h) \leq \rho + 100\delta +
           d(v,hg) + d(1,g)$. Since $d(1,h)=d(hgh^{-1},
           hg)=d(v,hgh^{-1}) + d(v,hg)$, one deduces 
           $d(v,hgh^{-1})\leq \rho + 100\delta  + d(1,g)$. 
           Together with the
           first bound obtained, this gives $d(1,hgh^{-1}) \leq  2\rho
           + 200\delta +  d(1,g)\leq K(g)$. Finally, the condition that quasi-prefix and suffix are disjoint is ensured by the fact that $hgh^{-1}$ is a sufficiently long element for the metric word $dist$ (this remark indeed applies for all the subsequent steps).  Hence  we proved the claim of Step 1. (For an illustration of the situation, see Figure  \ref{fig;claim12}, left side)

% FOR DANIEL
    \begin{figure}[hbt] 
          \begin{center}
              \includegraphics[width=5in]{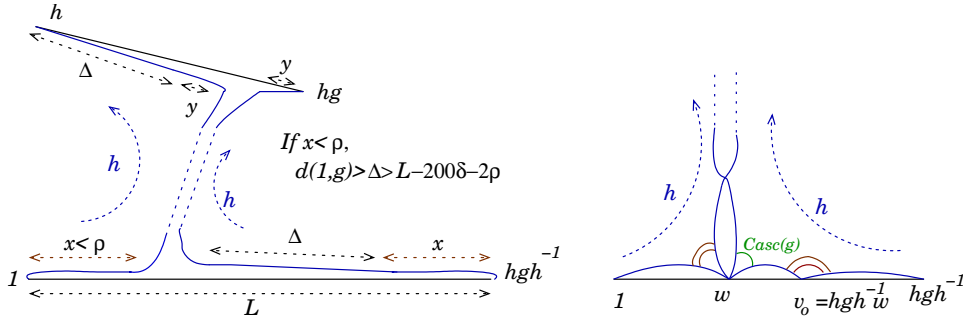}
              \caption{The situation in the first (left), and second (right) steps of Lemma \ref{lem;prefix_rh}.}
              \label{fig;claim12}
            \end{center}
          \end{figure}

% FOR FRANCOIS  
%    \begin{figure}[hbt] 
%            \begin{center}
%              \includegraphics{Claim1.eps}
%              \caption{The situation in the first (left), and second (right) steps of Lemma \ref{lem;prefix_rh}.}
%              \label{fig;claim12}
%            \end{center}
%          \end{figure}

           We now assume for the rest of the lemma that
           $d(1,hgh^{-1})\leq K(g)$, and we let $v_0$ be the vertex 
           closest to $1$ on $[1,hgh^{-1}]$ such that 
           $\Ang_{v_0} [1,hgh^{-1}] > K'(g)$.

        Let $K''(g) = (2 + \Theta)(Casc(g) +  100\delta) $. In the following claim, we will use that  
        $(K'(g)-\MaxAng[1,g]-50\delta)/2  \geq  K''(g) + 50\delta + Casc(g)>\Theta$, which is easily verified.

\vskip .5cm

           {\sc Step 2: } We claim that there are vertices $w\in [1,v_0]$ and $w' \in [v_0,
           hgh^{-1}]$ such that
           $\Ang_w[1,hgh^{-1}] \geq K''(g)    $ and 
        $\Ang_w[1,h] > \Theta$, and also 
        $\Ang_{w'}[1,hgh^{-1}] \geq K''(g) $ and 
           $\Ang_{w'}[hgh^{-1},hg] > \Theta$.

\vskip .5cm

           We prove the existence of $w$ first.
           In the quadrilateral $(1,h,hg,hgh^{-1})$, the triangle 
           inequality for  angles at $v_0$ easily gives  
           that one of $\Ang_{v_0}[1,h] $ and $\Ang_{v_0}[hgh^{-1},hg]$ must be
           greater than 
           $(K'(g)-\MaxAng[1,g]-50\delta)/2 \geq K''(g) + 50\delta + Casc(g)>
           \Theta$. 
           If  $\Ang_{v_0}[1,h]>\Theta$, we can choose $w=v_0$. 
    
           If on the contrary   
           $\Ang_{v_0}[1,h]\leq\Theta$, then  $\Ang_{v_0}[hgh^{-1},hg]>
           Casc(g)$. 
           Then (the contrapositive of) 
           Lemma \ref{lem;cascade}, on cascade effect, ensures that 
           $hg^{-1}h^{-1} v_0 \in [1,hgh^{-1}]$. 
           We define $w$ to be this vertex.  
           By translation, we see that $\Ang_w[1,h] = \Ang_{v_0}[hgh^{-1},hg] \geq 
           K''(g)+50\delta+Casc(g) > \Theta$. 
           Lemma  \ref{lem;cascade} applied to $w$ gives
           the bound  $\Ang_w[hgh^{-1},hg]\leq Casc(g)$ (we understand this
           angle as $0$ if $w \notin [hgh^{-1},hg]$).
           The triangle inequality for angles then gives $\Ang_w[1,hgh^{-1}]
           \geq \Ang_w[1,h] - 50\delta - Casc(g)$. Thus, $\Ang_w[1,hgh^{-1}]
           \geq K''(g)$. Thus $w$ satisfies the required properties.  (See Figure  \ref{fig;claim12}, right side).
           Symmetrically, by the same argument, one gets $w'$. 

\vskip .5cm

           In  the following, $w$ and $w'$ are chosen to satisfy the statement of Step
           2, and be closest to, respectively, $1$ and $hgh^{-1}$. It is still
           possible that they are both equal (to $v_0$). In any case, the
           angles at $w$ and $w'$ 
           implies that they are 
           vertices of infinite valence. 
 
\vskip .5cm

           {\sc Step 3: } If $d(1,w)\geq d(1,h_0)$, and $d(w',hgh^{-1})\geq d(1,h_0)$, then 
      any representative of 
           $hgh^{-1}$ in $\mathcal{L}$ is in $\mathcal{L}_{1,h_0}$.

\vskip .5cm

        As a general remark, $h_0 \neq w$, since the latter is of infinite
        valence.

        Let us first assume that $\MaxAng[1,h_0]<\Theta$. Then $h_0 \in
        [1,h]$, as well as $w$.  
        Since $d(1,h_0)\leq d(1,w)$, and $w\neq h_0$, the inequality is strict, and 
        $h_0$ is in the  bigone $(1,w)$, and by hyperbolicity,  
        one obtains an edge $e$ of the subsegment 
        $[1,w]$ of $[1,hgh^{-1}]$  so that 
        $h_0 \in \Cone_{50\delta,50\delta}(e)$. 
         By property of elements of $\mathcal{L}$, 
        for all representative of $hgh^{-1}$ 
        in $\mathcal{L}$, there is an edge $e'$ 
        on the path it defines in $\widehat{Cay}(H_2)$ such that
        $\Cone_{\epsilon,\epsilon} (e')$ contains $e$. One deduces that         
        $h_0\in \Cone_{\epsilon+150\delta,\epsilon+150\delta} (e')$, and 
        therefore that $h_0$ is a quasi-prefix of any representative of $hgh^{-1}$ 
        in $\mathcal{L}$.

        Let us assume now that there is $v$ in $[1,h_0]$ 
        adjacent to $h_0$, $d(1,v) <\rho$, such that $\Ang_v[1,h] \geq\Theta$. 
        Since $d(1,v)<d(1,h_0)\leq d(1,w)$, $v$ is in the bigone $(1,w)$, and by conical thinness, it is 
        on $[1,hgh^{-1}]$.

        Let $[1,v]_h$ and $[v,w]_h$ be the subsegments of $[1,h]$ from $1$ to $v$ and
        from $v$ to $w$, respectively, and define $[1,v]_{hgh^{-1}}$ and $[v,w]_{hgh^{-1}}$
        analogously as subsegments of $[1,hgh^{-1}]$. We can bound the angles  
        $\Ang_{v}([v, 1]_h, [v,1]_{hgh^{-1}}) \leq 50\delta$, and $\Ang_{v}([v, w]_h, [v,w]_{hgh^{-1}}) \leq 50\delta$.
       Thus the diagonal  in $Coset(v)$ between the entering point of $[1,hgh^{-1}]$ and its exiting point, 
        remains at distance $150\delta$ 
        from 
        the diagonal between the entering and exiting point of $[1,h]$. 
        In particular it passes at distance $150\delta$ from $h_0$,  
        and by Lemma \ref{lem;suff_Lh}, $h_0$ 
                is a quasi-prefix of any representative in $\mathcal{L}$ of 
                $(hgh^{-1})$, in this case also.  

        Similarly, one gets, studying the bigone $(w',hgh^{-1})$ that $h_0$ 
                is a quasi-prefix of any representative in $\mathcal{L}$ of 
                $(hgh^{-1})^{-1}$. As in in the end of Step 1, the fact that quasi-prefix and suffix are disjoint is ensured by the fact that $hgh^{-1}$ is a sufficiently long element for $dist$.  This completes Step 3.

\vskip .5cm

            We now assume that either $d(1,w) <d(1,h_0)$ or 
            $d(w',hgh^{-1}) <d(1,h_0)$.   By symmetry, one can assume, without
            loss of generality, that $d(1,w) <d(1,h_0)$. 
        Since $w\in [1,h]$ and $\Ang_w[1,h]>\Theta$, by the definition of $\mathcal{QP}$, 
        $w$ is adjacent to $h_0$, and $\Ang_w[1,h_0] \geq \Theta $.

\vskip .5cm            

              {\sc Step 4: }  
        If  $hgh^{-1} \in Stab(w)$, 
              then any
              representative in $\mathcal{L}$ of $hgh^{-1}$ 
              is in $\mathcal{L}_{2,h_0}$.

\vskip .5cm

              Recall that, by the choice of $K_{pre}(g)$,  $d(1,w) = d(w,hgh^{-1})=d(1,hgh^{-1})/2  >3$.   
        Thus, since   $\Ang_w[1,hgh^{-1}] \geq K''(g)> \Theta$, and  $\Ang_w[1,h_0] \geq \Theta$,
              Lemma \ref{lem;suff_prime} applies: 
              any representative in $\mathcal{L}$ 
              of $hgh^{-1}$ is in   $ \mathcal{L}_{2,h_0}$. 

\vskip .5cm

     From now on, we assume that 
     $w \neq hgh^{-1} w$.  
        Therefore, the vertices $1$,$w$,$hgh^{-1}w$, and $hgh^{-1}$ are all on $[1,hgh^{-1}]$, and appear in that order. 
                  We will now use that 
        $K''(g) -  Casc(g)-200\delta \geq \Theta(Casc(g)+100\delta) $, which follows from its definition.

\vskip .5cm

     {\sc Final step: }  In this setting ({\it i.e. } if $w \neq hgh^{-1} w$),  any representative in
        $\mathcal{L}$ of  $hgh^{-1}$ is in $\mathcal{L}_{1,h_0}$.

\vskip .5cm

        We first need to prove that $h_0$ is a quasi-prefix of every representative in  
        $\mathcal{L}$ of  $hgh^{-1}$.

        With the convention that $\Ang_w[x,z] =0$ when $w\notin [x,y]$, one has 
        $\Ang_w[h,gh]\leq \MaxAng[h,gh] \leq \MaxAng[1,g] \leq Casc(g)$. 
        We also need a bound on $\Ang_w[hg,hgh^{-1}]$: if $w\in[hg,hgh^{-1}]$, 
        since $w \neq hgh^{-1} w$, one checks that 
        $hg^{-1}h^{-1} w  \notin [1,hgh^{-1}]$, and Lemma \ref{lem;cascade} provides  
        $\Ang_w[hg,hgh^{-1}]\leq Casc(g)$. Since, on the other hand,   
        $\Ang_w[1,hgh^{-1}] \geq K''(g)$, by the triangle inequality for angles at $w$, 
        $\Ang_w[1,h] \geq K''(g) - 2 Casc(g) - 200\delta$ (see Figure \ref{fig;step6}).

        Now, in $Coset(w)$ one considers the diagonal between 
                the entering and exiting points of 
                $[1,h]$, and the diagonal between the entering and exiting 
                point of $[1,hgh^{-1}]$. 
                The vertex $h_0$ is on the first one, at distance (in  $\widehat{Cay}(H_2) \setminus \{w\}$) 
        at most 
                $\Theta$ from the entering point. 
                The two entering points are at most $50\delta$ apart 
                (in  $\widehat{Cay}(H_2) \setminus \{w\}$), and the two 
                exiting points are at most 
                $2(Casc(g)  +100\delta)$ apart. 
                Since the length of the second diagonal is greater than 
                $K''(g) -2(Casc (g)+ 100\delta) \geq \Theta \times 
                (Casc(g)  +100\delta)$,    
                the two diagonals remain at distance at most $100\delta$  
                in $\widehat{Cay}(H_2) \setminus \{w\}$, 
                for a length of at least $\Theta$. 
                Therefore, the second one has a point at distance at most 
                $100\delta$ from $h_0$, and by 
                Lemma \ref{lem;suff_Lh}, $h_0$ is a quasi-prefix of any 
                representative in 
                $\mathcal{L}$ of $hgh^{-1}$.

% FOR DANIEL
    \begin{figure}[hbt] 
            \begin{center}
              \includegraphics[width=5in]{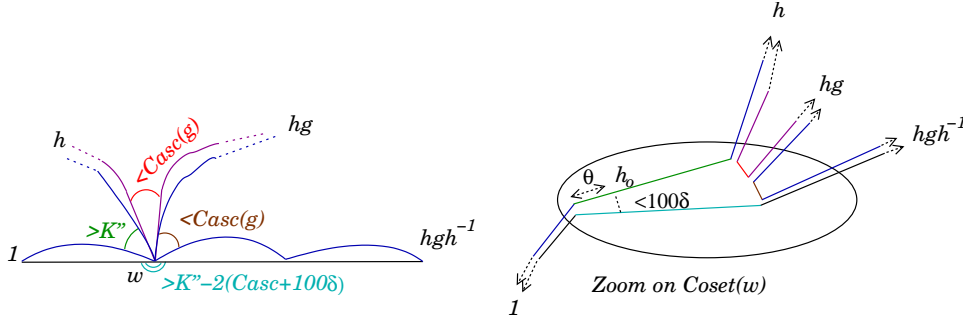}
              \caption{Angles and diagonals at the vertex $w$.}
              \label{fig;step6}
            \end{center}
          \end{figure}

% FOR FRANCOIS
%    \begin{figure}[hbt] 
%            \begin{center}
%              \includegraphics{Claim67.eps}
%              \caption{Angles and diagonals at the vertex $w$.}
%              \label{fig;step6}
%            \end{center}
%          \end{figure}

        We now need to prove that $h_0$ is a quasi-prefix of every representative in  
        $\mathcal{L}$ of  $hg^{-1}h^{-1}$.
        For that, we distinguish two cases. If $w\in [hgh^{-1}, hg]$, then $hgh^{-1}w$ 
                is in the bigone $(w,hgh^{-1})$, and the argument of Step 3 applies.
        If $w \notin [hgh^{-1}, hg]$, then $\Ang_w[1,h] 
                \geq \Ang_w[1,hgh^{-1}] - \MaxAng[h,hg] - 100\delta \geq K''(g)-Casc(g) -100\delta$. 
                        By translation, 
                $\Ang_{hgh^{-1}w}[hgh^{-1}, hg]\geq K''(g) - Casc(g) - 100\delta$.   
        
        We give bounds on angles at $hgh^{-1}w$:  $\Ang_{hgh^{-1}w}[h,hg] 
        \leq \MaxAng[1,g] \leq Casc(g)$, and by Lemma \ref{lem;cascade}, 
                $\Ang_{hgh^{-1}w}[1,h] \leq Casc(g)$. Therefore, by translation, 
                $\Ang_{hgh^{-1}w} [1,hgh^{-1}]\geq K''(g) - 3(Casc(g) + 100\delta)$.
        Now, by the same argument on diagonals in $Coset(hgh^{-1} w)$, as for $w$, one gets the claim. 
        
        This finishes the proof of the  lemma. 
        \qed \end{proof}

         \begin{cor} \label{coro;fini_rh}
           Given a compatible homomorphism $\phi: H_1 \to H_2$, 
           only finitely many conjugates of $\phi$ have an acceptable lift 
           satisfying $\Omega$.   

         \end{cor}

         One can reproduce {\it verbatim } the proof of Corollary
         \ref{coro;fini} from paragraph \ref{para;proof_hyp}, 
         with in mind that the word distance noted $d$ in
         the proof should be considered as the word 
         distance $dist$ here, and that calls to Lemma \ref{lem;prefix}  should be replace by calls to Lemma \ref{lem;prefix_rh}.

    \subsectionind{An existence result}\label{part;exists_rh}

      \subsubsection{The quantity $Q$}
        For an element
        $h\in H_2$, we recall that 
        $d(1,h)$ is the distance in the coned-off Cayley
        graph $\widehat{Cay}(H_2)$ between 
        the origin and the vertex associated to $h$. 
        We denote by
        $\theta(h)$ the minimum over all the geodesics of $\widehat{Cay}(H_2)$ 
        from $1$ to $h$ of the sum of the first and last angles between edges
        of the path. We endow pairs $\mu(h)= (d(1,h), \theta(h))$ with the
        lexicographical order. 
        Then, given a homomorphism $\phi: H_1 \to H_2$, 
        we set $Q(\phi) = \max \{ \mu(\phi(a)),  \mu(\phi(b)) \} $ (recall that $a$ and $b$ are two non-commuting elements of $H_1$ 
	chosen in advance.

        As a subset of $\mathbb{N}\times\mathbb{N}$ with lexicographical order,
        any   decreasing sequence is eventually stationary. In particular, one has the lemma:

        \begin{lem} \label{lem;mini}
          In every
          conjugacy class of homomorphism from $H_1$ to $H_2$, there is 
          a homomorphism realizing the minimum
          of $Q$ over the class.
        \end{lem}

	\subsubsection{Two lemmas for $\tilde{g} \in \mathcal{L}_{h}  $}
        Let $\eta'=\eta+\epsilon+100\delta$. Note that with this
        notation $\Theta \geq 
        2\eta' + 200\delta$.

        \begin{lem} \label{lem;pour_les_geod1_rh}
          Let $h\in \mathcal{QP}$, $\tilde{g} \in \mathcal{L}_{1,h}$, and
          $g$ its image in $H_2$. Let $[1,g]$ be a geodesic in
          $\widehat{Cay}(H_2)$. 

           Then there is a vertex  $w\in [1,g ]$ such that
          $d(w,h) \leq \eta'$.

          Moreover, if a geodesic segment $[1,h]$ contains a vertex $v$ with
          $\Ang_v[1,h] \geq \Theta$, then $v\in [1,g]$, $\Ang_v[1,g] \geq
          100\delta$ and
          $\Ang_v([v,h],[v,g]) < \Ang_v [1,g] - 100 \delta$.
        \end{lem}

\begin{proof}
          By definition of $\mathcal{L}_{1,h}$, there is a path  
	$p_{\tilde{g}}$ in $Cay(H_2)$ (not the coned-off : that is without the vertices of infinite valence) 	
		defined by  a shortest word representing 
	$\tilde{g}$, and that  
		contains an
          edge $e$ such that, in $\widehat{Cay}H_2$,  $h\in \Cone_{\eta,\eta}(e)$. 
      
	Now, the path 	$\hat{p}_{\tilde{g}}$ labeled by the normal form of  $\tilde{g}$ in $\widehat{Cay}H_2$ 
	(this time, possibly passing through vertices of infinite valence)
	 has a vertex at distance at most $1$ from $e$. This
          vertex $w_1$ is then at distance at most $\eta+1$ from $h$. Since $\tilde{g} \in \mathcal{L}$, the path labeled by its 
	normal form of    
           is a $(L'_1,L'_2)$-quasi-geodesic with end-points $1$ and $g$ (and containing $w_1$). 
          It is then contained in the
          $\epsilon$-neighborhood of $[1,g]$. Thus there is $w\in [1,g]$ at
          distance at most $\eta$ from $w_1$. One obtains $d(w,h)\leq
          \eta+\epsilon +1 < \eta'$, hence the first assertion.

          Assume   $[1,h]$ contains a vertex $v$ with
          $\Ang_v[1,h] \geq \Theta$.

		Let
          us consider a relative geodesic  
		from $1$ to $g$. 
          By the fellow
          traveling property for elements in $\mathcal{L}$ (see Proposition \ref{prop;stab_con}), 
          there is an edge
          $e'$ in this path such that $e\in \Cone_{\epsilon,\epsilon}(e')$,
          and therefore, $h\in \Cone_{\eta',\eta' }(e')$. 

          Let $v'$ be a vertex of $e'$. In the triangle $(v',h,1)$, the segment
          $[v',h]$ has maximal angle at most $\eta'$, whereas
          $\Ang_v[1,h]\geq \Theta$. This implies that $\Ang_v[1,v'] \geq
          \Theta - \eta'$. Since $v'$ is on the relative geodesic  
	from $1$ to
          $g$, the segment $[1,v']$ coincides with a segment $[1,g]$ up to the
          last edge of $[1,v']$. This ensures that $v\in[1,g]$, and that
          either $v'$
          is at distance $1$ from $v$, or $[1,g]$ coincides with $[1,v']$ until
          one edge after $v$. In this latter case, $\Ang_v[1,g] \geq \Theta -
          \eta'$, and $\Ang_v([h,v],[v,g]) \leq \Ang[1,g] -(\Theta + \eta')$,
          hence the result. 

          If now $d(v',v)=1$,  let $e_0$ be the first edge of a segment
          $[v,v']$.  
          Then 
          \[\Ang_v ([h,v],[v,g]) \leq \Ang_v([h,v],e_0) + 
          \Ang_v(e_0,[v,g]).\]
The first term is at most $\eta'$, and the
          second is at most $\Ang_v[1,g] - \Theta - \eta'$, thus
          $\Ang_v([h,v],[v,g]) \leq \Ang_v[1,g] - 200\delta$. Since the latter
          is positive,  $\Ang_v[1,g] \geq 200\delta$. 
          \qed \end{proof}

        \begin{lem} \label{lem;pour_les_geod2_rh}
           Let $h\in \mathcal{QP}$, $\tilde{g} \in \mathcal{L}_{2,h}$, and
          $g$ its image in $H_2$. Let $[1,g]$ be a geodesic in
          $\widehat{Cay}(H_2)$. 
          
          Assume that a geodesic segment $[1,h]$ contains a vertex $v$
          adjacent to $h$ with
          $\Ang_v[1,h] \geq \Theta$, then $v\in [1,g]$, $\Ang_v[1,g] \geq
          100\delta$, and $d(1,v)= d(v,g) \geq 2$.
        \end{lem}

\begin{proof}
        The element $\tilde{g}$ can be written $\tilde{g}=w_1ww_2^{-1}$, with $w_2$ representing the same element as $w_1$, and 
        $w_1$   and $w_1w$ representing  elements $h_1, h_2$  adjacent to $v$.
        and $\Ang_v[h_1,h_2] \geq 2\epsilon +100\delta$. 
	On the other hand, for a geodesic $[1,g]$, 
        one has that $\Ang_v[1,g] \geq \Ang_v[h_1,h_2] - \Ang_v([v,h_1 ]
        [v,1]) - \Ang_v([v,h_2 ] [v,g]) $. Therefore, $\Ang_v[1,g] \geq
        100\delta$. Since $\mathcal{L}_p\cap \mathcal{L}_{2,h} = \emptyset$,
        $d(1,v)= d(v,g) \geq 2$. 
        \qed \end{proof}

      \subsubsection{Possibilities of reduction}

        \begin{lem} \label{lem;reduc_reguliere_rh}
          Let $h\in\mathcal{QP}$
          Let $\tilde{g}$ be in 
          $\mathcal{L}_{h} \subset F$, and let
          $g\in H_2$ be its image in the quotient map $F \to H_2$.
          
          Then, either (i) $g \in \mathfrak{B}_{H_2}(2)$, the ball of radius $2$ about $1$ in $H_2$ for the distance $d$; or  
          (ii) $\mu(h^{-1}gh)< \mu(g)$.
        \end{lem}
\begin{proof}
          Let us assume that $g$ is not $\mathfrak{B}_{H_2}(2)$.

          First we assume that $\tilde{g} \in \mathcal{L}_{2,h}$. 
          Then by definition
          of $\mathcal{L}_{2,h}$, it is not in $\mathcal{L}_p$, and therefore,
          $d(1,g) \geq 3$. On the other hand, $h^{-1}gh$ is the image of an
          element of $\mathcal{L}_p$, by definition of  $\mathcal{L}_p$. 
	Therefore  $d(h^{-1}gh)=2$, hence the result.

          Now we assume that  $\tilde{g} \in \mathcal{L}_{1,h}$.  
          By Lemma \ref{lem;pour_les_geod1_rh}, there is a vertex $w\in [1,g]$
          at distance at most $\eta'$ from $h$, and similarly, there is $w'$
          on $[1,g]$  
          at distance at most $\eta'$ from $gh$.

          If $d(1,h) \geq \rho-1$, one computes  $d(h,gh)
          \leq 2\eta' + d(w,w') = 2\eta' +d(1,g) -d(1,w) - d(w',g)$.  Since
          $d(1,w)\geq d(1,h) -\eta'\geq \rho-1 - \eta'$, one gets $d(h,gh) \leq
          d(1,g) +3\eta' - \rho+1 < d(1,g)$.

          If on the contrary, $d(1,h)<\rho-1$, then there is $v$ adjacent to $h$
          such that $\Ang_v[1,h] \geq \Theta$. Then  $d(h,gh)\leq
          1+d(v,gv)+1$, which is at most $d(1,g)$, with equality only if
          $d(1,v)=1$, and both $v$ and $gv$ are on some geodesic $[h,gh]$. In
          this case, by
          Lemma \ref{lem;pour_les_geod1_rh}, one has $\Ang_v([v,h],[v,g]) < \Ang_v
          [1,g]$, therefore, $\Ang_v[h,gh] < \Ang_v
          [1,g]$, and similarly,  $\Ang_{gv}[h,gh] < \Ang_{gv}
          [1,g]$. This provides the inequality $\mu(h^{-1}gh)< \mu(g)$.
\qed \end{proof}

       We now state and prove the second main technical lemma of the section.

        \begin{lem}\label{lem;4cas_rh}
        Let $\phi:H_1 \to H_2$ be a homomorphism with an acceptable lift $\tilde{\phi}:
        \mathcal{B}_{H_1}(2) \to \mathcal{L}$ contradicting $\Omega$: there is $h\in
        \mathcal{QP}$ such that either both $\phi(a), \phi(b)$  are in $
        \mathcal{L}_{h} $, or all three $ \phi(b), \phi(ab), \phi(a^{-1}b)$ are in it, or all three  $ \phi(a), \phi(ba), \phi(b^{-1}a)$ are in it.
          
           Then $\mu(h^{-1} \phi
          (a) h) <Q(\phi)$ and similarly for $h^{-1} \phi (b) h$.

        \end{lem}
       \begin{proof}
		If $\phi(a) \in  \mathcal{L}_{h}$, then $\phi(a)$   is not in the $2$-ball of $H_2$. 
		Therefore  by Lemma \ref{lem;reduc_reguliere_rh}, 
           $\mu(h^{-1} \phi (a) h) <\mu(\phi (a)) \leq Q(\phi)$. Similarly for $\phi(b)$. 
		We now assume that  $\phi(a) \notin  \mathcal{L}_{h}$, and since $\Omega$ is 
		falsified by $\phi$, all three  $ \phi(b), \phi(ab), \phi(a^{-1}b)$ are in $  \mathcal{L}_{h} $.
    
          The discussion will be  about  the triangle 
          $(1,\phi(a), \phi(ab))$ in most cases except one for which it will be in    $(1,\phi(a^{-1}), \phi(a^{-1}b))$.

          We distinguish cases depending on whether $[1,h]$ has a final large angle or not.
          We first assume that a geodesic in $\widehat{Cay}(H_2)$ 
          from $1$ to $h$ has its last 
          angle greater than $\Theta$ at a vertex $v$ (adjacent to $h$),

          By Lemma \ref{lem;pour_les_geod1_rh}, if  $\varphi(ab)   \in \mathcal{L}_{1,h}$,    
        or Lemma \ref{lem;pour_les_geod2_rh}, if   
        $ \varphi(ab)  \in \mathcal{L}_{2,h}$, we have that  
        $\Ang_{v} [1,
          \phi(ab)]\geq 100\delta$, and similarly, 
          if $v_a=\phi(a)v$, 
          $v_a\in [\phi(a), \phi(ab)]$ and $\Ang_{v_a} [\phi(a),
          \phi(ab)]\geq 100\delta$.

          In particular, either $v$ is in 
          $[\phi(a),\phi(ab)]$  or in $[1,\phi(a)]$. We call this
          first situation case $(\alpha)$, and the second, case
          $(\beta)$. 
          
          The vertices $\phi(a),v,v_a,\phi(ab)$ are all on 
          the segment $[\phi(a),\phi(ab)]$. Let us consider case
          $(\alpha_1)$ when they appear in this order, and case 
          $(\alpha_2)$
          when $v$ and $v_a$ are inverted.
 
          We also subdivide case $(\beta)$ as follows. 
          Since $\Ang_{v_a} [\phi(a),
          \phi(ab)]\geq 100\delta$, either $v_a$ is in 
          $[1,\phi(ab)]$ or in $[1,\phi(a)]$. These two cases are called 
          $(\beta_1)$ and $(\beta_2)$, respectively. 
          (See Figure \ref{fig;4cas}).

% FOR DANIEL
          \begin{figure}[hbt]
                      \begin{center}
              \includegraphics[width=5in]{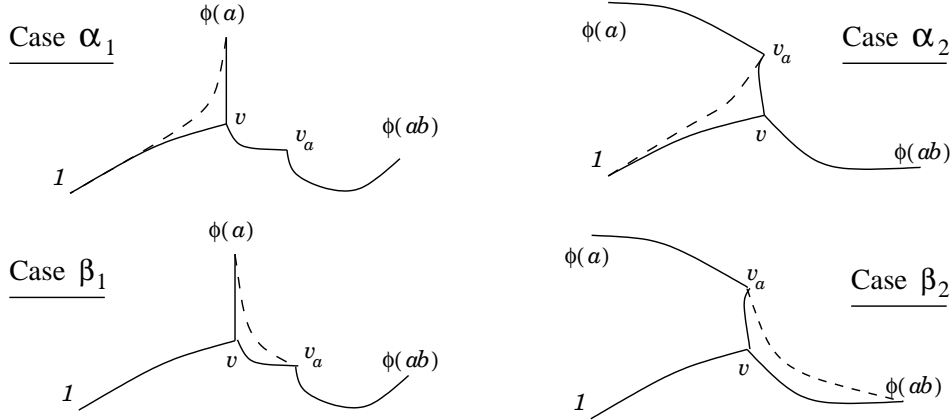}
              \caption{The first four cases of Lemma \ref{lem;4cas_rh}, when $[1,h]$
                has a large angle.}
                       \label{fig;4cas}
            \end{center}
          \end{figure}
  
  % FOR FRANCOIS
%          \begin{figure}[hbt]
%                      \begin{center}
%              \includegraphics{4cas.eps}
%              \caption{The first four cases of Lemma \ref{lem;4cas_rh}, when $[1,h]$
%                has a large angle.}
%                        \label{fig;4cas}
%            \end{center}
%          \end{figure}

          In all the cases, one has $d(1, h^{-1} \phi(a)h) =
          d(h,\phi(a)h) \leq 2+d(v,v_a)$.

          We now discuss the cases $\alpha_1$, and $\alpha_2$. In these cases,
          both $v$ and $v_a$ are on $[\phi(a),\phi(ab)]$, and  
          we compare $\mu(\phi^{h^{-1}}(a))$ to $\mu(\phi(b))$.

          In the case $\alpha_1$,
          $d(1,\phi(b)) = d(\phi(a),\phi(ab))=
          d(\phi(a),v)+d(v,v_a)+d(v_a,\phi(ab))$. This is strictly larger
          than $2+d(v,v_a) \geq d(1,\phi^{h^{-1}}(a))$ unless 
          $ d(\phi(a),v) = d(v_a,\phi(ab)) =1$, in which case, we have
          equality. 

          In this latter case, since $d(v_a,\phi(ab)) \geq
          d(\phi(a),v_a)$, we conclude that  $d(\phi(a),v)=
          d(\phi(a),v_a)=1$. Therefore, $v=v_a$ and $\phi(b) \in
          Stab(\phi(a)^{-1}v_a) = Stab(v)$, and all its representatives in
          $\mathcal{L}$  are in
          $\mathcal{L}_p$, which contradicts the fact that one of them 
          is in $\mathcal{L}_{h}$.

          Thus we conclude that, in case $(\alpha_1)$,  
          $\mu(h^{-1} \phi(a)h) < \mu(\phi(b))$.

         Case $(\alpha_2)$ is similar to the previous case, 
          the difference being that the vertices
          $v$ and $v_a$ are inverted on $[\phi(a),\phi(ab)]$.  
          We now have 
          \[d(1,\phi(b)) =
          d(\phi(a),v_a)+d(v_a,v)+d(v,\phi(ab)) \geq d(h,\phi(a)h),\]
          with equality only if $ d(\phi(a),v_a) = d(v,\phi(ab)) =1$.

          In this latter case, since  $ d(\phi(a),v_a)=d(1,v)$, the element
          $\phi(ab)$ is in $Stab(v)$, and we get the same contradiction
          as in case $(\alpha_1)$. Hence the same inequality  
          $\mu(h^{-1} \phi(a)h) < \mu(\phi(b))$.

          We now turn to case $(\beta)$.
          We first deal with $(\beta_2)$, in which the two vertices $v$ and $v_a$
          are on $[1,\phi(a)]$. We compare 
          $Q(h^{-1} \phi(a)h)$ to $Q(\phi(a))$.
          One has $d(1,\phi(a)) = d(1,v)+ d(v,v_a)+d(v_a,\phi(a))$, and this
          is larger than  $d(1,\phi^{h^{-1}}(a))$ unless 
          $ d(\phi(a),v_a) = d(1,v)  =1$.

      In this case, 
      first note that neither $ \varphi(ab) $, nor $ \varphi(b) $   
        is in $\mathcal{L}_{2,h}$.  They are therefore both in
        $\mathcal{L}_{1,h}$.

      If  $\Ang_v([v,\phi(a)],[v,\phi(ab)])\geq 50\delta$, 
      $v$ is in the third side of the triangle $(1,\phi(a),\phi(ab))$, 
      and we are in case $(\alpha_2)$. 

      If now  $\Ang_{v_a}([v_a,1],[v_a,\phi(ab)])\geq 50\delta$, then the triangle $(1,\phi(a^{-1}), \phi(a^{-1}b))$
	falls into   
      case $(\alpha_2)$. We get the conclusion that  $\mu(h^{-1} \phi (a^{-1}) h) < \mu(b)$.

      Otherwise, if both 
      \begin{eqnarray*}
      \Ang_v([v,\phi(a)],[v,\phi(ab)]) &<& 50\delta , \mbox{ and} \\ 
      \Ang_{v_a}([v_a,1],[v_a,\phi(ab)]) &<& 50\delta,
      \end{eqnarray*}
  then one can
      control $\Ang_v[h,  \phi(a) ] < \Ang_v[h,\phi(ab)] + 50\delta$,
      and by Lemma  \ref{lem;pour_les_geod1_rh},     
      \[        \Ang_v[h,  \phi(a) ] <\Ang_v[1,\phi(ab)] -100 \delta +50\delta.   \]
      Since $ \Ang_v[h,  \phi(ab) ] \leq  \Ang_v[1,\phi(a)] + 50\delta
      $, 
      one gets  
      \[ \Ang_v[h,\phi(a)h]  = \Ang_v[h,  \phi(a) ] <   \Ang_v[1,\phi(a)]
      , \] 
      and similarly   $\Ang_{v_a}[ \phi(a)h, h ] <
      \Ang_{v_a}[\phi(a),1] $. 
      Therefore, $\mu(h^{-1} \phi (a) h) < \mu(\phi(a))$.

      We end with case $(\beta_1)$, when it does not comply with case 
      $(\beta_2)$. In this situation, $v\in [1,\phi(a)]$ and 
      $v_a\notin [1,\phi(a)]$, and this means
      that $h$ is not a quasi-prefix for $\phi(a)^{-1}$. 
      In particular,   $\phi(a^{-1})$ is in situation $(\alpha_1)$ or  $(\alpha_2)$, for which
      we already have the result.

\vskip .5cm

          We now assume that no geodesic in $\widehat{Cay}(H_2)$ 
          from $1$ to $h$ has its last 
          angle greater than $\Theta$ (hence the vertex previously noted $v$ does not exist in this case). 
	  This situation is very similar to that
          which was studied in Subsection \ref{part;exists} when dealing with hyperbolic groups. 
          
          First,  since $h\in \mathcal{QP}$,
          $d(1,h) \geq \rho$.  Recall that in this case, $\mathcal{L}_{2,h}$ is
          empty, thus we can use Lemma \ref{lem;pour_les_geod1_rh}.

          By Lemma \ref{lem;pour_les_geod1_rh} there is 
          a vertex $v$ on a certain geodesic $[1,\phi(ab)]$ 
          such that $d(v,h)\leq \eta'$.

          Similarly, there is 
           a vertex $v_a$ lying on a geodesic
          $[\phi(a),\phi(ab)]$ such that  $d(v_a,\phi(a)h) \leq \eta'$.

          By hyperbolicity, $v$ lies at distance at most $\delta$ 
          from a vertex $v'$ 
          of $[1,\phi(a)] \cup [\phi(a),\phi(ab)]$.

           We again distinguish four cases. 
          The first dichotomy
          (cases $(\alpha)$ and $(\beta)$) 
          concerns whether $v'\in  [\phi(a),\phi(ab)]$ or 
          $v'\in  [1,\phi(ab)]$. 
          In case $(\alpha)$, both $v'$ and $v_a$ are in
          $[\phi(a),\phi(ab)]$, and we denote by $(\alpha_1)$ 
          the case where they
          appear on the segment in order  $(\phi(a), v',v_a, \phi(ab))$, 
          and by
          $(\alpha_2)$ the case when they appear in the order  
          $(\phi(a),v_a,v', \phi(ab))$.

          In case $(\beta)$,  $v'$ is on $[1,\phi(ab)]$, and we make a
          dichotomy on the position of $v_a$. If there is $v_a'$ in  $
          [1,\phi(a)]$ within a distance at most $2\eta'$ of $v_a$, 
          we say we are in case  $(\beta_2)$). Otherwise, we are in case $(\beta_1)$.

          We now treat case $(\alpha)$. One has $d(h,\phi(a)h)\leq
          2\eta' + d(v,v_a) \leq 2\eta' +d(v',v_a)
          +\delta$.

          In the case $\alpha_1$, on the segment  
          $[\phi(a),\phi(ab)]$  we
          have  
          \[    d(v_a,v')=  d(\phi(a),\phi(ab)) -  d(\phi(a), v')) -
          d(\phi(ab), v_a ).        \]

           By the  requirement in Definition \ref{def:Lh:RH} that quasi-prefix and quasi-suffix are disjoint, 
		and since $d(1,h \geq \rho$, we have  that 
           $d(\phi(ab), v_a ) \geq \rho - 2\eta'$,
          in such a way that 
          $d(v_a,v')\leq  d(\phi(a),\phi(ab)) - \rho +2\eta' $, and  
          $d(h,\phi(a)h)\leq
          4\eta' 
          +\delta - \rho + d(\phi(a),\phi(ab))$. 
          Since $\rho$ was chosen larger than 
          $20(\epsilon+\eta+10\delta) > 4\eta'+\delta$, one has  
          $d(1,h^{-1}\phi (a) h)< d(1,\phi(b))$.

          Case $(\alpha_2)$ is similar,  
          $d(v_a,v')=  d(\phi(a),\phi(ab)) -  d(\phi(a), v_a))
          - d(\phi(ab), v' )$, and we use  
          $d(\phi(a),v_a) \geq d(\phi(a),h)-\eta'\geq \rho-\eta'$ to get 
           $d(h,\phi(a)h)
          \leq 3\eta' +\delta  -\rho
          +d(\phi(a),\phi(ab))$. Once again 
          $d(1, h^{-1} \phi (a) h)< d(1,\phi(b))$.

          In  case $(\beta_1)$,  since
          $v_a$ is at distance greater than $2\eta'$ 
          from $[1,\phi(a)]$, $\phi(a)h$
          is at distance greater than $\eta'$ from     $[1,\phi(a)]$,      
           and by Lemma \ref{lem;pour_les_geod1_rh}, no representative in 
          $\mathcal{L}$ of $\phi(a^{-1})$  has 
          $h$ as quasi prefix.
	Therefore $\phi(a^{-1})$ in case $\alpha_1$ or
          $\alpha_2$. Hence we also get that  
          $d(1, h^{-1} \phi (a) h)< d(1,\phi(b))$.

          Let us now consider case $(\beta_2)$: both $v'$ and $v_a'$ are 
          on $[1,\phi(a)]$. We have $d(h,\phi(a)h)\leq
          2\eta' + d(v,v_a) \leq 4\eta' +d(v',v_a')
          +\delta$.

          If the vertices are in order $(1,v',v_a',\phi(a))$ on the segment
          $[1,\phi(a)]$, 
          then 
          \[d(1,\phi(a)) \geq d(1,v')+ d(v', v_a') \geq d(1,h)- \eta' 
          - \delta + d(v',v_a').\]     
          Hence, we get $d(h,\phi(a)h)\leq 5\eta' +2\delta  -\rho
          +d(1,\phi(a))$. Since $\rho \geq 20(\epsilon+\eta
          +10\delta)>5\eta'+2\delta$, 
          we find 
          $d(1,h^{-1} \phi (a) h)< d(1,\phi(a))$.

          Finally, if the vertices are in 
          order $(1,v_k',v',\phi(a))$ on the segment
          $[1,\phi(a)]$, one can  bound $d(v',v_a') \leq 2\eta' 
          +100\delta$.               
          This yields  $d(h,\phi(a)h)\leq  6\eta' 
          +101\delta$. Now $d(1,\phi(b)) \geq \rho -\eta' $. 
          Since
          $\rho > 20(\epsilon+\eta +10\delta)>7\eta'+101\delta$, 
          one gets the result  
          $d(1,h^{-1} \phi (a) h)< d(1,\phi(a))$. 

          This finishes the proof in all the cases.
\qed \end{proof}

        \begin{cor}\label{coro;exists_rh}
          For all conjugacy class of compatible monomorphism $H_1 \to H_2$,
          any acceptable lift of a homomorphism  minimizing  the 
          quantity $Q$ in its class,  satisfies $\Omega$. 
        \end{cor}

    This, with Corollary \ref{coro;fini_rh}, ends the proof of Proposition
    \ref{prop;Omega}.

\section{JSJ decompositions} \label{JSJ-RH}

In this subsection we investigate splittings of toral
relatively hyperbolic groups.  
The main purpose of this subsection is to prove Theorem 
\ref{JSJTheorem}, which describes a splitting of 
toral relatively hyperbolic groups which is canonical
enough for our needs.  For more information on JSJ decompositions,
see \cite{RipsSelaJSJ}, \cite{DunwoodySageev} or \cite{FujPap}.

The following three operations may be performed on a graph of groups to recover a graph of groups with an isomorphic fundamental group:

{\em 1: Conjugation:} Conjugate the entire graph of groups by some fixed $w$.

{\em 2: Sliding:}  The operation of sliding is a move corresponding to the relation
\[      (A \ast_{C_1} B) \ast_{C_2} D \cong (A \ast_{C_1} D) \ast_{C_2} B       ,       \]
if $C_1 \subseteq C_2$ (as subgroups of $A$).

{\em 3: Modifying boundary monomorphisms by conjugation:}  Suppose $H = A \ast_C = \langle A,t \mid t \alpha(c) t^{-1} = \omega(c), c \in C \rangle$.  For an element $a \in A$, let $\tau_a \in \text{Aut}(A)$ be conjugation by $a$ and replace $\alpha : C \to A$ by $\tau_a \circ \alpha$, and replace $t$ by $ta^{-1}$.  For an amalgamated free product, $H = A \ast_C B$, and $a \in A$, replace $\beta : C \to B$ by $\tau_a \circ \beta$ and replace $B$ by $aBa^{-1}$.

These three moves define an equivalence relation between graphs of groups, and any two primary JSJ decompositions will be equivalent under this relation (see Theorem \ref{JSJTheorem} below).

\begin{lem} \cite[Lemma 2.4]{Groves_RH1} \label{malnormal}
Let $\Gamma$ be a toral relatively hyperbolic group.  Then all maximal
abelian subgroups of $\Gamma$ are malnormal.  Each abelian subgroup of $\Gamma$ is contained in a unique maximal abelian subgroup.
\end{lem}

Given Lemma \ref{malnormal}, the proof of the following lemma is identical to that of \cite[Lemma 2.1]{SelaDio1} (and to the omitted proof of Lemma \ref{MakeAbEll} above). 

\begin{lem} \label{AbelianElliptic}
Let $\Gamma$ be a toral relatively hyperbolic group, let $M$ be a noncyclic maximal abelian subgroup of $\Gamma$ and let $A$ be an abelian subgroup of $\Gamma$.  Then
\begin{enumerate}
\item If $\Gamma = U \ast_A V$ then $M$ can be conjugated into either $U$ or $V$; and
\item If $\Gamma = U \ast_A$ then either $M$ can be conjugated into $U$, or $M$ can be conjugated to $M'$ and $\Gamma = U \ast_A M'$.
\end{enumerate}
\end{lem}

Using this lemma, if we have a splitting $\Gamma = U \ast_A$ where $A$ is a subgroup of a non-elliptic maximal abelian subgroup $M$, we can convert it into the amalgamated free product $\Gamma = U \ast_A M'$.  Thus we consider only splittings in which every noncyclic abelian subgroup of $\Gamma$ is elliptic. 

 It is also
worth remarking that we know (by \cite[Lemma 2.2]{Groves_MR}, for example) that all abelian subgroups of
toral relatively hyperbolic groups are finitely generated.  Since maximal abelian subgroups of
a toral relatively hyperbolic group are malnormal, and since we are only considering
primary splittings, there will not be two abelian vertex groups in any of our splittings
which are adjacent.

Recall the definition of an unfolded splitting from Definition \ref{d:unfolded}.  We need
an analog of \cite[Theorem 6.1]{RipsSelaJSJ} to see that there is not an infinite sequence of unfoldings
of abelian splittings of a in which all noncyclic abelian edge groups are elliptic.  In order
to make this proof work {\em verbatim}, we need the splittings to be acylindrical.

\begin{defn} \cite{SelaAcyl}
A splitting of a group $H$ is called {\em $k$-acylindrical} if for every element $h \in H \smallsetminus
\{ 1 \}$, the fixed set of $h$ in the Bass-Serre tree has diameter at most $k$.
\end{defn}

The proof of the following lemma from \cite{SelaDio1} also works in our setting.

\begin{lem} [cf. Lemma 2.3, \cite{SelaDio1}] \label{l:2acyl}
Let $\Gamma$ be a toral relatively hyperbolic group.  A splitting of $\Gamma$ in which
all edge groups are abelian and all noncyclic abelian groups are elliptic can be modified
(by modifying boundary monomorphisms by conjugations and sliding operations) to be
$2$-acylindrical.
\end{lem}

Given Lemma \ref{l:2acyl}, the following theorem is proved in exactly the same
way as \cite[Theorem 6.1]{RipsSelaJSJ}.

\begin{thm} [cf. Theorem 6.1, \cite{RipsSelaJSJ}] \label{t:unfolded}
Let $\Gamma$ be a toral relatively hyperbolic group.  There is not an infinite sequence
of unfoldings of abelian splittings of $\Gamma$ in which all noncyclic abelian subgroups
are elliptic.
\end{thm}

Theorem \ref{t:unfolded} shows the existence of unfolded splittings.

We now describe the construction of the primary JSJ decomposition of a toral
relatively hyperbolic group, which encodes all of the primary splittings.  This
construction is entirely analogous to the construction for limit groups in
\cite[$\S$2]{SelaDio1}.

Consider two elementary abelian splittings of a group $\Gamma$, where
noncyclic abelian subgroups are elliptic in both splittings.  Let the edge
groups be $A_1$ and $A_2$ and the associated Bass-Serre trees $T_1$ and
$T_2$.  The first splitting is called {\em elliptic} in $T_2$ if $A_1$ fixes a
point in $T_2$, and is {\em hyperbolic} otherwise. The two splittings are
called {\em elliptic--elliptic} if $A_1$ is elliptic in $T_2$ and $A_2$ is
elliptic in $T_1$, {\em elliptic--hyperbolic} if $A_1$ is elliptic in $T_2$ and 
$A_2$ is hyperbolic in $T_1$, and {\em hyperbolic--hyperbolic} if $A_1$ is
hyperbolic in $T_2$ and $A_2$ is hyperbolic in $T_1$.

The proof of \cite[Theorem 2.1]{RipsSelaJSJ} applies in this setting to give
the following result.

\begin{thm} \label{NoEllHyp}
Suppose that $\Gamma$ is a freely indecomposable toral relatively
hyperbolic group.  Then any pair of elementary abelian splittings of $\Gamma$
in which all noncyclic abelian subgroups are elliptic is either
elliptic--elliptic or hyperbolic--hyperbolic.
\end{thm}

Given Theorem \ref{NoEllHyp}, we can restrict our attention to 
hyperbolic--hyperbolic and elliptic--elliptic splittings. 
Since noncyclic abelian subgroups are elliptic in all the splittings
which we consider, a pair of hyperbolic--hyperbolic splittings must
both have cyclic edge groups.  Such a pair of splittings gives rise to 
a quadratically hanging subgroup of $\Gamma$.  

We now briefly recall the notions of {\em quadratically hanging},
{\em EMQH} 
and {\em CEMQ} subgroups from \cite[$\S$4]{RipsSelaJSJ}.
Since all of our groups are torsion-free, we restrict attention to
surface groups, rather than the more general $2$-orbifold groups.

\begin{defn} \label{QH}
\cite[Definition 4.3, p.69]{RipsSelaJSJ}
A subgroup $Q$ of $G$ is called {\em quadratically hanging (QH)}
if $G$ admits a cyclic splitting in which $Q$ is a vertex group, all
edge groups are cyclic puncture subgroups of $Q$, and $Q$ admits
one of the following presentations:
\begin{eqnarray*}
&\langle p_1,\ldots , p_m, a_1, \ldots , a_g, b_1, \ldots , b_g \mid
\prod_{k=1}^m p_k \prod_{j=1}^g [a_j,b_j] \rangle ,\\
&\langle p_1,\ldots , p_m, v_1, \ldots , v_g \mid \prod_{k=1}^m p_k
\prod_{j=1}^g v_j^2 \rangle.
\end{eqnarray*}

A \emph{socket subgroup} is a subgroup obtained from a QH  group by adding all maximal roots of  the punctures. 
They admit standard presentations: 

  \begin{eqnarray*}
&\langle q_1,\ldots , q_m, a_1, \ldots , a_g, b_1, \ldots , b_g \mid
\prod_{k=1}^m q_k^{n_k} \prod_{j=1}^g [a_j,b_j] \rangle ,\\
&\langle q_1,\ldots , q_m, v_1, \ldots , v_g \mid \prod_{k=1}^m q_k^{n_k}
\prod_{j=1}^g v_j^2 \rangle.
\end{eqnarray*}

\end{defn}
Denote the splitting of $G$ corresponding to $Q$ by $\Lambda_Q$,
and the surface with fundamental group $Q$ by $S_Q$.

\begin{defn} \cite[$\S$4]{RipsSelaJSJ}
A s.c.c. $l \subset S_Q$ which is not boundary parallel,
not the core of a M\"obius band in $S_Q$ and not null-homotopic
is called a {\em weakly essential s.c.c.}  If furthermore $l$ is such
that no connected component of $S_Q \smallsetminus l$ is a M\"obius
band then we call $l$ an {\em essential s.c.c.}.
\end{defn}

We always assume that the surface $S_Q$ associated to a quadratically
hanging subgroup has negative Euler characteristic and contains a pair
of intersecting weakly essential s.c.c.

\begin{defn} \label{EMQ}
\cite[Definition 4.4, p.70]{RipsSelaJSJ}
An {\em essential maximal quadratically hanging (EMQH) subgroup}
is a quadratically hanging subgroup $Q < G$ so that for every
essential cyclic splitting of $G = A \ast_C B$ (or $G = A \ast_C$), either
(i) $Q$ is conjugate to a subgroup of $A$ or $B$; or (ii) the edge group $C$ can be
conjugated in $Q$.  In case (ii) the given splitting of $G$ is
assumed to be inherited
from the one corresponding to $Q$ obtained by splitting over
an element corresponding to an essential s.c.c. on $S_Q$.
\end{defn}

We now recall the {\em essential quadratic decomposition} of 
a f.g. group $G$.  (The torsion-free assumption in the next theorems
is merely so that we can restrict our attention to surface groups, and
is not important.)

\begin{thm} \label{CEMQ}
\cite[Theorem 4.18, pp.83-84]{RipsSelaJSJ}
Let $G$ be a one-ended torsion-free finitely generated group which
is not a surface group, and let $G = A_1\ast_{C_1}B_1$ (or $G = 
A_1\ast_{C_1}$), $\ldots$, $G = A_p\ast_{C_p} B_p$ (or $G = 
A_p \ast_{C_p}$) be essential cyclic splittings of $G$.  Suppose that
there exists an integer $q$ and a function $f : \{ 1, \ldots , q\} \to
\{ 1, \ldots , p \}$ which is surjective and so that for $1 \le i \le q-1$
the group $C_{f(i)}$ is hyperbolic in the cyclic splitting along 
$C_{f(i+1)}$, and so that $C_{f(q)}$ is hyperbolic in the splitting
along $C_{f(1)}$.  Then there exists an EMQH subgroup $Q < G$ with
the following properties:
\begin{enumerate}
\item[(i)] All the splittings along $C_1, \ldots , C_p$ are obtained from
$\Lambda_Q$ by cutting the surface $S_Q$ along essential s.c.c. 
corresponding to conjugates of the cyclic subgroups $C_1, \ldots , C_p$.
\item[(ii)] If $Q_1$ is an EMQH subgroup of $G$ and each of the cyclic 
subgroups $C_1, \ldots , C_p$ can be conjugated into $Q_1$ then
$Q$ can be conjugated into $Q_1$ and the surface $S_Q$ is a 
subsurface of $S_{Q_1}$.
\end{enumerate}
\end{thm}

\begin{prop}
\cite[Proposition 4.20, p.86]{RipsSelaJSJ}
With the notation and assumptions above, let $Q_1$ and $Q_2$ be
EMQH subgroups of $G$ constructed according to Theorem 
\ref{CEMQ}.  Then either $Q_1$ is conjugate to $Q_2$ or
$Q_1$ is elliptic in $\Lambda_{Q_2}$ and $Q_2$ is elliptic
in $\Lambda_{Q_1}$.
\end{prop}

The EMQH subgroups constructed in Theorem \ref{CEMQ} will
be called {\em canonical essential maximal QH (CEMQ) subgroups}
of $G$.

\begin{thm} \label{EQD}
\cite[Theorem 4.21, p.87]{RipsSelaJSJ}
Let $G$ be a torsion-free finitely generated group with one end 
which is not a surface group.  There exists a (canonical) reduced
cyclic splitting of $G$, called the {\em essential quadratic 
decomposition} of $G$, with the following properties:
\begin{enumerate}
\item[(i)] Every CEMQ subgroup of $G$ is conjugate to a vertex group
in the essential quadratic decomposition.  In particular, there are only
finitely many conjugacy classes of CEMQ subgroups.  Every edge group
is a cyclic boundary subgroup of one of the CEMQ subgroups and 
every vertex with a non-CEMQ vertex group is adjacent only to
vertices stabilised by CEMQ subgroups.
\item[(ii)] An essential cyclic splitting $G = A \ast_C B$ or $G = A\ast_C$
which is hyperbolic in another essential elementary cyclic splitting is
obtained from the essential quadratic decomposition of $G$ by cutting
a surface corresponding to a CEMQ of $G$ along an essential s.c.c.
\item[(iii)] The edge group of any essential cyclic splitting 
$G = A\ast_C B$ or $G = A \ast_C$ can be conjugated into a vertex
group of the essential quadratic decomposition.  In case it can be 
conjugated into a vertex group which is not a CEMQ subgroup, the
given elementary cyclic splitting is elliptic-elliptic with respect to any
other elementary essential cyclic splitting of $G$.
\item[(iv)] The essential quadratic decomposition of $G$ is unique up
to sliding, conjugation and modifying boundary morphisms by
conjugation.
\end{enumerate}
\end{thm}

In general, the cyclic subgroups associated to the boundary components
of a CEMQ subgroup may not be maximal.  Therefore, we need the
following notion, which is due to Sela (see \cite[p.569]{SelaGAFA}):

\begin{defn}
Suppose that $\Gamma$ is a toral
relatively hyperbolic group.   A {\em canonical socket
subgroup} of $\Gamma$ is the subgroup obtained by adding
all maximal roots to the punctures of a CEMQ subgroup
of $\Gamma$.
\end{defn}

Given the essential quadratic decomposition, we continue by 
constructing the 
primary cyclic JSJ decomposition of $\Gamma$ in two steps.
First, turn each CEMQ subgroup into a canonical socket group,
by adding all roots of the puncture elements.  This makes
the splitting essential, and primary.  Second, refine this splitting
 by considering 
essential cyclic splittings
whose edge groups are elliptic in the quadratic decomposition, and 
so that every noncyclic abelian subgroup is elliptic.  This refinement
procedure terminates by generalised accessibility.

From the primary cyclic JSJ decomposition, we construct the primary 
JSJ  decomposition of $\Gamma$ by refining the non-QH and 
non-abelian vertex
groups in the primary cyclic JSJ decomposition, by considering 
non-cyclic
abelian decompositions of $\Gamma$ in which all noncyclic abelian
subgroups are elliptic.  Once again, Bestvina and Feighn's generalised
accessibility guarantees that this procedure terminates.

In general, the JSJ decomposition of a finitely presented group need
not be unique up to conjugacy, modifying boundary morphisms by
conjugation and sliding (see \cite{Forester}).  However, the
JSJ decomposition of a hyperbolic group is unique up to such
moves, see \cite[Theorem 1.7]{SelaGAFA}.  The proof from
\cite{SelaGAFA} applies directly in this case.  Thus the primary JSJ 
decomposition of $\Gamma$ is unique up to conjugation, modifying
boundary morphisms by conjugation and sliding.

In summary, we have the following

\medskip

{\noindent \bf Theorem  \ref{JSJTheorem} }
{\em
[cf. Theorem 2.7, \cite{SelaDio1}; see also Theorem 7.1, \cite{RipsSelaJSJ}]
Suppose $\Gamma$ is a freely indecomposable toral relatively hyperbolic group.  There exists
a reduced unfolded splitting of $\Gamma$ with abelian edge groups, which we call a {\em primary JSJ decomposition} of $\Gamma$, satisfying the following:
\begin{enumerate}
\item Every canonical socket subgroup of $\Gamma$ is conjugate to a vertex
group in 
the JSJ decomposition.  Every QH subgroup of $\Gamma$ can be conjugated into one of the CEMQ
subgroups of $\Lambda$.  Every vertex group in the JSJ decomposition which is not a socket subgroup
of $\Gamma$ is elliptic in any primary splitting of $\Gamma$;
\item A one edge primary splitting $\Gamma = D \ast_A E$ or $\Gamma = D \ast_A$ which is 
hyperbolic in another primary splitting is obtained from the primary JSJ decomposition of $\Gamma$ by cutting a 
surface  corresponding to a CMQ subgroup of $\Gamma$ along an essential s.c.c;
\item Let $\Theta$ be a one edge primary splitting $\Gamma = D \ast_A E$ or $\Gamma = D \ast_A$,
which is elliptic with respect to any other one edge primary splitting of $\Gamma$.  Then $\Theta$ is
obtained from the JSJ decomposition of $\Gamma$ by a sequence of collapsings, foldings and
conjugations;
\item If JSJ$_1$ is another JSJ decomposition of 
$\Gamma$ then JSJ$_1$ is obtained from the
JSJ decomposition by a sequence of slidings, conjugations and modifying boundary
monomorphisms by conjugations.
\end{enumerate}
}


\begin{thebibliography}{99} 

\bibitem{Adian} 
S. I.   \smalltextsc{Adian}, 
The unsolvability of certain algorithmic problems in the theory of groups, 
\textit{Trudy Moskov. Obsc.} {\bf 6} (1957), 231--298.

\bibitem{Alibegovic} 
E.  \smalltextsc{Alibegovi\'c}, 
A combination theorem for relatively hyperbolic groups, 
\textit{Bull. London Math. Soc.},  {\bf 37}  (2005),  no. 3, 459--466.

\bibitem{BGS} 
G.   \smalltextsc{Baumslag}, D. \smalltextsc{Gildenhuys} and R. \smalltextsc{Strebel}, 
Algorithmically insoluble problems about finitely presented soluble groups, Lie and associative algebras, I, 
\textit{J. Pure and Appl. Algebra} {\bf 39} (1986), 53--94.

\bibitem{Bel} 
I. \smalltextsc{Belegradek}, 
Aspherical manifolds with relatively hyperbolic fundamental groups, 
 {\it Geom. Dedicata } {\bf 129}  (2007), 119--144. 

\bibitem{Bestvina} 
M. \smalltextsc{Bestvina}, 
Degenerations of hyperbolic space, 
\textit{Duke Math. J.} {\bf 56} (1988), 143--161.


\bibitem{BFAccess} 
M. \smalltextsc{Bestvina} and M. \smalltextsc{Feighn}, 
Bounding the complexity of  simplicial group actions, 
\textit{Invent. Math.} {\bf 103} (1991), 449--469.
  
\bibitem{BF} 
M. \smalltextsc{Bestvina} and M. \smalltextsc{Feighn}, 
Stable actions of groups on real trees, 
\textit{Invent. Math.} {\bf 121} (1995), 287--321.



\bibitem{B_JSJ} 
B. H. \smalltextsc{Bowditch}, 
Cut points and canonical splittings of hyperbolic groups, 
\textit{Acta Math.} {\bf 180} (1998) 145--186.

\bibitem{B_rh}  
B. H. \smalltextsc{Bowditch },
Relatively Hyperbolic Groups, 
preprint (1999).

\bibitem{Bowditch_periph} 
B. H. \smalltextsc{Bowditch}, 
Peripheral splittings of  groups 
{\it Trans. Amer. Math. Soc.}  {\bf 353} (2001) 4057-4082. 
  
  \bibitem{BS} 
M. R. \smalltextsc{Bridson} and G. A. \smalltextsc{Swarup}, 
On Hausdorff-Gromov convergence and a theorem of Paulin, 
\textit{Enseign. Math.} {\bf 40} (1994), 267--289.
  
  \bibitem{Brown} 
K. S. \smalltextsc{Brown}, 
Cohomology of groups, 
Graduate  Texts in Mathematics, {\bf 87}, Springer, 1982.



\bibitem{Bumagin} 
I. \smalltextsc{Bumagin}, 
The Conjugacy Problem for Relatively Hyperbolic Groups, 
 \textit{Algebr. Geom. Topol.}  {\bf 4}  (2004), 1013--1040.
 

\bibitem{BKM} 
I. \smalltextsc{Bumagin}, O. \smalltextsc{Kharlampovich} and A. \smalltextsc{Miasnikov}, 
Isomorphism problem for finitely generated fully residually free groups,    
J. Pure Appl. Algebra  208  (2007),  no. 3, 961--977. 


\bibitem{CDP} 
M. \smalltextsc{Coornaert}, T. \smalltextsc{Delzant}, and A. \smalltextsc{Papadopoulos}, 
G\'eom\'etrie et th\'eorie des groupes, Les groupes hyperboliques de M. Gromov, 
Lecture Notes in Math. 1441, Springer, 1991. 


\bibitem{Dah_thesis} 
F. \smalltextsc{Dahmani}, 
Les groupes relativement hyperboliqes et leurs bords, 
 PhD Thesis, Strasbourg  (2003) 

\bibitem{DahComb} 
F. \smalltextsc{Dahmani}, 
Combination of convergence groups, 
\textit{Geom. Top.} {\bf 7} (2003), 933--963.


\bibitem{Dah_find} 
F. \smalltextsc{Dahmani}, 
Finding relative hyperbolic structures, 
\textit{Bull. London Math. Soc.} {\bf 40} 3, (2008) 395-404.


\bibitem{Dah_eq} 
\smalltextsc{F. Dahmani}, 
Existential questions in (relatively) hyperbolic groups, 
\textit{ Israel J. Math.} to appear.



\bibitem{DahG_freeprod} 
F. \smalltextsc{Dahmani} and D. \smalltextsc{Groves}, 
Detecting free splittings in relatively hyperbolic groups, 
\textit{Trans. Amer. Math. Soc.} to appear.

\bibitem{Dehn1} 
M. \smalltextsc{Dehn}, 
\"Uber unendliche diskontinuierliche Gruppen, 
\textit{Math. Ann.} {\bf 71} (1912), 413--421.

\bibitem{Dehn2} 
M. \smalltextsc{Dehn}, 
Papers on group theory and topology, Translated from the German and with introductions and an appendix by John Stillwell. With an appendix by Otto Schreier. Springer-Verlag, New York, 1987. viii+396 pp.  


\bibitem{DGH} 
V. \smalltextsc{Diekert}, C. \smalltextsc{Guti\'errez} and C. \smalltextsc{Hagenah},  
The existential theory of equations with rational constraints in free groups is PSPACE-complete.  
STACS 2001 (Dresden),  170--182, 
Lecture Notes in Comput. Sci., 2010, 
Springer, 2001.

\bibitem{Diekert_Muscholl} 
V. \smalltextsc{Diekert} and  A. \smalltextsc{Muscholl},  
Solvability of equations in free partially commutative groups is decidable.  
Automata, languages and programming, 543--554,
Lecture Notes in Comput. Sci., 2076,
Springer,  2001. 

\bibitem{DS} 
C. \smalltextsc{Dru\c{t}u} and M. \smalltextsc{Sapir}, 
Tree-graded spaces and asymptotic cones of groups,  
{\it Topology} {\bf 44} (2005), 959-1058. 

\bibitem{DS2} 
C. \smalltextsc{Dru\c{t}u} and M. \smalltextsc{Sapir}, 
Groups acting on tree-graded spaces and splittings of relatively hyperbolic groups, 
{\it Adv. Math.} {\bf 217} (2007), 1313-1367.


\bibitem{DunwoodySageev} 
M. \smalltextsc{Dunwoody} and M. \smalltextsc{Sageev}, 
JSJ-splittings for finitely presented groups over slender groups, 
\textit{Invent. Math.} {\bf 135} (1999), 25--44.

\bibitem{E+} 
D. \smalltextsc{Epstein}, J. \smalltextsc{Cannon}, D. \smalltextsc{Holt}, S. \smalltextsc{Levy}, M. \smalltextsc{Paterson} and W. \smalltextsc{Thurston}, 
Word processing in groups, 
Jones and Bartlett, Boston, 1992.

\bibitem{Farb} 
B. \smalltextsc{Farb}, 
Relatively hyperbolic Groups, 
\textit{Geom. Funct. Anal.} {\bf 8} (1998), no. 5, 810--840. 


\bibitem{Forester} 
M. \smalltextsc{Forester}, 
On uniqueness of JSJ decompositions of finitely generated groups, 
\textit{Comment. Math. Helv.} {\bf 78} (2003), 740--751.

\bibitem{FujPap} 
K. \smalltextsc{Fujiwara} and P. \smalltextsc{Papasoglu}, 
JSJ decompositions of finitely presented groups and complexes of groups,   
{\it Geom. and Funct. Anal.} {\bf 16}, no 1 (2006), 70-125.


\bibitem{Gerasimov} 
V. \smalltextsc{Gerasimov}, 
Detecting connectedness of the boundary of a hyperbolic group, 
preprint.

\bibitem{GL} 
R. \smalltextsc{Grigorchuk} and I. \smalltextsc{Lysenok},  
A description of solutions of quadratic equations in hyperbolic groups. 
\textit{Internat. J. Algebra Comput.} {\bf 2} (1992), no. 3, 237--274.

\bibitem{Grom} 
M. \smalltextsc{Gromov}, 
Hyperbolic groups.  
Essays in group theory,  
75--263, Math. Sci. Res. Inst. Publ., 8, Springer, New York, 1987. 


\bibitem{GrovesCAT(0)1} 
D. \smalltextsc{Groves}, 
Limits of certain CAT$(0)$ groups, I: Compactification,
\textit{Alg. and Geom. Top.} {\bf 5} (2005), 1325--1364.

\bibitem{Groves_RH1} 
D. \smalltextsc{Groves}, 
Limit groups for relatively hyperbolic groups, I: The basic tools, preprint. 

\bibitem{Groves_MR} 
D. \smalltextsc{Groves}, 
Limit groups for relatively hyperbolic groups, II: Makanin-Razborov diagrams,    
{\it Geom. Top. } {\bf 9} (2005),  2319--2358.



\bibitem{GrunSeg} 
F. \smalltextsc{Grunewald} and D. \smalltextsc{Segal}, 
Some general algorithms. II. Nilpotent groups, 
\textit{Ann. Math. (2)} {\bf 112} (1980), 585--617.

\bibitem{Guirardel} 
V. \smalltextsc{Guirardel}, 
Limit groups and groups acting freely on $\R^n$-trees,    
{\it Geom. Top. } {\bf 8} (2004),  1427--1470.

\bibitem{HK} 
G. C. \smalltextsc{Hruska} and B. \smalltextsc{Kleiner}, 
Hadamard spaces with  isolated flats. With an appendix by the authors and Mohamad Hindawi.
\textit{Geom. Top.} {\bf 9} (2005), 1501--1538.


\bibitem{Hummel} 
C. \smalltextsc{Hummel}, 
Rank one lattices whose parabolic isometries have no rational part, 
{\it Proc. Amer. Math. Soc.} {\bf 126}, 8, (1998), 2453--2458.

\bibitem{KM} 
O. \smalltextsc{Kharlampovich} and A. \smalltextsc{Miasnikov}, 
Effective JSJ decompositions, in \textit{``Groups, languages and algorithms"},
Contemp. Math. {\bf 378} (2005), 87--212.

\bibitem{KM2} 
O. \smalltextsc{Kharlampovich} and A. \smalltextsc{Miasnikov}, 
Elementary theory of free non-abelian  groups, 
 {\it J. Algebra } {\bf 302}  (2006),  no. 2, 451--552.

\bibitem{Levitt} 
G. \smalltextsc{Levitt}, 
Automorphisms of hyperbolic groups and graphs of groups,  
{\it Geom. Dedicata}  {\bf 114}  (2005), 49--70.

\bibitem{Levitt_rigid} G. \smalltextsc{Levitt}, 
Characterizing rigid simplicial actions on trees. {\it in  Geometric methods in group theory},  
27--33, Contemp. Math., {\bf 372} (2005).


\bibitem{LS} 
R. C. \smalltextsc{Lyndon} and P. E. \smalltextsc{Schupp}, 
Combinatorial group theory, 
\textit{Ergebnisse der Mathematik und ihrer Grenzgebiete} {\bf 89}, Springer-Verlag, Berlin, 1977.

\bibitem{Ly} 
I. \smalltextsc{Lysenok},  
Some algorithmic properties of hyperbolic groups.  
\textit{Izv. Akad. Nauk SSSR Ser. Mat.} {\bf 53} (1989), no. 4, 814--832, 912; translation in \textit{Math. USSR-Izv. } {\bf 35} (1990), no. 1, 145--163.

\bibitem{Makanin} 
G. S. \smalltextsc{Makanin}, 
Decidability of the universal and positive theories
of a free group. (Russian) \textit{Izv. Akad. Nauk SSSR Ser. Mat.} {\bf 48} (1984), 735--749; translation in \textit{Math. USSR-Izv.} {\bf 25} (1985), 75--88.

\bibitem{Miller} 
C. F. \smalltextsc{Miller III}, 
On group theoretic decision problems and their classification, 
\textit{Ann. Math Studies} {\bf 68}, Princeton University Press, Princeton, 1971.

\bibitem{Miller:dec} 
C. F. \smalltextsc{Miller III}, Decision problems for groups --
survey and reflections, in \textit{``Algorithms and classification in
combinatorial group theory (Berkeley, CA, 1989)"}, 
MSRI Publ. {\bf 23}
(1992), 1--59.

\bibitem{Osin} D. V. \smalltextsc{Osin}, 
Relatively hyperbolic groups: Intrinsic geometry, algebraic properties, and
algorithmic problems, 
\textit{Mem. Amer. Math. Soc.} {\bf 179} (2006).

\bibitem{Papa} 
P. \smalltextsc{Papasoglu}, 
An algorithm detecting hyperbolicity, 
Geometric and computational perspectives on infinite groups (Minneapolis, MN and New Brunswick, NJ, 1994), 193--200,
DIMACS Ser. Discrete Math. Theoret. Comput. Sci., 25,
Amer. Math. Soc., Providence, RI, 1996. 

\bibitem{Paulin} 
F. \smalltextsc{Paulin}, 
Topologie de Gromov \'equivariante, structures hyperbolique et arbres r\'eels, 
\textit{Invent. Math.} {\bf 94} (1988), 53--80.

\bibitem{Rabin} 
M. O. \smalltextsc{Rabin}, Recursive unsolvability of group
theoretic problems, 
\textit{Ann. Math.} {\bf 67} (1958), 172--194.

\bibitem{Rebbechi} 
D. \smalltextsc{Rebbechi}, 
Algorithmic Properties of Relatively Hyperbolic Groups. PhD Thesis,  ArXiv math.GR/0302245. 

\bibitem{RS} 
E. \smalltextsc{Rips} and Z. \smalltextsc{Sela}, 
Canonical representatives and equations in hyperbolic groups.  \textit{Invent. Math.}  {\bf 120}  (1995),  no. 3, 489--512. 

\bibitem{RipsSelaJSJ} 
E. \smalltextsc{Rips} and Z. \smalltextsc{Sela}, 
Cyclic splittings of finitely
  presented groups and the canonical JSJ decomposition,
  \textit{Ann. Math. (2)} {\bf 146} (1997), 53--104.

\bibitem{Romankov} 
V. \smalltextsc{Roman'kov}, 
Universal theory of nilpotent groups, 
\textit{Mat. Zametki} {\bf 25} (1979), 487--495, 635.




\bibitem{Segal} 
D. \smalltextsc{Segal}, 
Decidable properties of polycyclic groups, 
\textit{Proc. London Math. Soc.  (3)} {\bf 61} (1990), 497--528.


\bibitem{SelaIso} 
Z. \smalltextsc{Sela}, 
The isomorphism problem for hyperbolic groups I, 
\textit{Ann. Math. (2)} {\bf 141} (1995), 217--283.

\bibitem{SelaAcyl} 
Z. \smalltextsc{Sela}, 
Acylindrical accessibility for groups, 
\textit{Invent. Math.}  {\bf129} (1997), 527--565.

\bibitem{SelaGAFA} 
Z. \smalltextsc{Sela}, Structure and rigidity in (Gromov)  hyperbolic groups and discrete groups in rank $1$ Lie groups II,
\textit{Geom. Funct. Anal.} {\bf 7} (1997), 561--593.

\bibitem{SelaDio1} 
Z. \smalltextsc{Sela}, 
Diophantine geometry over groups, I: Makanin-Razborov diagrams, 
\textit{Publ. Math. Inst. Hautes \'Etudes Sci.} {\bf 93} (2003), 31--105.

\bibitem{SelaDio27} 
Z. \smalltextsc{Sela}, 
Diophantine geometry over groups, II--VII.  
Appeared in \textit{Israel J. Math.} and \textit{GAFA}.

\bibitem{SelaDio8} 
Z. \smalltextsc{Sela}, 
Diophantine geometry over groups, VIII: 
Elementary  Theory of Hyperbolic Groups, preprint (2002). 

\bibitem{Serre} 
J.-P. \smalltextsc{Serre}, 
Trees, 
Springer-Verlag, Berlin, 1980.

\end{thebibliography}
\end{document}